\documentclass[11pt]{article}
\usepackage{amsmath, amssymb}
\usepackage{mathtools}
\usepackage{algorithm}
\usepackage{algpseudocode}
\usepackage{fmtcount}
\usepackage{graphicx}
\usepackage{enumitem} 
\usepackage{subcaption}
\usepackage{afterpage}
\usepackage[noadjust]{cite}
\usepackage[font=small,labelfont=bf]{caption}
\usepackage[labelfont=bf]{caption}
\usepackage{setspace}
\allowdisplaybreaks
\usepackage[title,titletoc]{appendix}
\usepackage{epsfig}
\usepackage{graphicx}
\usepackage{hhline}
\usepackage{latexsym}
\usepackage{amssymb}
\usepackage{amsmath,amscd}
\usepackage{multirow} 
\usepackage{array}
\usepackage{caption}
\usepackage{subcaption}
\usepackage{authblk}
\usepackage{color}
\usepackage{natbib}
\usepackage{rotating}
\usepackage{graphicx,lscape}
\usepackage{amsmath, amsthm, amssymb}
\usepackage{graphicx}
\usepackage{epsfig}
\usepackage{graphics,color, graphicx}
\usepackage{latexsym}
\usepackage{fullpage}
\usepackage{booktabs,fixltx2e}
\usepackage[flushleft]{threeparttable}
\usepackage{amssymb,amsmath,amscd}
\usepackage{fancyvrb}
\usepackage{pdfpages}
\usepackage{fmtcount}
\usepackage{url}
\usepackage[margin=1in]{geometry}

\theoremstyle{plain}
\newtheorem{theorem}{Theorem}
\newtheorem{proposition}[theorem]{Proposition}
\newtheorem{corollary}[theorem]{Corollary}
\newtheorem{lemma}[theorem]{Lemma}
\newtheorem{definition}{Definition}
\newtheorem{remark}{Remark}
\newtheorem{assumption}{Assumption}

\usepackage{setspace}

\title{\bf Quantiles on global non-positive curvature spaces
	\medskip
}

\author[1]{Ha-Young Shin}
\author[2]{Hee-Seok Oh}
\affil[1]{
	Department of Statistics and Actuarial Science, Soongsil University
}
\affil[2]{
	Department of Statistics, Seoul National University
}
{
    \makeatletter
    \renewcommand\AB@affilsepx{: \protect\Affilfont}
    \makeatother

    \affil[ ]{Email}

    \makeatletter
    \renewcommand\AB@affilsepx{, \protect\Affilfont}
    \makeatother

    \affil[1]{hayoung.shin@gmail.com}
    \affil[2]{heeseok@stats.snu.ac.kr}
}

\begin{document}
	\maketitle
	
	\begin{abstract}
		We propose a notion of geometric quantiles on Hadamard spaces, or global non-positive curvature spaces. After providing some definitions and basic properties, including scaled isometry equivariance and a necessary condition on the gradient of the quantile loss function, we investigate asymptotic properties, such as strong consistency and joint asymptotic normality. We provide a detailed description of how to compute quantiles using a gradient descent algorithm in hyperbolic space. We detail several potential uses for these quantiles---defining measures of distributional characteristics, detecting outliers via a transformation-retransformation procedure, and testing the equality of distributions---and demonstrate their application through simulations and real biological data.
		\noindent
		
		\vspace{\baselineskip}
		
		
	\end{abstract}
	
	\pagenumbering{arabic}
	
	\section{Introduction} \label{intro}	
\begin{figure}
	\centering
	{\includegraphics[width=0.4\textwidth]{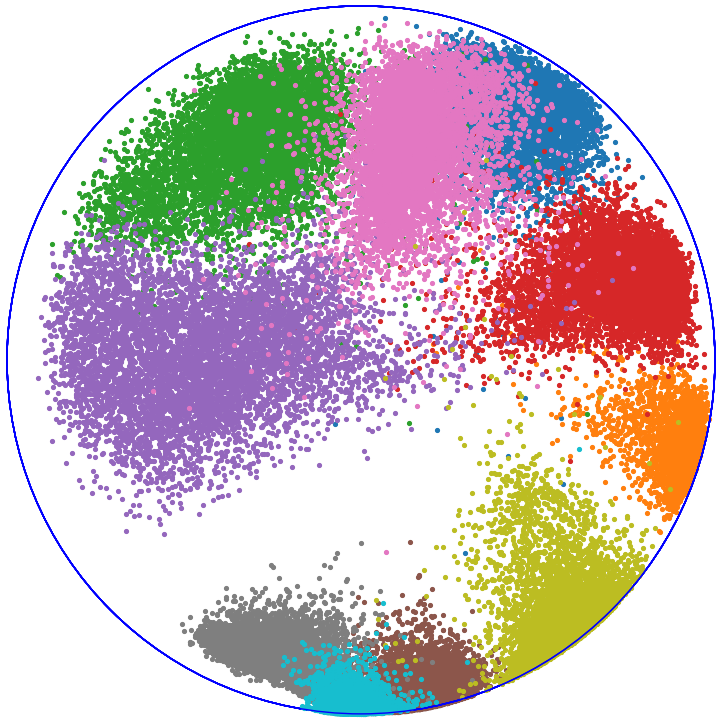}}{\includegraphics[width=0.5\textwidth]{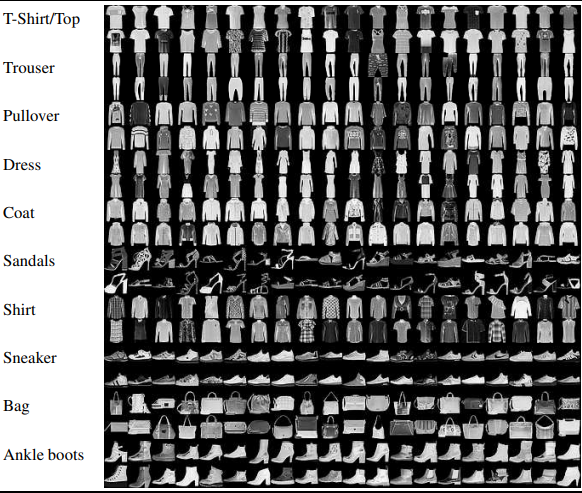}}
	 \vspace{-5mm}
    \caption{\small{On the left, the Fashion-MNIST data set embedded into the Poincar\'e disk, where distances increase rapidly to infinity near the disk's boundary. Colors represent different clothing classes. Key: blue=t-shirt, orange=trouser, green=pullover, red=dress, purple=coat, brown=sandals, pink=shirt, gray=sneaker, olive=bag, cyan=ankle boots. On the right, examples of images from each class in the data set [Reproduced from \citet{Xiao2017} under the terms of the Creative Commons Attribution 4.0 International License (CC BY 4.0)]. }}
    \label{fig:fashion}
\end{figure}

	%The median, quartiles, quintiles, deciles and percentiles are all examples of quantiles. Intuitively, the idea of a quantile is simple: Given data points $X_1,...,X_N$ and $\alpha\in(0,1)$, the $\alpha$-th quantile is the number which is greater than $100\alpha\%$ of the data and less than $100(1-\alpha)\%$ of the data. Given a random variable $X$, the formal definition of the $\alpha$-th quantile is
	%\begin{equation*}
	%q(\alpha;X)=\arg\min_{q\in\mathbb{R}}2E[|X-q|\{(1-\alpha)I(X\leq Q)+\alpha I(X>Q)\}],
	%\end{equation*}
	%where $I$ is the indicator function. That is, the quantile is the value that minimises a loss function that looks like a slanted check. Note that for $\alpha=0.5$, the $\alpha$-th quantile is the median, and the loss function becomes $E|X-q|$. Quantiles are useful because while the mean and median are useful ways of examining the centre of the data, quantiles allows one to see the shape of the rest of the data. Using a quantile loss function, one can perform quantile regression, just as least squares, or $L_2$, regression uses the squared loss function and median, or $L_1$, regression uses the absolute loss function. Quantile regression allows for an examination of the relationship between covariates and the response variable beyond the central tendency for which $L_2$ and $L_1$ regression are so useful.
		Many data sets are best understood as lying in a non-linear metric space. Among the most useful metric spaces are Hadamard spaces, also known as global non-positive curvature spaces. %f, which show up in fields as varied as phylogenetics, image processing, computational linguistics, and developmental biology
	Hilbert spaces such as $\mathbb{R}^n$ are the most basic examples. %Complete $\mathbb{R}$-trees are also Hadamard spaces, as are complete Euclidean buildings. 
	Phylogenetic tree spaces, the geometry of which was first studied by \citet{Billera2001}, are also Hadamard; phylogenetic trees represent evolutionary relationships, not only biological but also linguistic. Diffusion tensor imaging (DTI) data, used to probe the cellular structure of the brain, lie on the Hadamard manifold of symmetric positive definite (SPD) matrices. Hyperbolic spaces are the most commonly encountered non-Euclidean Hadamard spaces in data analysis as they have the property that volumes increase exponentially, as number of nodes in a tree, rather than polynomially as in Euclidean space, making it easier to embed data while preserving hierarchical relationships; %As such, data with an hierarchical structure have a natural home in hyperbolic space
see \citet{Klimovskaia2020} and \citet{Guo2022} for demonstrations of these advantages of hyperbolic embeddings over Euclidean methods like $t$-SNE. The ubiquity of such hierarchies is well-known in natural language processing, and \citet{Khrulkov2020} argued that natural hierarchies also exist in image data. Their methodology is used by \citet{Chien2021} to embed the Fashion-MNIST data set of \citet{Xiao2017}, which contains 70,000 images of clothing across 10 classes, into hyperbolic space; see Figure \ref{fig:fashion}. Many kinds of biological data can also be embedded in hyperbolic space; see Section \ref{experiments} for examples. Others have used hyperbolically embedded data to perform downstream tasks like deep learning and robust large-margin classification; see \citet{Mettes2024} and \citet{Weber2020}, respectively. See Section \ref{geomsupp} of the supplement for more information on hyperbolic spaces, including Poincar\'e disks. Besides hyperbolic spaces, \citet{Zhu2009} studied regression on SPD matrix spaces, and there has been research into statistics and probability on general Hadamard spaces; for examples, see \citet{Kostenberger2023}, \citet{Sturm2002}, and \citet{Yun2023}. %, and \citet{Zhang2016}. %\citet{Sturm2002} developed a theory of non-linear martingales and \citet{Sturm2003} studied probability theory on such spaces. \citet{Yun2023} generalised the notion of the median-of-means to these spaces and derived concentration inequalities for the median-of-means as an estimator of the population Fr\'echet mean, while \citet{Kostenberger2023} proved a strong law of large numbers for random elements with independent but not necessarily identically distributed values in Hadamard space under very weak conditions and applied this to the problem of robust signal recovery. \citet{Zhang2016} analysed first-order algorithms for convex optimisation on Hadamard manifolds.

	The mean easily generalises to $\mathbb{R}^n$ as the point that minimises the squared loss functions. As there is no obvious ordering of multivariate points, defining multivariate quantiles is less straightforward. \citet{Chaudhuri1996} provided an elegant solution. Let $\langle\cdot,\cdot\rangle_2$ and $\lVert \cdot\rVert_2$ denote the standard Euclidean inner product and norm, respectively. For a random $n$-vector $X$, \citet{Chaudhuri1996} defined the geometric $u$-quantile, where $u\in B^n(1):=\{v\in\mathbb{R}^n|\lVert v\rVert_2<1\}$, as the $p$ that minimises
	\begin{equation} \label{multiquantile}
	E[\lVert X-p\rVert_2+\langle u,X-p\rangle_2].
	\end{equation}
	Note that when $n=1$, $E[|X-p|+u(X-p)]=2E[|X-p|\{(1-\alpha)I(X\leq p)+\alpha I(X>p)\}]$ for $u=2\alpha-1$; that is, the familiar $\alpha$th univariate quantile, where $\alpha\in(0,1)$, is precisely the geometric quantile indexed by $2\alpha-1\in(-1,1)=B^1(1)$. %\citet{Chaudhuri1996} observed that sample geometric quantiles exist and are unique if the data are not collinear. They can be computed using first-order or Newton-Raphson-type methods and a Bahadur-type representation can be used to show consistency and asymptotic normality. %\citet{Girard2015} and \citet{Girard2017} proved several properties of extreme geometric quantiles, that is, quantiles corresponding to vectors with large norms close to 1.
	
	Although the generalisation from $\mathbb{R}^n$ to metric spaces is straightforward for the mean, the quantile remains a challenge. In (\ref{multiquantile}), both $X-p$ and $u$ can be considered vectors as $p$. This $u$ is in the tangent space of every point in $\mathbb{R}^n$, so we can index our quantiles with $u\in B^n(1)$. However, on a general Riemannian manifold, there is no canonical way to identify tangent vectors from different tangent spaces with each other. In this paper, we propose a notion of quantiles on Hadamard spaces by using the so-called boundary at infinity to provide a canonical sense of direction.

Quantiles provide a sense of a distribution's structure beyond what the mean or median can convey; similarly, quantile regression gives a richer description of the relationship between variables than ordinary least squares regression. This advantage of quantiles also allows one to measure the dissimilarity between distributions, test whether data sets come from different distributions, and, if so, find the loci of differences. One may also define a notion of tangent vector-valued ranks, closely related to quantiles, and use it for tests of location. Isoquantile contours can provide a sense of the shape of a data cloud and can therefore be used to define measures of spread, skewness, and kurtosis, as well as to detect outliers.

	Section \ref{geom} provides background on Hadamard spaces, and Section \ref{def} introduces our definition of quantiles on these spaces along with some basic properties. Section \ref{lst} investigates the asymptotic theory of sample quantiles on Hadamard manifolds, while Section \ref{hyp} explains how to calculate quantiles on hyperbolic spaces. Details on some of the applications mentioned above are provided in Section \ref{applications}. Experiments are conducted in Section \ref{experiments}. Section \ref{conc} concludes with remarks on future work.
	
	\section{Hadamard spaces and the boundary at infinity} \label{geom}

This section provides a brief overview of Hadamard spaces; for more detailed information, see Section \ref{geomsupp} of the supplement. Throughout this paper, $(M,d)$, usually shortened to $M$, will represent a Hadamard space, and we will denote all metrics by $d$. Intuitively, a Hadamard space can be thought of as a metric space in which triangles are at least as `thin' as Euclidean triangles. We give special attention to Hadamard manifolds; that is, Riemannian manifolds that are Hadamard spaces. All Riemannian manifolds in this paper are assumed to be $C^\infty$ with $C^\infty$ metrics (recall that a function is called $C^\infty$ if it is infinitely differentiable, and a manifold is called $C^\infty$ if its transition maps are $C^\infty$). Hadamard manifolds can also be defined as complete, simply-connected Riemannian manifolds of non-positive sectional curvature.
	
	In a metric space $M$, geodesic rays $c_1,c_2:[0,\infty)\rightarrow M$ are said to be asymptotic if there exists a constant $B$ such that $d(c_1(t),c_2(t))\leq B$ for all $t\geq0$. If $M$ is Hadamard, it can be shown that $d(c_1(t),c_2(t))\leq d(c_1(0),c_2(0))$ for all asymptotic $c_1, c_2$; see the proof of Proposition II.8.2 of \citet{Bridson1999}. In Euclidean space, unit-speed rays are asymptotic if and only if they are parallel, and in hyperbolic space, if their extensions to geodesic lines are limiting parallels. Letting two geodesic rays in $M$ be equivalent if they are asymptotic, we define the boundary at infinity of $M$, $\partial M$, as the set of equivalence classes of unit-speed geodesic rays. Crucially, if $M$ is a Hadamard space and $\xi\in\partial M$, there is exactly one geodesic ray in $\xi$ issuing from each point in $M$; thus, $\xi$ can define a direction on $M$. With an appropriately defined topology, the boundary of an $n$-dimensional Hadamard manifold $M$ is homeomorphic to $S^{n-1}:=\{p\in\mathbb{R}^n|\lVert p\rVert_2=1\}$ through the map taking $\xi$ to $\xi_x\in T_xM$, the unit vector in the tangent space $T_xM$ at $x$ that is the velocity of the unique geodesic ray in $\xi$ issuing from $x$, for any $x\in M$. The boundary of the Poincar\'e ball provides an intuitive visualisation of $\partial M$ for hyperbolic space.
	
\begin{figure}
	\centering
	{\includegraphics[width=0.32\textwidth]{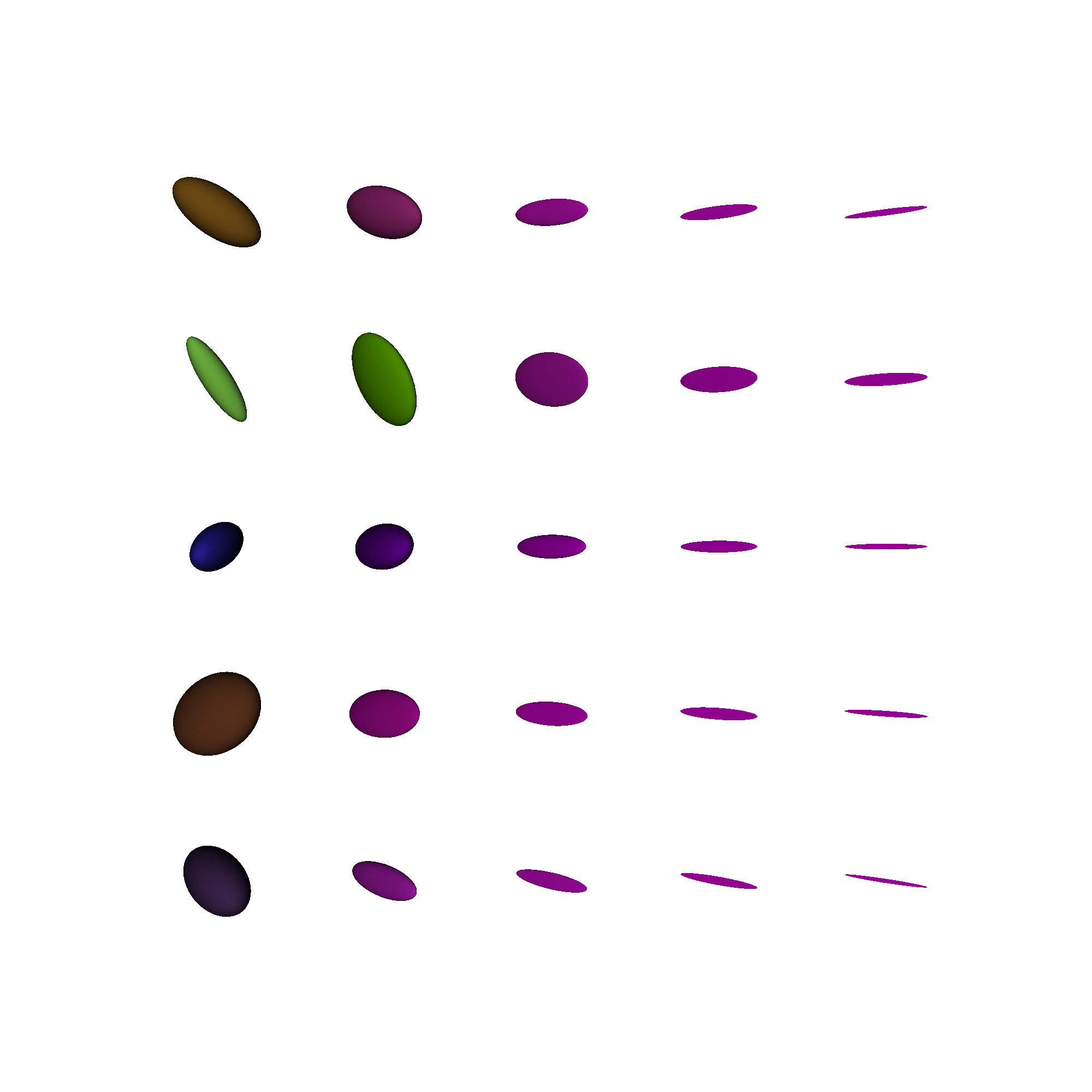}}
    %\hspace{0.03\textwidth}
	{\includegraphics[width=0.32\textwidth]{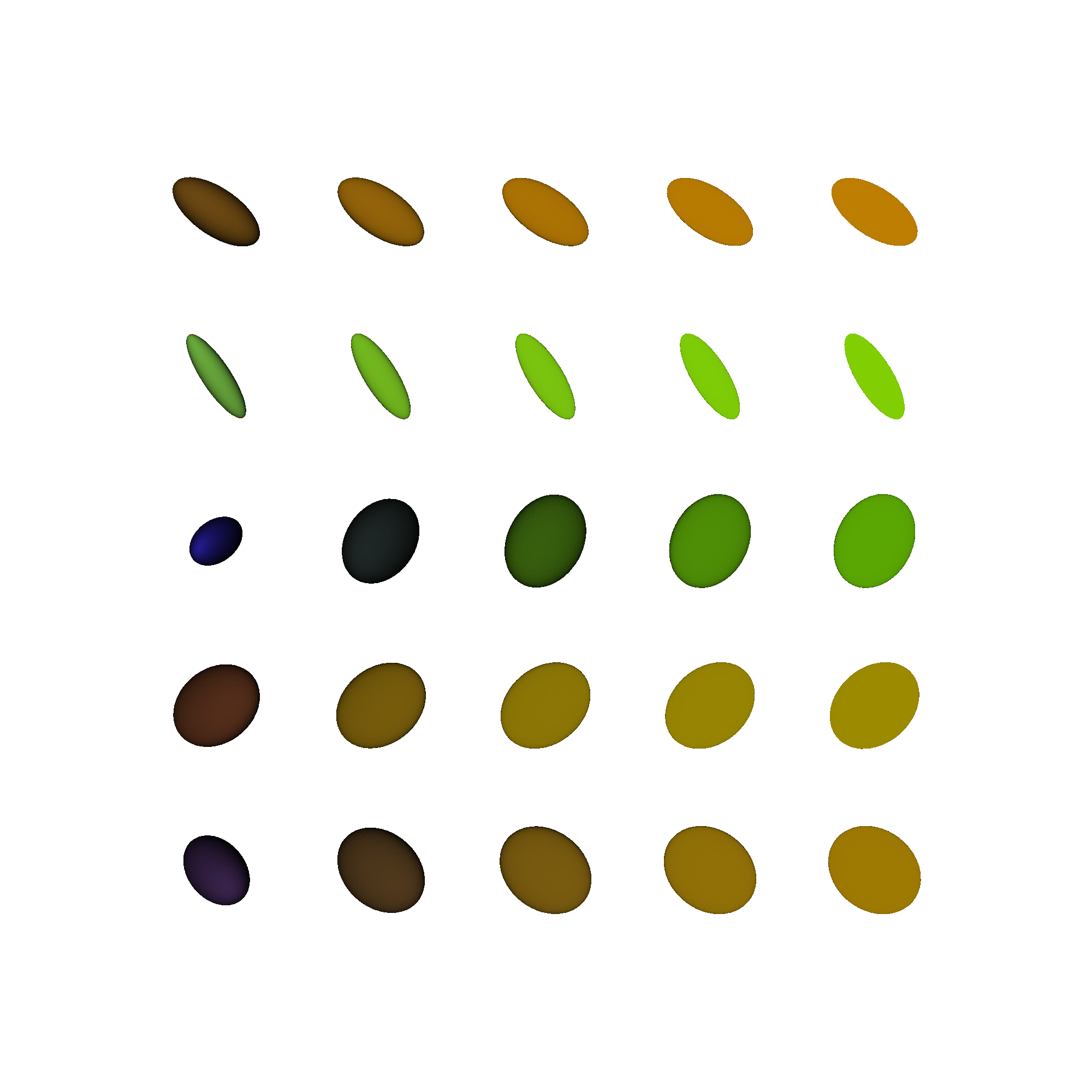}}
    %\hspace{0.03\textwidth}
    {\includegraphics[width=0.32\textwidth]{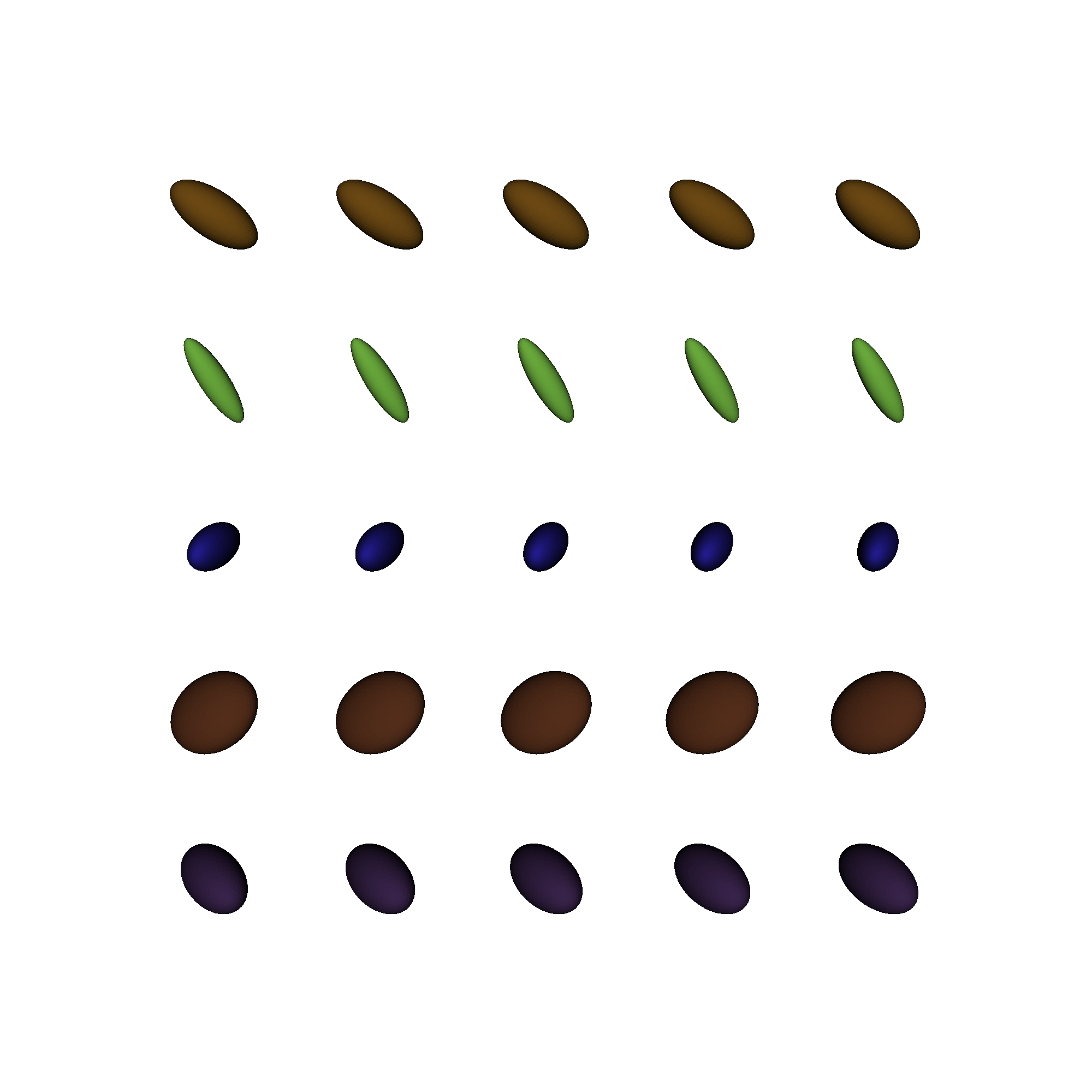}}
    \vspace{-5mm}
    \caption{\small{The $i$th row ($i=1,\ldots,5$) of the $j$th subfigure ($j=1,2,3$) represents values along the ray in $\xi_j$ issuing from $x_i$ at times $t=0,1,2,3,4$.}}
    \label{fig:xis}
\end{figure}

	In applications, a direction $\xi\in\partial M$ can often be given a concrete interpretation by observing the rays in $\xi$. Consider DTI, which explores brain structure via the diffusion of water molecules. The diffusion at a location in the brain is represented by a $3\times 3$ SPD matrix; the space of all such matrices, $\mathcal{P}_3$, is a 6-dimensional Hadamard manifold. Elements of $\mathcal{P}_3$ can be diagonalised and visualised as ellipsoids, with eigenvectors and eigenvalues determining the directions and lengths of the semi-axes, respectively. In Figure \ref{fig:xis}, DTI data are visualised in this way, resised to be on the same scale, and coloured according to fractional anisotropy, where brightness indicates eccentricity: Darker ellipsoids are closer to spheres. $\xi_1,\xi_2,\xi_3\in\partial \mathcal{P}_3$ and $x_1,\ldots,x_5\in \mathcal{P}_3$ were chosen, and the $i$th row of the $j$th subfigure represents the ray in $\xi_j$ issuing from $x_i$. Thus, $\xi_1$ can be interpreted as relatively greater diffusion in the left-right direction. $\xi_2$ indicates flattening in the direction perpendicular to the page, interpreted as relatively greater diffusion in the plane of the page. $\xi_3$ represents greater diffusion in all directions by the same factor; because the ellipsoids in Figure \ref{fig:xis} are scaled, the elements in any row appear identical.
	
	\section{Definition and basic properties} \label{def}
	
	Given a Hadamard space $(M,d)$ equipped with its metric topology and Borel $\sigma$-algebra $\mathcal{B}$, let $X$ be a random element in $M$, that is, a measurable map from some probability space $(\Omega,\mathcal{F},P)$ into $(M,\mathcal{B})$. We adopt in this paper the convention that if a function $f$ is 0 at $y$ and a function $g$ is infinite or undefined at $z$, $f(y)g(z)=0$. This is relevant in the following definition since $\angle_p(x,\xi)$ is undefined when $x=p$.
	\begin{definition} \label{quantile}
		For $\beta\in[0,1)$ and $\xi\in\partial M$, the \textit{$(\beta,\xi)$-quantile set} of $X$ is defined to be $q(\beta,\xi)=\arg\min_{p\in M}G^{\beta,\xi}(p)$, where 
		\begin{equation} \label{expect}
		G^{\beta,\xi}(p)=E[\rho(X,p;\beta,\xi)],
		\end{equation}
		$\rho(x,p;\beta,\xi)=d(p,x)+\beta d(p,x)\cos(\angle_p(x,\xi))$ and $\angle_p(x,\xi)$ is the Alexandrov angle at $p$ between the geodesic from $p$ to $x$ and the geodesic ray that is the unique member of $\xi$ issuing from $p$. An element of $q(\beta,\xi)$ is called a \textit{$(\beta,\xi)$-quantile} of $X$. If $q(\beta,\xi)$ is a singleton, in a slight abuse of notation, we use $q(\beta,\xi)$ to refer to its single element.
	\end{definition}
	
For definitions of terms like Alexandrov angle, refer to Section \ref{geomsupp} of the supplement. For any $\xi$, setting $\beta=0$ yields the geometric median. Proposition \ref{measurable} shows that (\ref{expect}) is well-defined; see Section \ref{proof_meas} of the supplement for a proof.
	
	\begin{proposition} \label{measurable}
		For fixed $\beta$, $\xi$, and $p$, the map $x\mapsto\rho(x,p;\beta,\xi)$ is measurable. %Therefore, $\rho(X,p;\beta,\xi)$ is a non-negative random variable, and it makes sense to take its expected value.
	\end{proposition}

	\begin{definition}
		Given data points $X_1,\ldots,X_N\in M$, the \textit{sample $(\beta,\xi)$-quantile (set)} is defined to be $\hat{q}_N(\beta,\xi)=\arg\min_{p\in M}\hat{G}^{\beta,\xi}_N(p)$, where $\hat{G}^{\beta,\xi}_N(p)=\frac{1}{N}\sum_{i=1}^N\rho(X_i,p;\beta,\xi)$.
	\end{definition}
	
	For some $\sigma>0$, we define a $\sigma$-scaled isometry to be a bijection $g:(M,d)\rightarrow (M',d')$  between metric spaces such that $d'(g(x),g(y))=\sigma d(x,y)$ for all $x,y\in M$ (an isometry is a 1-scaled isometry). For any geodesic $\gamma$ into $M$, define $g\gamma$ by $g\gamma(t)=g(\gamma(t/\sigma))$, its domain being the set of all $t$ such that $t/\sigma$ is in the domain of $\gamma$. Note that if $g$ is a $\sigma$-scaled isometry, $g^{-1}$ is a $\sigma^{-1}$-scaled isometry. 
 
 Throughout this paper, for a set $A$ and a function $f$, let $f(A):=\{f(a)|a\in A\}$. Quantiles on Hadamard spaces are equivariant to scaled isometries. %in the following sense.
	
	\begin{proposition} \label{equivariance}
		Given Hadamard spaces $M'$ and $M$ and $\sigma$-scaled isometry $g:M\rightarrow M'$, the $(\beta,g\xi)$-quantile set of $g(X)$, where $g\xi:=[g\gamma|\gamma\in\xi]\in\partial M'$, is $g(q(\beta,\xi))$.
	\end{proposition}
	
	See Section \ref{proof_meas} of the supplement for a proof. This generalises Facts 2.2.1 and 2.2.2 of \citet{Chaudhuri1996}. This is an advantage of the geometric approach because other non-geometric attempts to define multivariate quantiles exist but lack this equivariance; for example, coordinate-wise medians are not rotation equivariant.
	
	\subsection{Quantiles on Hadamard manifolds}
	%\begin{remark}
	Recall from Section \ref{geom} that all Riemannian manifolds are assumed to be $C^\infty$. Let $M$ be a Hadamard manifold. For $v_1,v_2\in T_pM$, $\langle v_1,v_2\rangle$ and $\lVert v_1\rVert$  denote the Riemannian inner product of $v_1$ and $v_2$ and Riemannian norm of $v_1$, respectively. On Riemannian manifolds, the Alexandrov and Riemannian angles coincide. Thus,
	\begin{equation}\begin{aligned} \label{riemquantile}
	G^{\beta,\xi}(p)=E[\lVert \log_p(X)\rVert+\langle \beta\xi_p,\log_p(X)\rangle].
	\end{aligned}\end{equation}
	Here, we have used the fact that on $M$, the domain of the inverse exponential map $\log_p$ is all of $M$ for all $p\in M$. Making the natural identification between $S^{n-1}$ and $\partial \mathbb{R}^n$ and letting $u=\beta\xi$, (\ref{riemquantile}) and (\ref{multiquantile}) are equal when $M=\mathbb{R}^n$.
	
	\begin{remark} \label{c2}
		The vector field $p\mapsto (p,\xi_p)$ is called the radial field in the direction of $\xi$ (\citet{Heintze1977}, \citet{Shcherbakov1983}). It is the negative gradient of the Busemann function at $\xi$ and is thus orthogonal to the level sets of this function, called horospheres. Under regularity conditions on $M$ that, roughly speaking, prevent the curvature from changing too much or too quickly, the $C^2$ differentiability of radial fields can be ensured; recall that a function is $C^2$ if it is twice continuously differentiable. 
		%For example, defining $\Delta<0$ and $\delta<0$ to be the $\inf$ and $\sup$, respectively, of all the sectional curvatures on $M$, Proposition 1.5 of \citet{Green1974} shows that if the covariant differential (denoted $\nabla$) of the curvature tensor $R$ is bounded and $0>\delta>4\Delta$, radial fields are $C^2$. Hyperbolic space clearly satisfies these regularity conditions, and in Proposition \ref{hyper}, we show that the radial fields are $C^\infty$ on this space. 
		For example, defining $\delta<0$ to be the $\sup$ of all the sectional curvatures on $M$, \citet{Shcherbakov1983} showed that if the curvature tensor $R$ and its covariant differential are bounded and $\delta<0$, the radial fields are $C^2$. Hyperbolic space clearly satisfies these regularity conditions, and in Proposition \ref{hyper}, we show that the radial fields are $C^\infty$ on this space.
		%This is a stronger result than that of \citet{Green1974} because the boundedness of the sectional curvatures implies the boundedness of $R$ (the second theorem in \citet{Karcher1970} gives an explicit bound for the norm of $R$ in terms of $\delta$ and $\Delta$). 
		However, these conditions do not completely characterise the Hadamard manifolds on which radial fields are $C^2$. For example, Euclidean space $\mathbb{R}^n$ does not satisfy these conditions, even though its radial fields are constant. Another example is the space of SPD matrices, for which $\delta$ can be shown to be 0 (see the explicit expression for $R$ in Proposition 2.3 of \citet{Dolcetti2018}), but which is known to have $C^\infty$ radial fields (see Corollary 2.4 of \citet{Shin2026}).
	\end{remark}
	
	\begin{proposition} \label{basic}
		On a Hadamard manifold $M$, assume that $G^{\beta^*,\xi^*}(p^*)$ is finite for some $\beta^*\in[0,1),\xi^*\in\partial M, p^*\in M$. Then, for every $(\beta,\xi)\in[0,1)\times \partial M$, (a) $G^{\beta,\xi}$ is finite and continuous on $M$, and (b) the $(\beta,\xi)$-quantile set is nonempty and compact.
	\end{proposition}
	A proof is provided in Section \ref{proof_prop3.4} of the supplement. Based on this result, the following condition is assumed in the rest of this paper.

  \begin{assumption} \label{assump}
For some $\beta^*\in[0,1),~\xi^*\in\partial M$, and $p^*\in M$, ~$G^{\beta^*,\xi^*}(p^*)<\infty$.
  \end{assumption}
	
	The following theorem, proved in Section \ref{proof_thm3.1} of the supplement, is used in the proof of the joint asymptotic normality of sample quantiles, but it is interesting in its own right as it provides a necessary condition for quantiles involving the gradient of $\rho$ by generalising Theorem 2.1.2 of \citet{Chaudhuri1996}. Denote by $\nabla \rho(x,p;\beta,\xi)$ the Riemannian gradient of the map $p'\mapsto\rho(x,p';\beta,\xi)$ at $p'=p\neq x$. Letting $I$ denote an indicator function, we adopt the notation $E[X;X\in A]:=E[XI(X\in A)]$. %Proofs for the following results are provided in Section \ref{proof_thm3.1} of the supplement.
 
	\begin{theorem} \label{qgrad}
		Let $M$ be a Hadamard manifold and $X$ be an $M$-valued random element with bounded support. If $q\in q(\beta,\xi)$, then $\lVert E[\nabla \rho(X,q;\beta,\xi);X\neq q]-P(X=q)\beta\xi_{q}\rVert\leq P(X=q)$. Thus, if $P(X=q)=0$, then
		%\begin{equation}
		$E[\nabla \rho(X,q;\beta,\xi);X\neq q]=0$. 
		%\end{equation}
	\end{theorem}
	
	%\begin{corollary} \label{sampleqgrad}
	%	Let $M$ be a Hadamard manifold and $X_1,\ldots,X_N$ be points in $M$. If $\hat{q}_N\in \hat{q}_N(\beta,\xi)$, then $\lVert \sum_{i:1\leq i\leq N,X_i\neq \hat{q}_N}\nabla \rho(X_i,\hat{q}_N;\beta,\xi)-L_{\hat{q}_N}\beta\xi_{\hat{q}_N}\rVert\leq L_{\hat{q}_N}$, where $L_{\hat{q}_N}:=\sum_{i=1}^N I(X_i=\hat{q}_N)$. Thus if $L_{\hat{q}_N}=0$, then
		%\begin{equation}
	%	$\sum_{i=1}^N\nabla \rho(X_i,\hat{q}_N;\beta,\xi)=0$. 
		%\end{equation}
	%\end{corollary}
	
We conclude this section with a remark on dimensionality. Any dimension reduction technique from a Hadamard manifold $M$ to a lower-dimensional Hadamard manifold $M'$ can easily be incorporated into our framework, as we can calculate the quantiles of the dimension-reduced data in $M'$. Any dimension reduction technique onto a complete connected totally geodesic submanifold preserves ``Hadamardness'' (and hyperbolicity) in this way; crucially, the principal geodesic analysis of \citet{Fletcher2004}, a generalisation of PCA and among the most important dimension reduction techniques on manifolds, is an example of such a technique.

	\section{Asymptotic properties on Hadamard manifolds} \label{lst}
	 In this section, $X,X_1,X_2,\ldots$ are i.i.d. $M$-valued random elements. %We also assume $G^{\beta,\xi}(q^*)$ is finite for some $q^*\in M$.
	
	\subsection{Strong consistency}
	In the following theorem and its proof, provided in Section \ref{proof_thm4.1} of the supplement, we adopt the notation $d(p,A):=\inf_{p'\in A}d(p,p')$ for $p\in M$, $A\subset M$.
	\begin{theorem} \label{slln}
 Let $M$ be a Hadamard manifold. If $X$ has a unique $(\beta,\xi)$-quantile, then any measurable choice from the sample $(\beta,\xi)$-quantile set of $X_1,\ldots,X_N$ converges almost surely to the $(\beta,\xi)$-quantile of $X$. %If $X$ has a unique $(\beta,\xi)$-quantile, then with probability 1, any sequence of random elements where the $N$-th element is in the sample $(\beta,\xi)$-quantile set of $X_1,\ldots,X_N$ for all $N\in\mathbb{Z}_+$ converges to the $(\beta,\xi)$-quantile of $X$.
	\end{theorem}
 \begin{remark}\label{unique}
% Regarding the uniqueness of the quantile, we have verified empirically that the quantile loss functions, unlike those for the mean and median (\citet{Karcher1970}, \citet{Yang2010}), need not be convex as a function of $p$. However, based on this empirical evidence, we suspect that $\rho(x,\cdot;\beta,\xi)$ may be quasiconvex, and strictly so when the distribution is not supported on a geodesic. We have also empirically verified that it is not necessarily quasiconvex for $\beta\in(-1,0)$, which is part of the reason that $\beta$ has been restricted to $[0,1)$.
Strict convexity of $G^{\beta,\xi}$ guarantees the uniqueness of its local minimum, and thus the $(\beta,\xi)$-quantile. \citet{Chaudhuri1996} showed that $G^{\beta,\xi}$ is strictly convex in the Euclidean case if the distribution is not supported on a line. \citet{Yang2010} showed on Hadamard manifolds that $G^{0,\xi}$ is strictly convex if the distribution is not degenerate; that is, not supported on a geodesic. We have verified that $G^{\beta,\xi}$ need not be convex in general, but the strict convexity is obviously a much stronger condition than the uniqueness of local, let alone global, minima. We have never observed a non-degenerate case in which even the local minima are non-unique, and we suspect that quantiles, such as medians, are, in fact, unique in such cases.
 \end{remark}

	\vspace{-5mm}

	\subsection{Joint asymptotic normality}\label{asymp}
 
	In this section, we investigate the joint asymptotic normality of sample quantiles on Hadamard manifolds. All proofs are provided in Section \ref{proof_thm4.2} of the supplement. 
	
	Any $C^\infty$ global chart $\phi$ for an $n$-dimensional manifold $M$ is a diffeomorphism between $M$ and an open subset of $\mathbb{R}^n$. For example, the inverse exponential map at any $p\in M$ can define such a $\phi$. %In this section and the associated proofs in Section \ref{proof_thm4.2} of the supplement, $\lVert \cdot\rVert_2$ and $\langle\cdot,\cdot\rangle_2$ denote the $L_2$ norm and standard Euclidean inner product, respectively, in $\phi(M)$, and $\lVert \cdot\rVert$ and $\langle\cdot,\cdot\rangle$ denote the Riemannian norm and inner product, respectively. We also make reference to the $\sup$, or $L_\infty$, norm of a vector with $\lVert\cdot\rVert_\infty$, and the Frobenius norm of a matrix with $\lVert\cdot\rVert_F$.
	%In this section and the associated proofs in Section \ref{proof_thm4.2} of the supplement, in addition to the already-introduced Euclidean inner product $\langle \cdot,\cdot\rangle_2$ and norm $\lVert \cdot\rVert_2$, and Riemannian inner product $\langle \cdot,\cdot\rangle$ and norm $\lVert \cdot\rangle_2$, we also make reference to the $\sup$, or $L_\infty$, norm of a vector with $\lVert\cdot\rVert_\infty$, and the Frobenius norm of a matrix with $\lVert\cdot\rVert_F$.
	This section makes use of both the Euclidean inner product $\langle \cdot,\cdot\rangle_2$ and norm $\lVert \cdot\rVert_2$, and Riemannian inner product $\langle \cdot,\cdot\rangle$ and norm $\lVert \cdot\rVert$. 
  
	By standard abuse of notation, this section identifies $M$ with its image under $\phi$ so that all points $p$ in $M$ and tangent vectors in $T_pM$ are identified with their $n\times 1$ local coordinate representations in the chart $\phi$ and can thus be thought of as vectors in $\mathbb{R}^n$. This identifies the tangent bundle $TM$ with $\phi(M)\times\mathbb{R}^n$. Recall the so-called Riemannian metric, which defines an inner product on $T_pM$ at each point $p\in M$. Under a chart, the Riemannian metric at $p$ can be represented by an $n\times n$ symmetric positive definite matrix $g_p$, which defines the inner product at $p$ by $\langle v_1,v_2\rangle=v_1^Tg_pv_2$ for all $v_1,v_2\in T_pM$. Thus, the standard inner product in $\mathbb{R}^n$ can be defined by letting $g_p$ be the $n\times n$ identity for all $p$. Then, define $\Psi(x,p;\beta,\xi)\in T_pM$, for $x,p\in M$, by
	\vspace{-6mm}
	\begin{align*}
	\Psi(x,p;\beta,\xi)=(\Psi^1(x,p;\beta,\xi),\ldots,\Psi^n(x,p;\beta,\xi))^T 
	=\begin{cases}
	g_p\nabla \rho(x,p;\beta,\xi) &\text{ if $x\neq p$}, \\
	-\beta g_p\xi_p &\text{ if $x=p$},
	\end{cases}
	\end{align*} 
	where $\nabla \rho(x,p;\beta,\xi)$ is as defined in Theorem \ref{qgrad}. Let $df_p$ denote the Riemannian differential of a real-valued function $f:M\rightarrow\mathbb{R}$ at $p\in M$. For $x\neq p$ and $v\in T_pM$, 
	\begin{equation}\begin{aligned} \label{equiv}
	\langle v,\Psi(x,p;\beta,\xi)\rangle_2=\langle v,\nabla \rho(x,p;\beta,\xi)\rangle=d(\rho(x,\cdot;\beta,\xi))_p(v),
	\end{aligned}\end{equation}
	and therefore $\Psi(x,p;\beta,\xi)$ is the Euclidean gradient at $p\neq x$ of $\rho(x,\cdot;\beta,\xi)$. Assuming this function is $C^2$ when $p\neq x$, denote its (Euclidean) Hessian by $D\Psi(x,p;\beta,\xi)$ and, for $r,r'=1,\ldots,n$, the $(r,r')$-entry of this matrix by $D_{r'}\Psi^r(x,p;\beta,\xi)$. $D\Psi(x,p;\beta,\xi)$ is the Jacobian of $\Psi(x,\cdot;\beta,\xi)$ at $p\neq x$. Let $D\Psi(x,p;\beta,\xi)$ be 0 when $p=x$.
	
	Given $(\beta_1,\xi_1),\ldots,(\beta_K,\xi_K)\in[0,1)\times\partial M$, denote $\Psi(x,p;\beta_k,\xi_k)$, $D\Psi(x,p;\beta_k,\xi_k)$ and $D_{r'}\Psi^r(x,p;\beta_k,\xi_k)$ by the shorthands $\Psi_k(x,p)$, $D\Psi_k(x,p)$ and $D_{r'}\Psi_k^r(x,p)$, respectively ($k=1,\ldots,K$). We use $q_k$ and $\hat{q}_{k,N}$ to denote the quantile $q(\beta_k,\xi_k)$ and some measurable selection from the sample quantile set $\hat{q}_N(\beta_k,\xi_k)$, respectively. In addition, we define, $\lambda_k(Y):=E[\Psi_k(X,Y)|Y]$, where $Y$ may be fixed or random.

	\begin{theorem} \label{clt}
\sloppy Let $M$ be an $n$-dimensional Hadamard manifold, $n\geq2$, and $(\beta_1,\xi_1),\ldots,(\beta_K,\xi_K)\in[0,1)\times\partial M$ be such that the functions $p\mapsto \langle \beta_k(\xi_k)_p,\log_p(x)\rangle$ ($k=1,\ldots,K$, $x\in M$) are twice continuously differentiable at $p\neq x$.		
  
  Suppose $X$ satisfies the following for each $k=1,\ldots,K$: (I) $q_k$ exists uniquely, $\lambda_k(q_k)=0$, and there is some neighbourhood $Q_k\subset M$ around $q_k$ in which the density of $X$ exists and is bounded, (II) $E[\lVert\Psi_k(X,q_k)\rVert_2^2]<\infty$, and (III) there exists a $d_1>0$ and neighbourhood $U\subset M$ of $q_k$ that is bounded in the Riemannian metric such that $\{q:\lVert q-q_k\rVert_2\leq d_1\}\subset U$ and $E[\sup_{q:\lVert q-q_k\rVert_2\leq d_1}\lvert D_{r'}\Psi_k^r(X,q)\rvert^2;X\not\in \bar{U}]<\infty$ for all $r,r'=1,\ldots,n$. Then, $D_{r'}\Psi_k^r(X,q_k)$ is integrable for each $k=1,\ldots,K$ and $r,r'=1,\ldots,n$. Additionally, suppose the following is satisfied for each $k=1,\ldots,K$: (IV) $E[D\Psi_k(X,q_k)]$ is nonsingular.
   
		\sloppy Then, $\sqrt{N}([\hat{q}_{1,N}^T,\ldots,\hat{q}_{K,N}^T]^T-[q_1^T,\ldots,q_K^T]^T)\rightsquigarrow N(0,\Lambda^{-1}\Sigma\Lambda^{-1})$ as $N\rightarrow\infty$, where $\Lambda$ is the $Kn\times Kn$ block diagonal matrix whose $k$-th $n\times n$ diagonal block ($k=1,\ldots,K$) is $E[D\Psi_k(X,q_k)]$, and $\Sigma$ is the covariance matrix of $[\Psi_1(X,q_1)^T,\ldots,\Psi_K(X,q_K)^T]^T$.
	\end{theorem}

See Remark \ref{c2} for a discussion on Hadamard manifolds on which radial fields are $C^2$, and, therefore, all maps $p\mapsto\langle \beta_k(\xi_k)_p,\log_p(x)\rangle$, are $C^2$ at $p\neq x$. If $K=1$ and $\beta_1=0$, the $C^2$ condition is satisfied on any Hadamard manifold. 

We can show that $\sup_{q:\lVert q-q_k\rVert_2\leq d_1}\lvert D_{r'}\Psi_k^r(X,q)\rvert^2 I(X\not\in \bar{U})$ from (III) is indeed a random variable and its expected value is therefore well-defined. For this and other remarks on Theorem \ref{clt}, see Section \ref{remarks} of the supplement.

% \begin{remark} \label{lsc}
    % We can show that $h_{r,r'}(\eta,w,X)$ is indeed a random variable for fixed $\eta\in(0,\infty)$ and $w\in S^{n-1}$ and its expected value is, therefore, well-defined. Define a function $H$ by
%			\begin{equation}\begin{aligned} \label{H}
%			H(x,q)=\lvert D_{r'}\Psi_k^r(x,q)-D_{r'}\Psi_k^r(x,q_k)\rvert.
%			\end{aligned}\end{equation}
%			For a fixed $x=x_0\not\in q_k+[0,\eta]\times w$, this is continuous as a function of $q$ on the compact set $q_k+[0,\eta]\times w$, so there exists some $s_0\in[0,\eta]$ for which $H(x_0,q_k+s_0w)=\sup_{q:q_k+[0,\eta]\times w}H(x_0,q)=h_{r,r'}(\eta,w,x_0)$. For a fixed $q\in q_k+[0,\eta]\times w$, $H$ is also continuous as a function of $x$ at $x_0$, and hence, we can say that
%			\begin{align*}
%			\lim\inf_{x\rightarrow x_0}h_{r,r'}(\eta,w,x)&\geq\lim\inf_{x\rightarrow x_0}H(x,q_k+s_0w)\\
%			&=\lim_{x\rightarrow x_0}H(x,q_k+s_0w)=H(x_0,q_k+s_0w)=h_{r,r'}(\eta,w,x_0);
%			\end{align*}
%			that is, $h_{r,r'}(\eta,w,x)$ is lower semi-continuous as a function of $x$ outside of $q_k+[0,\eta]\times w$ and 0 in it, making it measurable as a function of $x$. The lower semi-continuity of $\sup_{q:\lVert q-q_k\rVert_2\leq d_1}\lvert D_{r'}\Psi_k^r(x,q)\rvert^2 I(x\not\in \bar{U})$ as a function of $x$ can be shown in a similar manner.
% \end{remark}

The next result provides simple, easy-to-understand conditions that guarantee joint asymptotic normality. %, based on the observation that if $X$ has bounded support according to the Riemannian metric, we need not worry about (II) and (III). The requirement on $\lambda$ in (I) is also no longer present.
	
	\begin{corollary} \label{bdd}
Let all terms be as in the first paragraph of Theorem \ref{clt} and $X$ have bounded support according to the Riemannian metric. Suppose that for each $k=1,\ldots,K$, $q_k$ exists uniquely, and there is some neighbourhood around $q_k$ in which the density of $X$ exists and is bounded. Then, $D_{r'}\Psi_k^r(X,q_k)$ is integrable for each $k=1,\ldots,K$ and $r,r'=1,\ldots,n$. Additionally, suppose (IV) in Theorem \ref{clt} holds. Then, the conclusion of Theorem \ref{clt} follows.
	\end{corollary}

    The following proposition also provides conditions under which we can specify the convergence rate of sample quantiles.
    
\begin{proposition}\label{speed}
    If $\sqrt{N}(\hat{q}_{k,N}-q_k)$ converges weakly to some random $n$-vector as $N\rightarrow\infty$, then $d(\hat{q}_{k,N},q_k)=O_p(1/\sqrt{N})$.
\end{proposition}
	
\section{Computing quantiles on hyperbolic spaces} \label{hyp}	
 The most obvious method for computing quantiles on a Hadamard manifold $M$ is a gradient descent algorithm. Given $(\beta,\xi)\in[0,1)\times\partial M$, a data set $X_1,\ldots,X_N\in M$ and an initial estimate $\hat{q}_{[0]}$ for the sample quantile, the $(j+1)$th estimate for non-negative integer $j$ is $\hat{q}_{[j+1]}=\exp_{\hat{q}_{[j]}}(-\alpha_j\sum_{i=1}^N\nabla \rho(X_i,\hat{q}_{[j]};\beta,\xi))$, where $\hat{q}_{[j]}$ is the $j$th estimate and $\alpha_j$ is a learning rate that may be updated. The algorithm terminates when the update sises are sufficiently small. In the Euclidean case, $\nabla \rho(x,p;\beta,\xi)$, $p,x\in M$, is quite straightforward: $-(x-p)/\lVert x-p\rVert-\beta\xi$, assuming $p\neq x$. We might hope that the gradient in the manifold case, again assuming $p\neq x$, is the obvious analogue $-\log_p(x)/d(p,x)-\beta\xi_p$; in fact, it is significantly more complicated. In this section, we provide an explicit expression for $\nabla \rho(x,p;\beta,\xi)$ when $M$ is a hyperbolic space; the proof is provided in Section \ref{proof_thm5.1} of the supplement. For more details on hyperbolic spaces, including Poincar\'e balls, see Section \ref{geomsupp} of the supplement. 
	
	%The gradient of $\rho(x,p;\beta,\xi)$ is quite straightforward in the Euclidean case: $-(x-p)/\lVert x-p\rVert-u$, assuming $x\neq p$. We might hope that the gradient in the manifold case, again assuming $x\neq p$, is the obvious analogue $-\log_p(x)/d(p,x)-\beta\xi_p$; in fact, it is significantly more complicated.
	
	\begin{theorem}\label{thm5.1}
		Let $M$ be a hyperbolic space of sectional curvature $\kappa<0$. $\nabla \langle \xi_p,\log_p(x)\rangle$ denotes the gradient of the map $p'\mapsto \langle \xi_{p'},\log_{p'}(x)\rangle$ at $p'=p$. For $p=x$, $\nabla \langle \xi_p,\log_p(x)\rangle=-\xi_p$, and for $p\neq x$,
		\begin{equation*}\begin{aligned}
		&\nabla \langle \xi_p,\log_p(x)\rangle=-\sqrt{-\kappa}d(p,x)\bigg(\coth(\sqrt{-\kappa}d(p,x))-\bigg\langle\xi_p,\frac{\log_p(x)}{d(p,x)}\bigg\rangle\bigg)\xi_p  \\
		&\quad-\bigg((1-\sqrt{-\kappa}d(p,x)\coth(\sqrt{-\kappa}d(p,x)))\bigg\langle \xi_p ,\frac{\log_p(x)}{d(p,x)}\bigg\rangle+\sqrt{-\kappa}d(p,x)\bigg) \frac{\log_p(x)}{d(p,x)}.
		\end{aligned}\end{equation*}
	\end{theorem}

		 For $p\neq x$, the gradient of $d(p,x)$ with respect to $p$ is $-\log_p(x)/d(p,x)$, and thus, $\nabla \rho(x,p;\beta,\xi)=-\log_p(x)/d(p,x)+\beta\nabla \langle \xi_p,\log_p(x)\rangle$.

\begin{remark}
	Although $\coth(\sqrt{-\kappa}d(p,x))$ and $1/d(p,x)$ both tend to infinity as $p$ approaches $x$, $\nabla \langle \xi_p,\log_p(x)\rangle$ does not, as explained in the last paragraph of the proof of Theorem \ref{thm5.1}; instead, it approaches $-\xi_x$ as $p\rightarrow x$, and by Theorem \ref{thm5.1}, it equals $-\xi_x$ when $p=x$. 
	%Thus to avoid any potential numerical instability, in our algorithm we simply approximate $\nabla \langle \xi_p,\log_p(x)\rangle\rho(x,p;\beta,\xi)$ by $-\xi_p$ for values of $p$ extremely close to $x$. 
	As for the gradient of $p\mapsto d(p,x)$, this function is not differentiable when $p=x$, but for reasons of symmetry, it is natural to treat this gradient as 0 when $p=x$, and for reasons of stability, also when $p$ is very near $x$. Thus, our algorithm avoids potential instability by updating using $-\xi_p$ when $p$ is extremely close to $x$.
\end{remark}
\begin{remark}
	It is clear from Theorem \ref{thm5.1}, the expressions for $\xi_p$ in Section \ref{proof_thm5.1} and $\log_p$, and the Minkowski pseudo-inner product in Section \ref{geomsupp} of the supplement, that calculating $\nabla \langle \xi_p,\log_p(x)\rangle$ requires $O(n+1)=O(n)$ operations. These and the expression for $\exp_p$ in Section \ref{geomsupp} show that each iteration of the algorithm requires $O(nN)$ operations, so denoting the number of iterations by $T$, the complexity of the algorithm is $O(nNT)$, and seems to scale linearly with dimension and sample sise. As for the complexity of the number of iterations $T$ until convergence, in general, gradient descent algorithms do not seem to scale with $n$ or $N$.
\end{remark}

	%\begin{remark}
	%	Because $\coth$ goes to infinity with its argument, the norm of this gradient goes to infinity as $x$ goes to infinity along a geodesic ray beginning at $p$ as long as $\log_p(x)/d(p,x)\neq\pm\xi^p$. So unlike in Euclidean spaces, on Hadamard manifolds there is in general no guarantee that $\lVert \nabla \rho(x,p;\beta,\xi)\rVert$ is contained in $[1-\beta,1+\beta]$.
	%\end{remark}
	
	\section{Discussion of some possible applications} \label{applications}

 \subsection{Isoquantile contours, outlier detection, and a transformation-retransformation procedure} \label{tr}
Let $M$ be a Hadamard manifold. For a fixed $\beta\in[0,1)$, we call the set of all $(\beta,\xi)$-quantiles, $\xi\in\partial M$, the $\beta$-isoquantile contour. Isoquantile contours can give a sense of the shape of a distribution, but \citet{Girard2017} observed that as $\beta\rightarrow 1$, Euclidean geometric quantiles are most distant from the geometric median in the directions where the data vary the least. This problem can be addressed using the transformation-retransformation (TR) procedure of \citet{Chakraborty2001}. %Here, we give brief outlines of this procedure and a similar one on Hadamard manifolds, which can be used for outlier detection.

For iid $X_1,\ldots,X_N\in\mathbb{R}^n$, $N\geq n+1$, \citet{Chakraborty2001} suggested choosing $i_0,\ldots,i_n$ satisfying $1\leq i_0<\cdots<i_n\leq N$ so that the $n\times n$ matrix $A:=[X_{i_1}-X_{i_0},\ldots,X_{i_n}-X_{i_0}]$ is invertible. Then, the TR $u$-quantile with being $\lVert u\rVert_2<1$, is $Aq^*$, where $q^*$ is the $\Vert u\rVert_2(A^{-1}u/\lVert A^{-1}u\rVert_2)$-quantile of the transformed data set $A^{-1}X_1,\ldots,A^{-1}X_N$. This procedure calculates the isoquantile contour of the transformed data set for a given $\beta=\lVert u\rVert_2$ and retransforms it to the original coordinate system. Ideally, the $i_0,\ldots,i_n$ are chosen so that the covariance of $A^{-1}X_i$, treating $A$ as constant, is as isotropic as possible; this is done by minimising the ratio of the arithmetic mean to the geometric mean of the eigenvalues of $V^{-1}$ (that is, $\text{tr}(V^{-1})/(\det(V^{-1})^{1/n}$), where $V^{-1}=A^T\hat\Sigma^{-1}A$ and $\hat\Sigma$ is an estimate of the covariance of the distribution underlying the data $X_1,\ldots,X_N$.

\sloppy For $X_1,\ldots,X_N$ on an $n$-dimensional Hadamard manifold $M$, we might define an invertible matrix $A:=[\log_{X_{i_0}}(X_{i_1}),\ldots,\log_{X_{i_0}}(X_{i_n})]$, where vectors in $T_{X_{i_0}}M$ are identified with their representations in some orthonormal basis, calculate the $\beta$-isoquantile contour $Q^*(\beta)$ of $\exp_{X_{i_0}}(A^{-1}\log_{X_{i_0}}(X_1)),\ldots,\exp_{X_{i_0}}(A^{-1}\log_{X_{i_0}}(X_N))$ for a given $\beta$, and define the TR $\beta$-isoquantile contour as $\{\exp_{X_{i_0}}(A\log_{X_{i_0}}(q)):q\in Q^*(\beta)\}$. For simplicity, we can replace $X_{i_0}$ by the median $m_1$ of $X_1,\ldots,X_N$, which, as shown by \citet{Yang2010}, is unique if the data are not supported on a geodesic. Finally, $i_1<\cdots<i_n$ should be chosen to minimise $\text{tr}(A^T\hat\Sigma^{-1}A)/(\det(A^T\hat\Sigma^{-1}A)^{1/n})$, where $\hat\Sigma$ estimates the covariance of the distribution underlying $\log_{m_1}(X_1),\ldots,\log_{m_1}(X_N)$, treating $m_1$ as constant.

TR isoquantile contours can be used for outlier detection. An obvious method is to label as outliers all points on the outside of the TR $\beta_0$-isoquantile contour for some $\beta_0$ close to 1, but this method can be self-defeating as extreme quantiles can themselves be sensitive to outliers. Instead, consider a common rule of thumb with univariate data which designates points outside $[Q_{1/4}-1.5(Q_{3/4}-Q_{1/4}),~Q_{3/4}+1.5(Q_{3/4}-Q_{1/4})]$ as outliers, where $Q_{\tau}$, $\tau\in(0,1)$, is the univariate $\tau$-quantile. A similar approach (equivalent when $\lvert Q_{1/4}-Q_{1/2}\rvert=\lvert Q_{3/4}-Q_{1/2}\rvert$) could instead use $[Q_{1/2}-3\lvert Q_{1/4}-Q_{1/2}\rvert,~Q_{1/2}+3\lvert Q_{3/4}-Q_{1/2}\rvert]$. Considering quantiles for $\beta=0.5$ to be analogues of quartiles and with the same notation as in the previous paragraph, one generalisation would designate as outliers points outside $\{\exp_{m_1}(A\log_{m_1}(p)):p\in O\}$, where $O:=\{\exp_{m_1^*}(3\log_{m_1^*}(q)):q\in Q^*(0.5)\}$, and $m_1^*$ and $Q^*(0.5)$ are the median and 0.5-isoquantile contour, respectively, of the transformed data set. 

When using TR isoquantile contours for outlier detection, the covariance estimate $\hat\Sigma$ should be robust; otherwise, the entire method could become sensitive to outliers. When outliers are absent and the goal is simply to identify the shape of the distribution, the sample or plug-in covariance estimate is sufficient.

\subsection{Measures of distributional characteristics} \label{measures}

The centrality, dispersion, skewness, and kurtosis of a univariate distribution are typically measured using the first four moments. Quantiles can be used to define alternative measures of these distributional characteristics. Again, $Q_{\tau}$, $\tau\in(0,1)$, is the $\tau$-quantile for univariate data. The median $Q_{0.5}$ measures centrality, and $\delta_0(\beta):=Q_{(1+\beta)/2}-Q_{(1-\beta)/2}$ and $\gamma_0(\beta):=(Q_{(1+\beta)/2}+Q_{(1-\beta)/2}-2Q_{1/2})/\delta_0(\beta)$, $\beta\in(0,1)$, are the standard quantile-based measures of dispersion and skewness, respectively, typically with $\beta=1/2$. Kurtosis (or tailedness) can be measured by $\kappa_0(\beta,\beta'):=\delta_0(\beta')/\delta_0(\beta)$ for $0<\beta<\beta'<1$.  \citet{ShinOh2024} generalised these to multivariate distributions, thereby presenting quantile-based alternatives to existing moment-based measures of multivariate distributional characteristics, such as those of \citet{Mardia1970}, and also introduced a measure for spherical asymmetry. Here, we generalise these quantile-based measures to distributions on Hadamard manifolds with isoquantile contours and investigate them further in Section \ref{sexperiments}. 

For simplicity, assume the uniqueness of the geometric median and, for fixed $\beta,\beta'$ satisfying $0<\beta<\beta'<1$, all $(\beta,\xi)$- and $(\beta',\xi)$-quantiles. The median $m_1$ is the obvious measure of centrality. Define a bijection $f:S_{m_1}^{n-1}\rightarrow \partial M$, where $S_{m_1}^{n-1}$ is the unit sphere in $T_{m_1}M$, by $f(\xi_{m_1})=\xi$ and let $SA(n-1)$ be the surface area of the unit $(n-1)$-sphere. \citet{ShinOh2024} defined two measures each for dispersion, skewness, and kurtosis, and one for spherical asymmetry; we generalise these to
\begin{equation*}\begin{gathered}
\delta_1(\beta):=\sup_{v\in S_{m_1}^{n-1}}{\lVert \log_{m_1}(q(\beta,f(v)))-\log_{m_1}(q(\beta,f(-v)))\rVert}, \\
\delta_2(\beta):=\frac{1}{SA(n-1)}\int_{S_{m_1}^{n-1}}{\lVert \log_{m_1}(q(\beta,f(v)))-\log_{m_1}(q(\beta,f(-v)))\rVert} dv, \\
\gamma_1(\beta):=\frac{\sup_{v\in S_{m_1}^{n-1}}\lVert \log_{m_1}(q(\beta,f(v)))+\log_{m_1}(q(\beta,f(-v)))\rVert}{\delta_1(\beta)},\\
\gamma_2(\beta):=\frac{(1/SA(n-1))\int_{S_{m_1}^{n-1}}\log_{m_1}(q(\beta,f(v))) dv}{\delta_2(\beta)},\\
\kappa_1(\beta,\beta'):=\frac{\delta_1(\beta')}{\delta_1(\beta)}, ~~
\kappa_2(\beta,\beta'):=\frac{\delta_2(\beta')}{\delta_2(\beta)}, ~~
\alpha(\beta):=\log\bigg(\frac{\sup_{\xi\in\partial M}d(m_1,q(\beta,\xi))}{\inf_{\xi\in\partial M}d(m_1,q(\beta,\xi))}\bigg).
\end{gathered}\end{equation*}
By transforming the data as in Section \ref{tr} to have roughly isotropic covariance, $\alpha(\beta)$ can also be used to measure elliptical asymmetry. Because $\gamma_0(\beta)$ can be positive or negative, its generalisation to Hadamard manifolds should be a vector as $\gamma_2(\beta)$ is; $\gamma_1(\beta)$ generalises $\lvert\gamma_0(\beta)\rvert$ rather than $\gamma_0(\beta)$ itself. These measures can be used to test properties of distributions, such as whether a distribution is skewed or spherically symmetric, or whether one distribution is more spread out than another.

\subsection{Testing the equality of distributions, measuring distributional dissimilarity, and finding loci of differences} \label{testing}

\citet{Koltchinskii1997} demonstrated that geometric quantiles uniquely determine a multivariate distribution; even if that is not true in general on Hadamard spaces, we can test whether different data sets are generated by the same distribution by comparing corresponding quantiles. Given samples of sises $N$ and $N'$, denote their underlying distributions by $F$ and $F'$, and denote random elements following these distributions by $X$ and $X'$. Denote the respective sample means by $\hat m_N$ and $\hat m'_{N'}$ and, for $(\beta,\xi)\in[0,1)\times\partial M$, sample $(\beta,\xi)$-quantiles by $\hat q_N(\beta,\xi)$ and $\hat q_N'(\beta,\xi)$ and population $(\beta,\xi)$-quantiles by $q(\beta,\xi)$ and $q'(\beta,\xi)$. For real data, $\lvert \hat m_N-\hat m'_{N'}\rvert$ is a common statistic for testing $H_0:F=F'$. This generalises to $T_0=d(\hat m_N,\hat m'_{N'})$; alternatively, one could use the distance between the sample medians $T_1=d(\hat q_N(0,\xi),\hat q_N'(0,\xi))$. When different distributions have similar centres, $T_0$ and $T_1$ may not reliably reject $H_0$. On the other hand, given a positive integer $K$, $l>0$, and $(\beta_1,\xi_1,w_1),\ldots,(\beta_K,\xi_K,w_K)\in[0,1)\times \partial M\times[0,\infty)$, a new test statistic
\begin{align*}
T_2((\beta_1,\xi_1,w_1),\ldots,(\beta_K,\xi_K,w_K)):=\Big[\sum_{k=1}^Kw_kd(\hat q_N(\beta_k,\xi_k),\hat q'_{N'}(\beta_k,\xi_k))^l\Big]^{1/l}
\end{align*}
should distinguish $F$ and $F'$ if they differ at only one of the quantiles, regardless of how well quantiles represent distributional shape. We can consider $T_2$ as allowing for comparisons between distributions at higher resolutions than $T_0$ and $T_1$.

We also note a connection between $T_2$ and the Wasserstein distance. Let $[0,1)\times\partial M/\sim$ denote the quotient space, where $\sim$ is the equivalence relation defined by $(0,\xi)\sim (0,\xi')$ for all $\xi,\xi'\in\partial M$. Denote the element of this space containing $(\beta,\xi)$ by $[\beta,\xi]$. A population version of $T_2$, $[\sum_{k=1}^Kw_kd(q(\beta_k,\xi_k),q'(\beta_k,\xi_k))^l]^{1/l}$, can be used as a measure of the dissimilarity between distributions, and this measure can be made arbitrarily precise by adding quantiles. More generally, if we define a probability measure $\mu$ on $[0,1]\times\partial M/\sim$, then $[\int d(q(\beta,\xi),q'(\beta,\xi))^l \mu(d[\beta,\xi])]^{1/l}$ can be used to define a (pseudo-)metric between distributions on $M$. In particular, if $M=\mathbb{R}$ and $\mu$ is the uniform distribution on the interval $(-1,1)$, this quantity is precisely the $l$-Wasserstein distance on $\mathbb{R}$. Thus, in $\mathbb{R}$, both the population version of $T_2$ and the Wasserstein distance are instances of this pseudo-metric.

$T_2$ can also determine whether distributions differ at specific loci by restricting the quantile indices to a region of interest $R\subset[0,1)\times\partial M$. For example, if $R=\{(\beta_0, \xi_0)\}$, rejection indicates that the distributions differ at the $(\beta_0,\xi_0)$-quantile. Letting $R=\{(\beta,\xi):\beta\in R^*\}$, we can test whether the distributions differ at certain isoquantile contours determined by $R^*\subset[0,1)$; smaller values of $\beta$ correspond to the central masses of the distributions, and larger ones to the peripheries.

%$T_2$ can be used in permutation tests, but we also expect that Theorem \ref{clt} can be used to find the asymptotic distribution of $T_2$; we leave this for future research. 
When the samples are independent and consist of {\it i.i.d.} elements, another ``high resolution" quantile-based test for $H_0$, this time using asymptotic approximations, follows from Theorem \ref{clt}. Let $N'\equiv N'(N)$ be a deterministic sequence of positive integers satisfying $N'/N\rightarrow \eta$ as $N\rightarrow\infty$ for $\eta\neq 0$, and let $\hat q_{k,N}$, $\Lambda$, $\Sigma$ be as defined in Section \ref{asymp}, and $\hat q'_{k,N}$, $\Lambda'$, $\Sigma'$ be the equivalents for $X'$. Denote $[\hat q_{1,N}^T,\ldots,\hat q_{K,N}^T]^T$ and $[(\hat q'_{1,N'})^T,\ldots,(\hat q'_{K,N'})^T]^T$ by $\hat q_{N}$ and $\hat q'_{N'}$, respectively. Then if the conditions of Theorem \ref{clt} hold for both $X$ and $X'$, Slutsky's theorem implies that for any consistent estimators $\hat\Lambda_N$, $\hat\Lambda'_{N'}$, $\hat\Sigma_N$, $\hat\Sigma'_{N'}$ of $\Lambda$, $\Lambda'$, $\Sigma$, $\Sigma'$, respectively, under $H_0$,
%Suppose $N=N'$. For all $k=1,\ldots,K$ and positive integers $N$, let $q_k$, $q'_k$, $\hat q_{k,N}$ and $\hat q'_{k,N}$ denote the images of $q(\beta_k,\xi_k)$, $q'(\beta_k,\xi_k)$, $\hat q_N(\beta_k,\xi_k)$ and $\hat q'_{N}(\beta_k,\xi_k)$, respectively, under some global $C^\infty$ chart $\phi$. Then defining $\Psi_k$ and $D\Psi_k$ as in Section \ref{asymp}, let $\Lambda$ and $\Lambda'$ be the $Kn\times Kn$ block diagonal matrices whose $k$-th $n\times n$ diagonal blocks ($k=1,\ldots,K$) are $E[D\Psi_k(X,q_k)]$ and $E[D\Psi_k(X',q'_k)]$, respectively; also let $\Sigma$ and $\Sigma'$ be the covariance matrices of $[\Psi_1(X,q_1)^T,\ldots,\Psi_K(X,q_K)^T]^T$ and $[\Psi_1(X',q'_1)^T,\ldots,\Psi_K(X',q'_K)^T]^T$, respectively. Finally, denote $[\hat q_{1,N}^T,\ldots,\hat q_{K,N}^T]^T$ and $[(\hat q'_{1,N})^T,\ldots,(\hat q'_{K,N})^T]^T$ by $\hat q_{1:K,N}$ and $\hat q'_{1:K,N}$, respectively. Then by Theorem \ref{clt} and Slutsky's theorem, for any consistent estimators $\hat\Lambda_N$, $\hat\Lambda'_{N}$, $\hat\Sigma_N$ and $\hat\Sigma'_{N}$ of $\Lambda$, $\Lambda'$, $\Sigma$ and $\Sigma'$, respectively,
\begin{align*}
	&T_5((\beta_1,\xi_1),\ldots,(\beta_K,\xi_K)) \\
    =&N(\hat q_{N}-\hat q'_{N'})^T\Big(\hat\Lambda_N^{-1}\hat\Sigma_N\hat\Lambda_N^{-1}+\frac{1}{\eta}(\hat\Lambda'_{N'})^{-1}\hat\Sigma'_{N'}(\hat\Lambda'_{N'})^{-1}\Big)^{-1}(\hat q_{N}-\hat q'_{N'})\rightsquigarrow \chi_{Kn}^2.
\end{align*}
as $N\rightarrow\infty$. Denote the analogous test statistics for the sample means or medians, based on the covariance of the gradient and the average Hessian for the respective loss functions, by $T_3$ and $T_4$, respectively; under $H_0$, $T_3\rightsquigarrow\chi^2_n$ and $T_4\rightsquigarrow\chi^2_n$ as $N\rightarrow\infty$.
%Thus, we can reject $H_0$ at significance level $\alpha$ when $T_5>\chi_{Kn,1-\alpha}^2$, the $(1-\alpha)$-quantile of the $\chi^2$-distribution with $Kn$ degrees of freedom. As with permutation tests, we can find specific loci of differences by restricting the quantile indices to a region of interest in $[0,1)\times\partial M$. 
Then $T_5$ allows for finer comparisons between $F$ and $F'$ than $T_3$ and $T_4$.

\subsection{Other applications}\label{otherapplications}
Univariate quantiles and ranks are closely related. \citet{Chaudhuri1996} defined a notion of multivariate rank, which can be further generalised to a tangent vector-valued rank of $p\in M$ that describes the position of $p$ relative to the data cloud. These tangent vector-valued ranks could be used to construct tests of location.

While ordinary least squares regression estimates the conditional mean of the response variable given the predictors, this is only one aspect of the relationship between the response and predictors. By replacing the squared loss function with other loss functions, we can gain a more comprehensive view. The absolute loss function yields the conditional median, which is another conditional centre, whereas quantile loss functions return conditional quantiles. The foundational framework of univariate quantile regression dates back to \cite{Koenker1978}, while \cite{Chakraborty2003} developed multivariate geometric quantile regression by using the geometric quantile loss functions of \cite{Chaudhuri1996}. In the case of Hadamard space-valued responses, one can further generalise \cite{Chakraborty2003} by using our quantile loss functions, adapted to account for predictors, and investigate the entire conditional distribution to find deeper insights into conditional heteroscedasticity or tail behavior that cannot be captured by existing Fr\'echet regression methods; that is, methods for least squares regression with metric space-valued responses.

%However, this is only one aspect of the relationship between response and predictors. By replacing the squared-distance loss in metric spaces with other loss functions, we can gain a more comprehensive view of this relationship. The absolute distance loss yields the conditional median, another measure of central tendency, whereas in Hadamard spaces, quantile loss functions can provide a sense of the shape of the response distribution, along with other measures discussed in this section. In particular, if the variation in the response depends on the values of the predictors, or if the association between the predictors and the conditional mean response is weak, quantile regression could reveal useful patterns. The advantages of quantile regression for univariate Euclidean responses are well known, and \citet{Chakraborty2003} developed the methodology of multivariate geometric quantile regression.

	%All of the approaches mentioned in this section can also be used after performing the aforementioned transformation-retransformation procedure. 
	
	\section{Numerical experiments} \label{experiments}
	
The code for the experiments in this section, implemented in Python, can be found at \url{https://github.com/hayoungshin1/Quantiles-on-Hadamard-spaces}.

\subsection{Illustration of measures of distributional characteristics} \label{sexperiments}

%This section is very similar to Section 3.1 of \citet{ShinOh2024}, translated to a hyperbolic context. 
To explore the measures defined in Section \ref{measures}, we generated 16 data sets in the hyperbolic plane, four for each of dispersion, skewness, kurtosis, and spherical asymmetry. For a given characteristic, each of the data sets has an associated level $\nu\in\{0,1,2,3\}$, and the data sets are generated so that a plausible measure of that characteristic should decrease in sise as $\nu$ increases. For each data set, 2-vectors $w_j^*\ (j=1,\ldots,300)$, were generated to have geometric median $(0,0)^T$, ensuring that the median of the $X_j$'s is $p$, where $X_j=\exp_p((0,(w_j^*)^T)^T)$ and $p=(1,0,0)^T$.

First, we generated $N=300$ vectors $v_1,\ldots,v_N$ from a bivariate normal distribution $\mathcal{N}(0,I/4)$, and for each of the 16 data sets, we created a new vector $w_j$ from $v_j=(v_j^1,v_j^2)^T$, where for the dispersion data sets $w_j=(v_j^1, 4(2^{-\nu})v_j^2)^T$, for the skewness data sets $w_j=(1-\nu/3)((v_j^1)^2,(v_j^2)^2/2)^T+(\nu/3)(v_j^1, v_j^2/2)^T$, for the kurtosis data sets $w_j=(1-\nu/3)((v_j^1)^3,(v_j^2)^3/2)^T+(\nu/3)(v_j^1, v_j^2/2)^T$, and for the spherical asymmetry data sets $w_j=(1-\nu/3)((v_j^1)^3,(v_j^2)^3)^T+(\nu/3)(v_j^1, v_j^2)^T$. Defining $\hat{w}$ to be the Euclidean geometric median of $w_1,\ldots,w_N$, these were centred by letting $w_j^*=w_j-\hat{w}$ for each $j$, and 0 was concatenated to the beginning of each $w_j^*$; the resulting 3-component vectors can be considered elements of $T_p\mathbb{H}^2$, and the data point is given by $X_j:=\exp_p((0,(w_j^*)^T)^T)$. Subtracting $\hat{w}$ ensures that the Euclidean geometric median of the $w_j^*$'s is 0 and the median of the $X_j$'s is $p$ in every data set; this makes comparisons between contours for different data sets easier.

It is clear that the dispersion data sets become less dispersed as $\nu$ increases. The rationale for the skewness and kurtosis data sets is that the element-wise square of a non-skewed centred real data set is positively skewed, and the element-wise cube of a centred real data set has fatter tails than the original data set. Finally, cubing each coordinate for a spherically symmetric data set gives a data set that is no longer spherically symmetric but still has isotropic covariance. A common criticism of geometric quantiles is that the isoquantile contours do not follow the shape of the distribution when the covariance is not isotropic, especially for extreme quantiles; this was detailed in Section \ref{tr}. To show that this need not concern us, we have considered non-isotropic data. For additional exploration of the data sets, including visualisations and alternative measures of the characteristics, that show that these data sets serve our purposes with respect to $\nu$, see Section \ref{details} in the supplement.

%We are measuring the characteristics not of the generating distributions themselves but the empirical distributions based on the generated samples.

Recalling that the median of each data set is $p$, the values of $\xi$ used were those for which $\xi_p$ is in $\{(0,\cos(2\pi k/24),\sin(2\pi k/24))^T\}_{k=1,\ldots,24}\subset T_p\mathbb{H}^2$. Then, each $\sup$ over $S_p^1$ is approximated by the maximum over these 24 values, and each of the integrals (divided by the surface area) over $S_p^1$ by their mean. The results for each of these measures are given in Table \ref{noextable2}. For $\kappa_1$ and $\kappa_2$ to reasonably reflect tailedness, $\beta'$ should be quite large, while there is more freedom in $\beta$. We chose $\beta'=0.8$ and a somewhat small value of $\beta$ at 0.2. For the other measures, we used $\beta=0.5$.

\begin{table}[!h]
    \centering
    \caption{\small{Estimates for the seven measures for each $\nu\in\{0,1,2,3\}$.}}
    {\small
\begin{tabular}{ |c||c|c|c|c|c|c|c| } 
\hline
$\nu$ & $\delta_1(0.5)$ & $\delta_2(0.5)$ & $\gamma_1(0.5)$ & $\lVert\gamma_2(0.5)\rVert$ & $\kappa_1(0.2,0.8)$ & $\kappa_2(0.2,0.8)$ & $\alpha(0.5)$\\
\hline
0 & 2.785 & 2.267 & 0.392 & 0.160 & 29.219 & 29.070 & 0.663 \\ 
1 & 1.507 & 1.333 & 0.346 & 0.141 & 8.282 & 8.685 & 0.271 \\ 
2 & 0.894 & 0.870 & 0.155 & 0.068 & 5.970 & 6.351 & 0.155 \\ 
3 & 0.723 & 0.645 & 0.035 & 0.012 & 5.072 & 5.466 & 0.123 \\
\hline
\end{tabular}
    }
    \label{noextable2}
\end{table}

\begin{figure}
	\centering
		\includegraphics[width=0.3\linewidth]{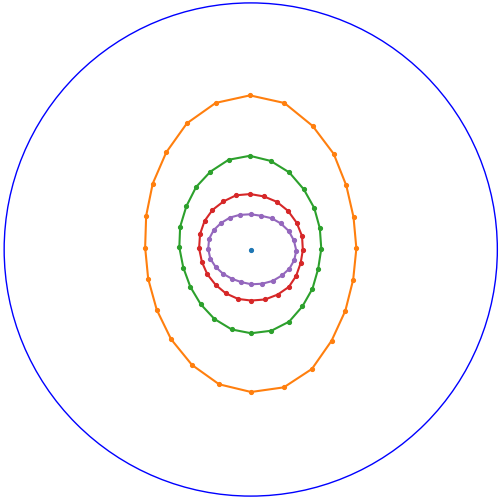}
        %\hspace{0.2\textwidth}
		\includegraphics[width=0.3\linewidth]{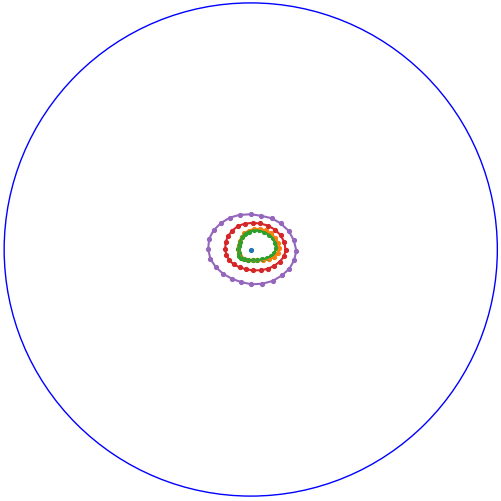}
		\vspace{-5mm}
 \caption{\small{Contours for the dispersion (left) and skewness (right) data sets visualised in the Poincar\'e disk. The orange, green, red, and purple contours represent the $\nu=0,1,2,3$ data sets; the blue point in the centre is their common median.}}
	\label{fig:simquantiles}
\end{figure}

Each of the measures decreases as $\nu$ increases. We performed the same experiments with $\beta$ or $\beta'$ set to $0.98$ to address concerns about extreme geometric quantiles. %which need not follow the shape of the distribution for large $\beta$ nor be contained in the convex hull of the distribution, even when the support is compact.
Even without a TR procedure, concerns about geometric quantiles do not seem to affect how our measures perform, and they should not. For example, if a data set is transformed to preserve its median but becomes less dispersed, then it is sufficient, though not necessary, for the contours to be pulled toward the median, and their specific shapes do not matter; see Figure \ref{fig:simquantiles}, which also shows reductions in skewness as $\nu$ increases. Section \ref{details} of the supplement includes visualisations of the 0.98-isoquantile contours, along with images for the other two characteristics.
\vspace{-1mm}

	\subsection{Outlier detection with single-cell RNA sequencing data} \label{rexp}
	
	\citet{Olsson2016} performed single-cell RNA sequencing, which is useful for `delineating hierarchical cellular states, including rare intermediates and the networks of regulatory genes that orchestrate cell-type specification', a central task in modern developmental biology. These hierarchical structures make the data ideal for hyperbolic embedding. The data, consisting of 8 cell types, have been embedded into the Poincar\'e disk by \citet{Chien2021} following the procedure in \citet{Klimovskaia2020}. We imagined a scenario in which our goal was to study one of these cell types, with 95 observations, but the data set had been corrupted by data from another cell type, consisting of 14 observations, which were treated as outliers.

 \begin{figure}[!t]
		\centering
	{\includegraphics[width=0.24\textwidth]{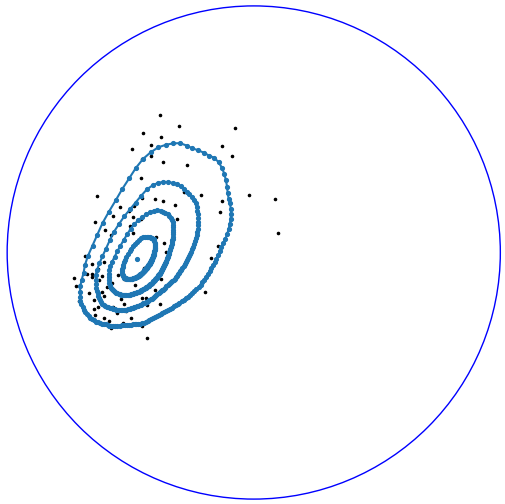}}
 {\includegraphics[width=0.24\textwidth]{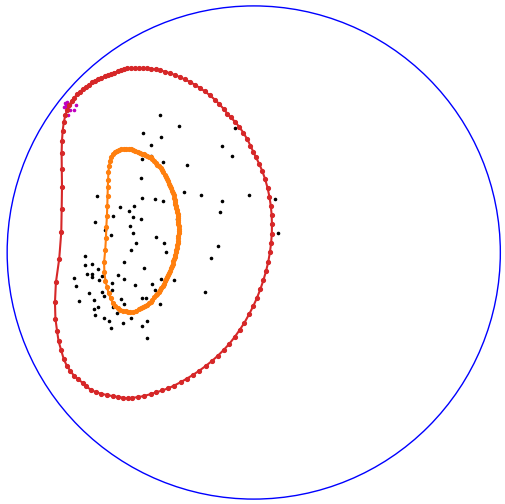}}
 {\includegraphics[width=0.24\textwidth]{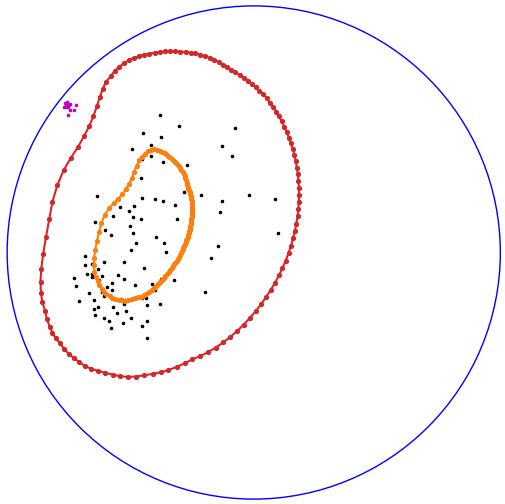}}
 {\includegraphics[width=0.24\textwidth]{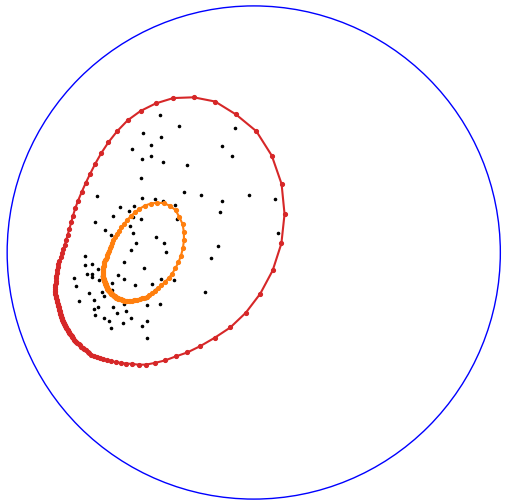}}
 \vspace{-5mm}
\caption{\small{The original data (black) and outliers (magenta) are embedded into the Poincar\'e disk. The TR median and isoquantile contours (blue) for the uncorrupted data are shown in the first subfigure. The next three subfigures show the TR 0.5-isoquantile (orange) contour and the outlier detection fence (red) for the following settings: a non-robust covariance estimate on the corrupt data; a robust covariance estimate on the corrupt data; a robust covariance estimate on the non-corrupt data.}}
		\label{fig:realdata}
	\end{figure}

To discern the shape of the non-corrupt data cloud, we used the procedure from Section \ref{tr}. We calculated the TR median ($\beta=0$) and the TR $\beta$-isoquantile contour for $\beta=0.2,0.4,0.6,0.8$, corresponding to univariate deciles, by using 120 different values of $\xi$. We then applied the outlier detection method from Section \ref{tr} to the corrupted data set. When the plug-in covariance estimate was used, some outliers and non-outliers were misidentified. However, with the robust minimum covariance determinant estimator of \cite{Rousseeuw1984}, every observation was correctly identified; additionally, in the original non-corrupt data set, all points were correctly designated as non-outliers. Results are shown in Figure \ref{fig:realdata}.

\subsection{Testing with single-cell gene expression data} \label{pexp}
\begin{figure}
	\centering
	{\includegraphics[width=0.24\textwidth]{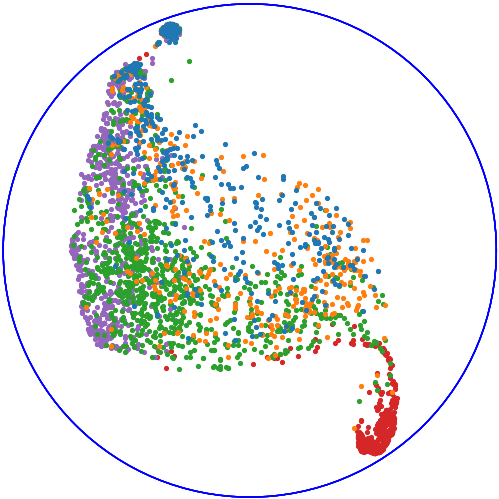}}
    %\hspace{0.03\textwidth}
	{\includegraphics[width=0.24\textwidth]{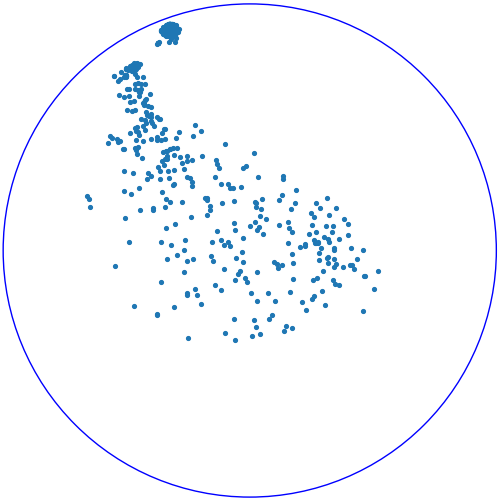}}
    %\hspace{0.03\textwidth}
    {\includegraphics[width=0.24\textwidth]{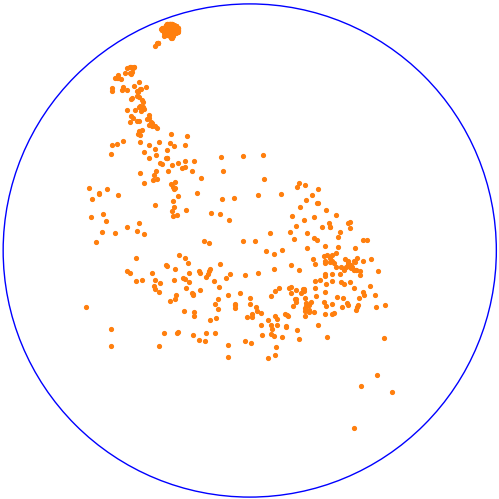}}
    {\includegraphics[width=0.24\textwidth]{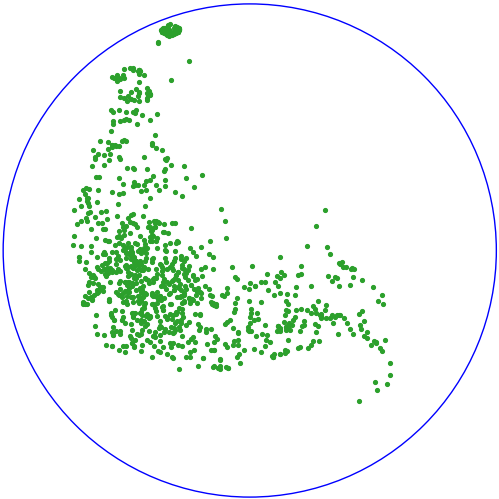}}
    \vspace{-5mm}
    \caption{\small{Left: Data coloured by stage. Centre-left: PS data. Centre-right: NP data. Right: HF data.}} %Key: blue=PS, orange=NP, green=HF, red=4SG, purple=4SFG-}
    \label{fig:moignard}
\end{figure}

\citet{Moignard2015} mapped blood development with single-cell gene expression analysis, which had `already been successfully applied to study the earliest stages of preimplantation mouse and human development, to identify lineage commitment and transcriptional regulatory events in blood, and ... to probe the emergence of hematopoietic stem cells from the hemogenic endothelium of the dorsal aorta.' Their data were embedded into the Poincar\'e disk by \citet{Mishne2023}. Each cell belongs to one of five stages of development; we focus on the three earliest, primitive streak (PS), neural plate (NP), and head fold (HF), from embryonic day 7.0, 7.5, and 7.75, respectively. As might be expected given their chronological proximity, these stages exhibit similarities (though not equality) in their distributions, as shown in Figure \ref{fig:moignard}. We compared the distributions for consecutive stages: PS vs. NP and NP vs. HF.

The data sets for these three stages are large, with 552, 624, and 1005 cells, so we considered the resulting empirical distributions as true population distributions. Our first goal was to estimate and compare the power of the tests using $T_0$, $T_1$, and $T_2$ from Section \ref{testing} by repeatedly sampling with replacement from both populations. Let $N=N'=120$, $l=1$ and $K=9$, and let $w_1=\cdots=w_9=1/9$, $\beta_1=0$, $\beta_2=\cdots=\beta_5=0.5$, and $\beta_6=\cdots=\beta_9=0.9$, and $\xi_k$ satisfy $(\xi_k)_{p}=(0,\cos(2\pi (k-2)/4),\sin(2\pi (k-2)/4))^T$ for $k=2,\ldots,9$ where $p=(1,0,0)^T$. We performed permutation tests with 1000 permutations, and Table \ref{proptable} was produced based on 1000 samplings for each of the two comparisons. In addition, identifying all points with their images under $\log_p$, we performed analogous permutation tests with $T_2$ calculated with the Euclidean geometric quantiles in $T_{p}\mathbb{H}^2$ using the same values of $\beta_k$ and $(\xi_k)_{p}$ ($k=1,\ldots,K$).

 \begin{table}[!h]
    \centering
    \caption{\small{Proportion of samplings (as percentages) in which the $p$-value for the given permutation test statistic exceeded the considered value, for both PS vs. NP and NP vs HF cases.}}
    {\small
\begin{tabular}{ |c||c|c|c|c|c||c|c|c|c|c| } 
\hline
& \multicolumn{5}{c||}{PS vs. NP} &  \multicolumn{5}{c|}{NP vs. HF} \\
\hline
Test statistic & 0.1 & 0.05 & 0.01 & 0.005 & 0.001 & 0.1 & 0.05 & 0.01 & 0.005 & 0.001 \\ \hline
 $T_0$ & 2.5 & 6.3 & 16.7 & 23.3 & 37.5 & 0 & 0.1 & 0.7 & 0.9 & 1.9 \\
$T_1$ & 6 & 10.4 & 27.7 & 35.9 & 53.4 & 0 & 0 & 0.1 & 0.1 & 0.4 \\ 
$T_2$ & 0.3 & 0.6 & 3.6 & 5.7 & 11.8 & 0 & 0 & 0 & 0 & 0.2 \\
Euclidean $T_2$ & 2.7 & 5.6 & 17.6 & 25.3 & 40.9 & 0 & 0 & 0.2 & 0.4 & 0.8 \\
\hline
\end{tabular}
    }
    \label{proptable}
\end{table}

Table \ref{proptable} shows that the hyperbolic version of $T_2$ showed much greater power than $T_0$, $T_1$, and the Euclidean version at all significance levels. For example, at the 5\% level in the PS/NP case, $T_0$, $T_1$, and the Euclidean $T_2$ led to Type II errors 6.3\%, 10.4\%, and 5.6\% of the time, respectively, while the hyperbolic $T_2$ did so only 0.6\% of the time. The contrast is less stark, but still clear, in the NP/HF case. The claim that the PS/NP distinction is more difficult than the NP/HF one is valid based on Figure \ref{fig:moignard}, as the former pair appears more similar than the latter pair.

 %\begin{table}[!h]
%    \centering
%    \caption{Means and quartiles across runs of the $p$-value for the considered test statistic.}
%    {\small
%\begin{tabular}{ |c||c|c|c|c| } 
%\hline
% & mean & 1st quartile & median & 3rd quartile \\ \hline
%$T_0$ & 0.0239 & 0.0002 & 0.0021 & 0.0126 \\ 
%$T_1$ & 0.0028 & 0.0000 & 0.0000 & 0.0005 \\
%\hline
%\end{tabular}
%    }
%    \label{meantable}
%\end{table}

     \begin{table}[!t]
    \centering
    \caption{\small{Mean $p$-values (as percentages) for permutation test statistics across the 1000 samplings, for both PS vs. NP and NP vs HF cases.}}
    {\small
\begin{tabular}{ |c||c|c| } 
\hline
Test statistic & PS vs. NP&  NP vs. HF \\
\hline
$T_0$ & 1.2834 & 0.0245 \\
$T_2((0,\xi_1,w_1),\ldots,(0.9,\xi_9,w_9))$ & 0.1890 & 0.0011 \\
Euclidean $T_2((0,(\xi_1)_p,w_1),\ldots,(0.9,(\xi_9)_p,w_9))$ & 1.1953 & 0.0085 \\ \hline
$T_1$ & 2.3827 & 0.0035 \\
$T_2((0.5,\xi_2,w_2))$ & 0.6074 & 0.0245 \\
$T_2((0.5,\xi_3,w_3))$ & 3.9903 & 0.0285 \\
$T_2((0.5,\xi_4,w_4))$ & 0.3247 & 0 \\
$T_2((0.5,\xi_5,w_5))$ & 0.1745 & 0.1181 \\
$T_2((0.9,\xi_6,w_6))$ & 0.0647 & 15.8351 \\
$T_2((0.9,\xi_7,w_7))$ & 20.2686 & 0.0362 \\
$T_2((0.9,\xi_8,w_8))$ & 0.2725 & 2.0032 \\
$T_2((0.9,\xi_9,w_9))$ & 0.1237 & 12.3611 \\ \hline
$T_2((0.5,\xi_2,w_2),\ldots,(0.5,\xi_5,w_5))$ & 0.2614 &0.0007  \\
$T_2((0.9,\xi_6,w_6),\ldots,(0.9,\xi_9,w_9))$ & 0.0548 & 0.0585 \\
 \hline
\end{tabular}
    }
    \label{meantable}
\end{table}

Mean $p$-values for $T_2((0,\xi_1,w_1),\ldots,(0.9,\xi_9,w_9))$ across the samplings were much lower than those for $T_0$, $T_1$, and the Euclidean $T_2$ in both columns of Table \ref{meantable}, further supporting the claim that quantiles are helpful in distinguishing distributions. Table \ref{meantable} also contains mean $p$-values for tests using only some of the quantile indices, in order to detect loci of distributional differences as described in Section \ref{testing}. For example, in the PS/NP case, $(0.9,\xi_6)$ and $(0.9,\xi_9)$ gave very small mean $p$-values, so we are particularly confident that the distributions differ at those quantiles. The mean $p$-value for $T_2((0.9,\xi_6,w_6),\ldots,(0.9,\xi_9,w_9))$ was smaller than those of $T_2((0.5,\xi_1,w_1),\ldots,(0.5,\xi_5,w_5))$ or $T_1=T_2((0,\xi_0))$, giving more support for the idea that the PS and NP distributions differ at their peripheries. Since $\xi_6$ and $\xi_9$ correspond to the rightward and downward directions, respectively, in Figure \ref{fig:moignard}, this aligns with our expectations based on the major difference between the PS and NP distributions—the scattering of points in the bottom-right of the third subfigure. 

Table \ref{meantable} shows that even when the distributions are difficult to distinguish at a single quantile index (such as $(\beta_7,\xi_7)$ for the PS/NP comparison or $(\beta_6,\xi_6)$ and $(\beta_9,\xi_9)$ in the NP/HF one), by using multiple indices, $T_2$ easily distinguishes the distributions if they are sufficiently different at other loci. This illustrates the idea from Section \ref{testing} of examining the distributions at a ``high resolution''.

We conducted similar comparisons using the asymptotic normality-based statistics $T_3$, $T_4$, and $T_5$. The results for the PS/NP case are presented in Tables \ref{proptable2} and \ref{meantable2}. The NP/HF case has been excluded because the $p$-values were so small that the comparison was uninteresting, with every entry in the equivalent of Table \ref{proptable2} equal to 0. As in the permutation tests, the ``high resolution'' statistic outperforms the ``low resolution'' ones and the Euclidean version. There is also broad agreement between Tables \ref{meantable} and \ref{meantable2} on the most significant loci of differences, indicated by the smallest mean $p$-values; for example, the indices and order of the 5 most significant loci of differences are the same in both tables. %both indicate that the greatest evidence for distributional difference exists at the $(0.5,\xi_6)$-quantile, followed by, in order, the $(0.9,\xi_9)$-, $(0.5,\xi_5)$-, $(0.9,\xi_8)$- and $(0.5,\xi_4)$-quantiles.
Though the two tables are not in complete accord, some disagreements are to be expected.
 \begin{table}[!h]
    \centering
    \caption{\small{Proportion of samplings (as percentages) in which the $p$-value for the given asymptotic normality-based test statistic exceeded the considered value, for the PS vs. NP case.}}
    {\small
\begin{tabular}{ |c||c|c|c|c|c| } 
\hline
Test statistic & 0.1 & 0.05 & 0.01 & 0.005 & 0.001 \\ \hline
$T_3$ & 3.3 & 6.5 & 20.8 & 28.8 & 49.4 \\
$T_4$ & 29.8 & 43.7 & 68.6 & 76.6 & 86.7 \\ 
$T_5$ & 0 & 0.1 & 0.7 & 1 & 1.8 \\
Euclidean $T_5$ & 0.8 & 1.3 & 2.3 & 2.8 & 4.3 \\
\hline
\end{tabular}
    }
    \label{proptable2}
\end{table}

     \begin{table}[!h]
    \centering
    \caption{\small{Mean $p$-values (as percentages) for asymptotic normality-based test statistics across the 1000 samplings, for the PS vs. NP case.}}
    {\small
\begin{tabular}{ |c||c| } 
\hline
$T_3$ & 1.3849 \\
$T_5((0,\xi_1),\ldots,(0.9,\xi_9))$ & 0.0243 \\
Euclidean $T_5((0,(\xi_1)_p),\ldots,(0.9,(\xi_9)_p))$ & 0.3304 \\ \hline
$T_4$ & 8.8929 \\
$T_5((0.5,\xi_2))$ & 1.6318 \\
$T_5((0.5,\xi_3))$ & 1.3457 \\
$T_5((0.5,\xi_4))$ & 0.7288 \\
$T_5((0.5,\xi_5))$ & 0.4443 \\
$T_5((0.9,\xi_6))$ & 0.1684 \\
$T_5((0.9,\xi_7))$ & 1.4683 \\
$T_5((0.9,\xi_8))$ & 0.4658 \\
$T_5((0.9,\xi_9))$ & 0.3974 \\ \hline
$T_5((0.5,\xi_2),\ldots,(0.5,\xi_5))$ & 0.4941  \\
$T_5((0.9,\xi_6),\ldots,(0.9,\xi_9))$ & 0.7604 \\
 \hline
\end{tabular}
    }
    \label{meantable2}
\end{table}

%The mean $p$-values for $T_0$, $T_1$ and $T_2$ were 0.0067, 0.0238 and 0.0019, respectively, and the proportion of runs in which the $p$-value for $T_1$ was greater than that for $T_2$ is 0.631, while for the reverse it was only 0.005 (these do not add to 1 because in the remaining cases the two $p$-values are equal, usually 0; the true $p$-values are unlikely to be equal, and running more permutations would give a more precise picture). For $T_0$ and $T_2$, the analogous figures were 0.288 and 0.073.  This and Table \ref{proptable} support the claim that quantiles are helpful tools in distinguishing distributions.
    
 %\begin{table}[!h]
%    \centering
%    \caption{Means and quartiles across runs of the $p$-value for the considered test statistic.}
%    {\small
%\begin{tabular}{ |c||c|c|c|c| } 
%\hline
% & mean & 1st quartile & median & 3rd quartile \\ \hline
%$T_0$ & 0.0239 & 0.0002 & 0.0021 & 0.0126 \\ 
%$T_1$ & 0.0028 & 0.0000 & 0.0000 & 0.0005 \\
%\hline
%\end{tabular}
%    }
%    \label{meantable}
%\end{table}

	\section{Conclusion} \label{conc}
	
	%In this paper, we generalise the notion of multivariate geometric quantiles to global non-positive curvature spaces with a bound at infinity. After showing some basic properties, notably scaled isometry equivariance and a necessary condition on Hadamard manifolds for the gradient of the quantile loss function at quantiles, we investigated the large-sample properties of sample quantiles on Hadamard manifolds, specifically a strong law of large numbers and joint asymptotic normality. We also provided an explicit expression for the gradient of the quantile loss function and other details necessary to perform the gradient descent algorithm on hyperbolic spaces. We defined some measures of dispersion, skewness, kurtosis and spherical asymmetry using quantiles, described how to do outlier detection using a transformation-retransformation procedure and defined a measure of distributional dissimilarity and described how to test the equality of distributions with it. Finally, we demonstrated these applications using simulated and real data sets.
	
	Considering the plethora of work on multivariate geometric quantiles, there are many avenues for fruitful research that generalise these studies to Hadamard spaces. Potential areas of research include geometric expectiles and the measure of distributional dissimilarity in Section \ref{testing}. Although this paper has focused on theory, more work is needed on the many potential applications and case studies with real-world data, such as diffusion tensors, to further demonstrate its practical utility. We also intend to devise faster algorithms for computation on hyperbolic and other Hadamard spaces, including SPD matrix space. Finally, efforts should be made to extend geometric quantiles to non-Hadamard spaces; questions to consider include whether our framework can be applied locally to regions of non-positive curvature, and examining how positive curvature might affect quantiles.

    In particular, we would like to further explore quantile regression, briefly introduced in Section \ref{otherapplications}. Regression methods, including local linear methods, could be developed for both Euclidean and non-Euclidean predictors, as has been done in the Fr\'echet regression case \citep{Petersen2019, Im2025}.

    There are also alternative definitions of multivariate quantiles besides that of \cite{Chaudhuri1996} which can also be generalised to Hadamard spaces, such as the depth-based notions of \cite{Hallin2010}, \cite{Kong2012}, and \cite{Chernozhukov2017}, the last of which also connects with optimal transport. For example, \cite{Hallin2010} indexed quantiles of $\mathbb{R}^n$-valued data $X_1,\ldots,X_N$ with a vector $u\in B^n(1)$ too, but instead defined the $u$-quantile as the $(n-1)$-dimensional hyperplane that minimises 
    \begin{align}\label{displacement}
    \sum_{i=1}^N \lvert D_i\rvert\{(1-\lVert u\rVert_2)I(D_i\leq 0)+\lVert u\rVert_2I(D_i>0)\},
    \end{align}
    where $D_i$ is the displacement in the direction $u/\lVert u\rVert_2$ of $X_i$ from the hyperplane ($i=1,\ldots,N$); that is, it is hyperplane that results from performing univariate $\lVert u\rVert_2$-quantile regression, treating $u/\lVert u\rVert_2$ as the vertical direction. This could be repurposed for Hadamard manifold-valued data $X_1,\ldots,X_N$ by defining, for $(\beta,\xi)\in[0,1)\times \partial M$, the $(\beta,\xi)$-quantile to be the minimiser, from among a certain class of submanifolds of codimension 1 that split the manifold in two, of (\ref{displacement}), where each $D_i$ is now defined as the displacement in the direction $\xi$ of $X_i$ from the submanifold. \cite{Hallin2010} then performed multivariate quantile regression by calculating their quantile hyperplanes in the Cartesian product of the predictor and response spaces; because the product of Hadamard spaces is also Hadamard, this notion of quantile regression can also be generalised to  Hadamard manifold-valued responses with Euclidean (or more generally Hadamard) predictors.

    %Existing multivariate quantile approaches predominantly rely on directional orientations or halfspace depth (e.g., \cite{Hallin2010}; \cite{Kong2012}), which, conceptually, align with Chaudhuri’s quantiles by using $L_1$-optimization. Recently, some studies have used optimal transport to define vector quantiles and Monge-Kantorovich depth \citep{Chakraborty2003, Carlier2016}.
	
	There is plenty of interesting work to do on Hadamard space statistics beyond quantiles. There has been very little work on persistent homology on Hadamard spaces. Because persistence diagrams are constructed from the spatial arrangement of points, persistent homology can yield very different results for data embedded in hyperbolic space, where distances increase more rapidly than in Euclidean space. Thus, we suspect that hierarchical data can be better understood through persistent homology when embedded in hyperbolic space.	
		
    \section*{Funding} 
\small{This research was supported by the National Research Foundation of Korea (NRF) funded by the Korean government (RS-2026-25468563), and by the Basic Science Research Program and G-LAMP Program through the NRF grant funded by the Ministry of Education (RS-2021-NR060140; RS-2025-25441317).}	

\section*{Data Availability} 
The data sets used in Section 7 can be downloaded at \url{https://github.com/hayoungshin1/Quantiles-on-Hadamard-spaces}.

	\begin{appendix}
\renewcommand{\thesection}{\Alph{section}}
\renewcommand{\theequation}{\thesection.\arabic{equation}}
\setcounter{equation}{0}

\renewcommand{\thefigure}{\thesection.\arabic{figure}}
\renewcommand{\thetable}{\thesection.\arabic{table}}

\makeatletter
\@addtoreset{figure}{section}
\@addtoreset{table}{section}
% supplement: every environment gets its own counter
\let\theorem\relax     \let\endtheorem\relax
\let\proposition\relax \let\endproposition\relax
\let\corollary\relax   \let\endcorollary\relax
\let\lemma\relax       \let\endlemma\relax
\let\definition\relax  \let\enddefinition\relax
\let\remark\relax      \let\endremark\relax
\let\assumption\relax  \let\endassumption\relax
\makeatother

\newcounter{sthm}[section]   \renewcommand{\thesthm}{\thesection.\arabic{sthm}}
\newcounter{sprop}[section]  \renewcommand{\thesprop}{\thesection.\arabic{sprop}}
\newcounter{scor}[section]   \renewcommand{\thescor}{\thesection.\arabic{scor}}
\newcounter{slem}[section]   \renewcommand{\theslem}{\thesection.\arabic{slem}}
\newcounter{sdef}[section]   \renewcommand{\thesdef}{\thesection.\arabic{sdef}}
\newcounter{srem}[section]   \renewcommand{\thesrem}{\thesection.\arabic{srem}}
\newcounter{sass}[section]   \renewcommand{\thesass}{\thesection.\arabic{sass}}

\theoremstyle{plain}
\newtheorem{theorem}[sthm]{Theorem}
\newtheorem{proposition}[sprop]{Proposition}
\newtheorem{corollary}[scor]{Corollary}
\newtheorem{lemma}[slem]{Lemma}
\newtheorem{definition}[sdef]{Definition}
\newtheorem{remark}[srem]{Remark}
\newtheorem{assumption}[sass]{Assumption}
\makeatother

		\section*{Supplement}

Throughout this document, labels of the form A.x, B.x, and C.x allude to entities in this supplement, while those of the form x allude to those in the main manuscript.

\section{Geometric preliminaries} \label{geomsupp}

This section is largely based on \citet{Bridson1999}, particularly chapters I.1., II.1., II.2, II.3, II.4, and II.8, with slight modifications. Throughout this section, $(M,d)$, usually shortened to $M$, will represent a metric space. Also, we will refer to all metrics as $d$ unless disambiguation is necessary. An isometry $f:M_1\rightarrow M_2$ is a bijection between two metric spaces such that $d(f(x),f(y))=d(x,y)$ for all $x,y\in M_1$. Given an interval $I\subset \mathbb{R}$, a geodesic is a map $\gamma:I\rightarrow M$ for which there exists some constant $u_\gamma\geq0$ (called the speed of the geodesic) and, for all $t\in I$, some neighbourhood $J_t\subset I$ such that $t_1,t_2\in J_t$ implies $d(\gamma(t_1),\gamma(t_2))=u_\gamma\lvert t_2-t_1\rvert$; a minimal geodesic is a geodesic for which the aforementioned $J_t$ can be all of $I$. If $0\in I$, we say that $\gamma$ issues from $\gamma(0)$; if $I=[0,l]$, that $\gamma$ is a geodesic from $\gamma(0)$ to $\gamma(l)$; and if $I=[0,\infty)$, that $\gamma$ is a geodesic ray. A geodesic segment is the image of a geodesic. In $\mathbb{R}^n$, (minimal) geodesics are lines, line segments, and rays, while in $S^n$, geodesic segments are arcs of great circles, but only those arcs of length less than or equal to $\pi$ are minimal geodesic segments. 
	
	%Note that our definition of geodesics slightly differs from that of \citet{Bridson1999}; they would call our geodesics linearly reparameterised local geodesics, while we would call theirs unit-speed minimal geodesics. We have adjusted the definitions so that on connected Riemannian manifolds, which are metric spaces, our geodesics are equivalent to geodesics in the Riemannian sense and $u_\gamma=\lVert\dot\gamma(0)\rVert$, where $\lVert\cdot\rVert$ denotes the norm according to the Riemannian metric. In the rest of this paper (excluding Sections \ref{asymp} and \ref{hyp} and their associated proofs in the appendices), the term `geodesic' will, as in \citet{Bridson1999}, refer to a unit-speed minimal geodesic unless otherwise stated.

 Note that our definition of geodesics slightly differs from that of \citet{Bridson1999}; we would call their geodesics unit-speed minimal geodesics. In the rest of this section, the term `geodesic' will refer to a unit-speed minimal geodesic unless otherwise stated. However, in the rest of the paper, it will be more convenient to use our definition as it coincides with the differential geometric definition of a geodesic on connected Riemannian manifolds, which are metric spaces, with $u_\gamma=\lVert\dot\gamma(0)\rVert$, where $\lVert v\rVert$ denotes the Riemannian norm of $v\in T_pM$ according to the Riemannian metric in $T_pM$, the tangent space at $p\in M$. Thankfully, the shift between the two definitions is not difficult, as on Hadamard spaces they turn out to be equivalent modulo speed.
	
	$M$ is called a geodesic space if there exists a geodesic from $x$ to $y$ for all $x,y\in M$, and uniquely geodesic if such a geodesic is unique. Recall that a proper metric space is one whose subsets are compact if and only if they are closed and bounded, or equivalently, all closed balls are compact. The Hopf-Rinow theorem shows that if $M$ is a complete and connected Riemannian manifold, it is both proper and geodesic.

	A notion of an angle between two geodesics can be defined using so-called comparison triangles. Given $p,q,r\in M$, a comparison triangle $\overline{\Delta}(p,q,r)$ is a triangle in $\mathbb{R}^2$, the Euclidean plane with its standard metric, with vertices $\overline{p}$, $\overline{q}$ and $\overline{r}$ such that $d(p,q)=d(\overline{p},\overline{q})$, $d(q,r)=d(\overline{q},\overline{r})$ and $d(p,r)=d(\overline{p},\overline{r})$. By the triangle inequality, such a triangle always exists and is unique up to isometry. Denote the interior angle of $\overline{\Delta}(p,q,r)$ at $\overline{p}$ by $\overline{\angle}_p(q,r)$. Then, given geodesics $\gamma_1$ and $\gamma_2$ in $M$ with $\gamma_1(0)=\gamma_2(0)$, the Alexandrov angle between $\gamma_1$ and $\gamma_2$ is 
	\begin{align*}
	\angle(\gamma_1,\gamma_2):=\lim\sup_{t_1,t_2\rightarrow 0}\overline{\angle}_{\gamma_1(0)}(\gamma_1(t_1),\gamma_2(t_2))=\lim_{\epsilon\rightarrow 0}\sup_{t_1,t_2\in(0,\epsilon)}\overline{\angle}_{\gamma_1(0)}(\gamma_1(t_1),\gamma_2(t_2)).
	\end{align*}
	The angle between two geodesic segments with a common endpoint is then defined as the angle between geodesics issuing from the common point whose images are those segments.
	
	A geodesic triangle $\Delta\subset M$ consists of three points $p,q,r\in M$ called vertices and three geodesic segments called sides joining each pair of vertices. Denote these sides by $[p,q],[q,r]$ and $[p,r]$, and the corresponding sides of a comparison triangle $\overline{\Delta}$ by $[\overline{p},\overline{q}],[\overline{q},\overline{r}]$ and $[\overline{p},\overline{r}]$. A point $\overline{x}\in[\overline{p},\overline{q}]$ is called a comparison point for $x\in[p,q]$ if $d(p,x)=d(\overline{p},\overline{x})$, and comparison points on the other sides are defined similarly. Then, $\Delta$ is said to satisfy the CAT(0) inequality if for all $x,y\in\Delta$ and corresponding comparison points $\overline{x},\overline{y}\in\overline{\Delta}$,  $d(x,y)\leq d(\overline{x},\overline{y})$. If $M$ is a geodesic space in which every geodesic triangle satisfies the CAT(0) inequality, $M$ is called a CAT(0) space. Intuitively, triangles in a CAT(0) space are at least as thin as their Euclidean comparison triangles. A metric space $M$ is said to be of non-positive curvature if it is locally CAT(0); that is, for each $x\in M$, there exists some $r_x>0$ such that the ball $B(x,r_x)$ is a CAT(0) space under the induced metric. CAT(0) spaces are uniquely geodesic spaces and are simply connected. A metric generalisation of the Cartan-Hadamard theorem shows that a complete, simply connected geodesic space that has non-positive curvature (i.e., locally CAT(0)) is a CAT(0) space (i.e., globally CAT(0)). A complete CAT(0) space is called a Hadamard space or, alternatively, a global non-positive curvature space. Alternatively, a Hadamard space is defined as a nonempty, complete metric space $M$ in which, for all $p_1,p_2\in M$, there exists some $m_{p_1,p_2}\in M$ for which $d(m_{p_1,p_2},p_3)^2+d(p_1,p_2)^2/4\leq (d(p_1,p_3)^2+d(p_2,p_3)^2)/2$ for all $p_3\in M$.
	
    Recall that on a Riemannian manifold $M$, the Riemannian gradient $\mathrm{grad }f(p)$ at $p\in M$ of a differentiable real-valued function $f:M\rightarrow\mathbb{R}$ is defined by $\langle\mathrm{grad }f(p), v\rangle=df_p(v)$ for all $v$ in the tangent space $T_pM$, where for $v_1,v_2\in T_pM$, $\langle v_1,v_2\rangle$ denotes the Riemannian inner product of $v_1$ and $v_2$ according to the Riemannian metric in $T_pM$, and $df_p$ is the differential of $f$ at $p$. 
 
 Importantly, a $C^\infty$ connected Riemannian manifold is a metric space of non-positive curvature in the above sense if and only if its sectional curvatures are non-positive, and the Riemannian angle between two geodesics is equal to the Alexandrov angle. This and the aforementioned generalisation of the Cartan-Hadamard theorem imply that complete, simply connected Riemannian manifolds of non-positive sectional curvature are Hadamard spaces, and they are called Hadamard manifolds. Because CAT(0) spaces can be shown to be simply connected, Hadamard manifolds can alternatively be defined as Riemannian manifolds that are Hadamard spaces. Examples include hyperbolic spaces and the spaces of symmetric and positive definite matrices. 
 
 If $M$ is a Hadamard manifold and $p\in M$, the exponential map at $p$, $\exp_p:T_pM\rightarrow M$, is a diffeomorphism defined on all of the tangent space $T_pM$ by the Cartan-Hadamard theorem. It is defined by $\exp_p(v)=p$ if $v=0$ and $\gamma_1(\lVert v\rVert)$ otherwise, where $\gamma_1:\mathbb{R}\rightarrow M$ is the unique geodesic satisfying $\gamma_1(0)=p$, $\gamma_1'(0)=v/\lVert v\rVert$. Then, the inverse exponential map $\log_p:M\rightarrow T_pM$ is also a diffeomorphism and is defined by $\log_p(q)=\lVert v\rVert\gamma_2'(0)$, where $\gamma_2:\mathbb{R}\rightarrow M$ is the unique geodesic between $p$ and $q$ issuing from $p$. For example, when $M=\mathbb{R}^n$, $\exp_p(v)=p+v$ and $\log_p(q)=q-p$.

    A hyperbolic space is a Hadamard manifold of constant negative sectional curvature. There is a scaled isometry between hyperbolic spaces of the same dimension and different curvatures, so without loss of generality, we mainly deal with hyperbolic spaces of sectional curvature -1. As discussed in the introduction, because hierarchical data can be effectively embedded into hyperbolic spaces, they have become popular across many fields.
    
    Unlike the sphere, it is impossible to embed hyperbolic space in Euclidean space without distortions, so it is harder to visualise and understand intuitively. There are, however, many popular ways to embed $n$-dimensional hyperbolic space in Euclidean space with distortions. We will use two of the most famous. The first is the Poincar\'e ball (also called the Poincar\'e disk when $n=2$),
	$
	\mathbb{B}^n=\{(q^1,\ldots,q^n)^T\in\mathbb{R}^n|(q^1)^2+\cdots+(q^n)^2<1\}.
	$
    In this embedding, the entire space is contained in the interior of an $n$-dimensional unit ball, while distances approach infinity as one approaches the boundary of the ball, which is identified with the boundary at infinity $\partial\mathbb{H}^n$. When $n=2$, the Poincar\'e ball (or disk) is especially useful for visualisation purposes. The other embedding is the hyperboloid,
	$
	\mathbb{H}^n=\{(p^1,\ldots,p^{n+1})^T\in\mathbb{R}^{n+1}|\langle p,p\rangle_{\textbf{M}}=-1,\ p^1>0\},
	$
	where $\langle\cdot,\cdot\rangle_{\textbf{M}}$ is the Minkowski (pseudo-)inner product defined by $\langle p,q\rangle_{\textbf{M}}=-p^1q^1+\sum_{j=2}^{n+1}p^jq^j$ for $p=(p^1,\ldots,p^{n+1})^T,~q=(q^1,\ldots,q^{n+1})^T\in\mathbb{R}^{n+1}$. When endowed with appropriate metrics, these two embeddings are isometric via the function $h:\mathbb{H}^n\rightarrow \mathbb{B}^n$ defined by $h((p^1,p^2,\ldots,p^{n+1})^T)=(p^2,\ldots,p^{n+1})^T/(p^1+1)$.
		
	The hyperboloid is more convenient for doing calculations. In this embedding, the tangent space at $p\in\mathbb{H}^n$ is characterised by $T_p\mathbb{H}^n=\{v\in\mathbb{R}^{n+1}|\langle p,v\rangle_{\textbf{M}}=0\}$. Even though $\langle \cdot,\cdot\rangle_{\textbf{M}}$ is not positive definite, its restriction to $T_p\mathbb{H}^n$ is, and we can define a norm on the tangent space by $\lVert v\rVert_{\textbf{M}}=\sqrt{\langle v,v\rangle_{\textbf{M}}}$ for $v\in T_p\mathbb{H}^n$. The exponential map is then given by $\exp_p(v)=\cosh(\lVert v\rVert_{\textbf{M}})p+\sinh(\lVert v\rVert_{\textbf{M}})v/\lVert v\rVert_{\textbf{M}}$. For $p,q\in \mathbb{H}^n$, the inverse exponential map is given by $\log_{p}(q)=\cosh^{-1}( -\langle p,q\rangle_{\textbf{M}})(q+\langle p,q\rangle_{\textbf{M}} p)/\lVert q+\langle p,q\rangle_{\textbf{M}} p\rVert_{\textbf{M}}$, meaning $d(p,q)=\cosh^{-1}(-\langle p,q\rangle_{\textbf{M}})$.

\section{Experimental details} \label{details}

\begin{figure}[!t]
	\centering
		\includegraphics[width=0.20\linewidth]{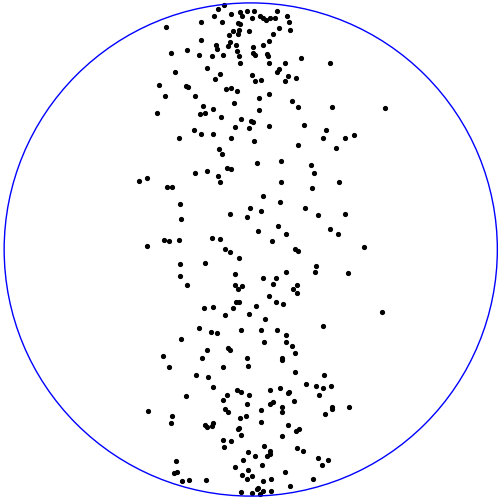}
		\includegraphics[width=0.20\linewidth]{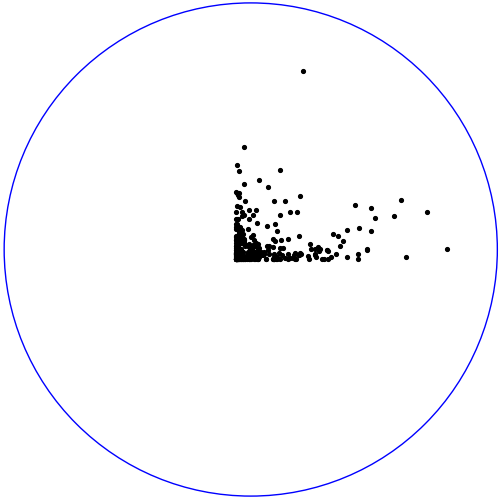}
		\includegraphics[width=0.20\linewidth]{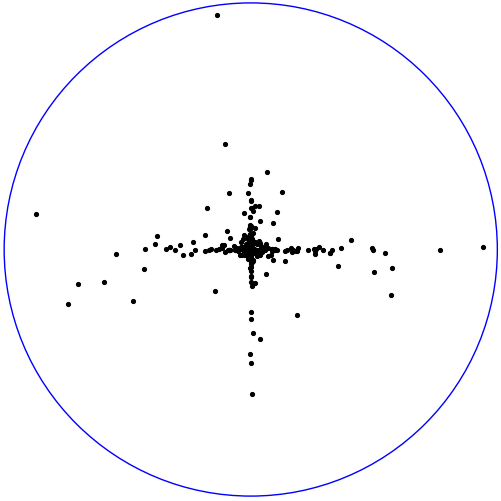}
		\includegraphics[width=0.20\linewidth]{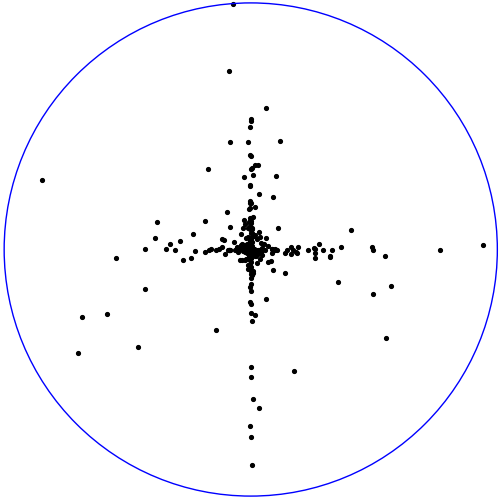}
		\includegraphics[width=0.20\linewidth]{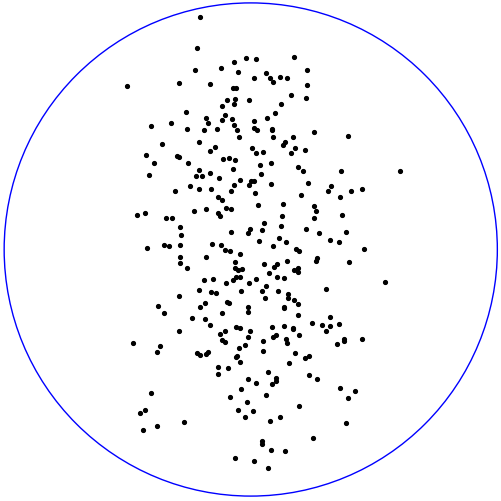}
		\includegraphics[width=0.20\linewidth]{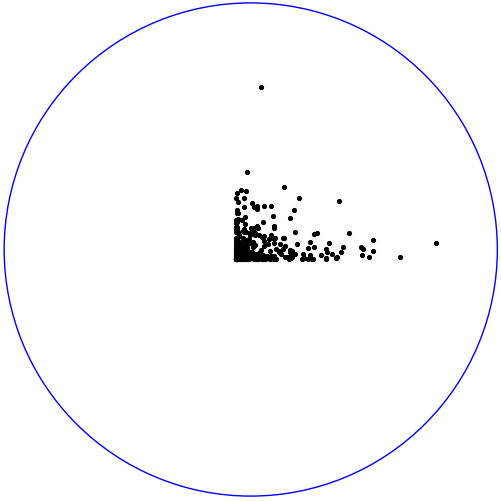}
		\includegraphics[width=0.20\linewidth]{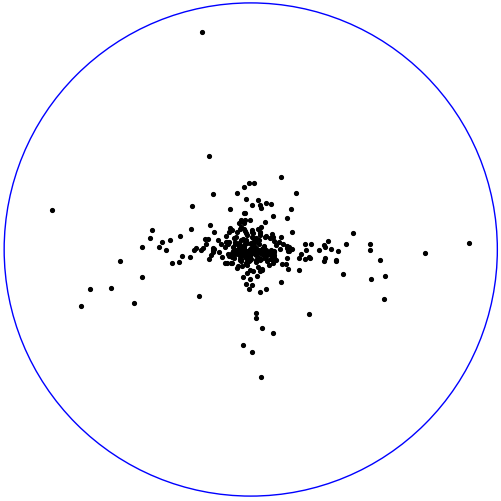}
		\includegraphics[width=0.20\linewidth]{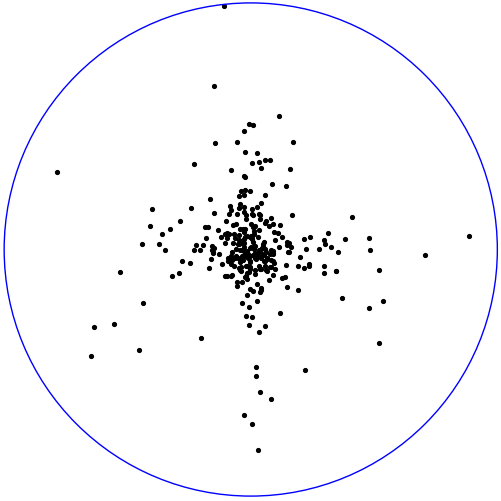}
		\includegraphics[width=0.20\linewidth]{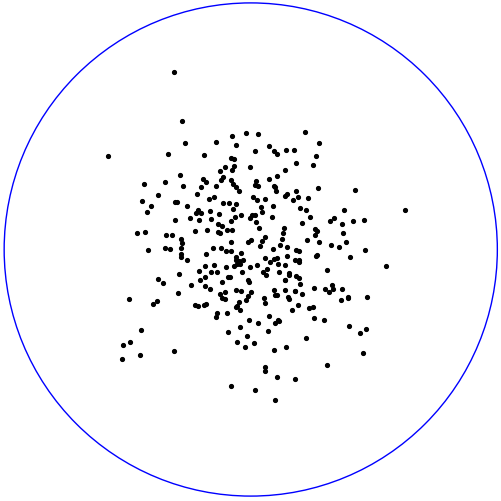}
		\includegraphics[width=0.20\linewidth]{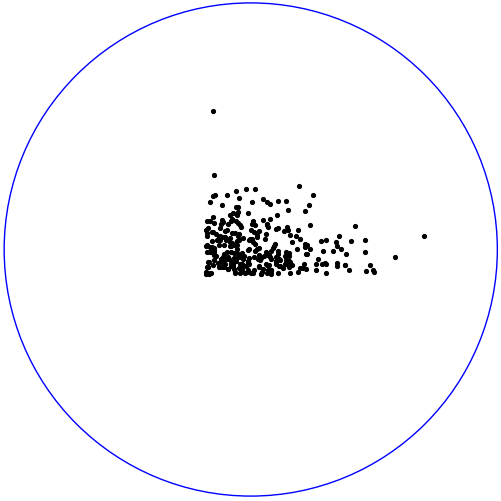}
		\includegraphics[width=0.20\linewidth]{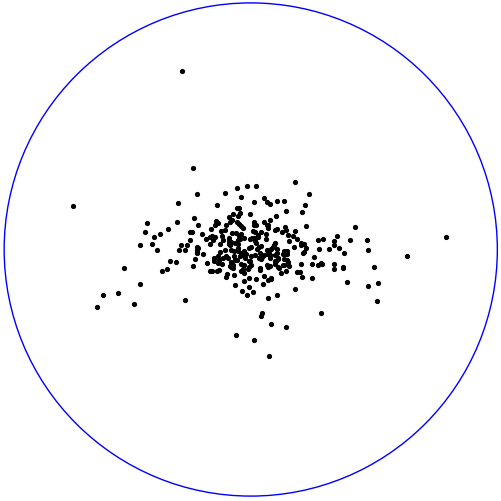}
		\includegraphics[width=0.20\linewidth]{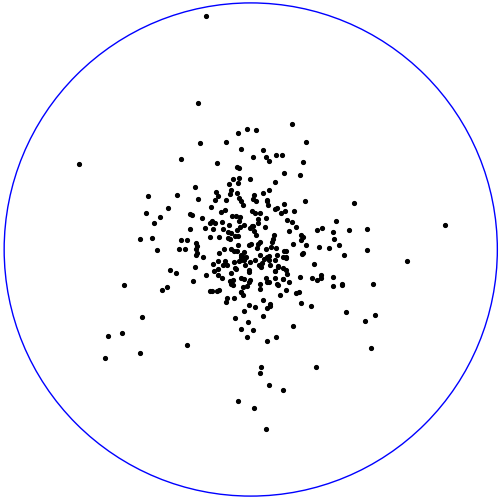}
		\includegraphics[width=0.20\linewidth]{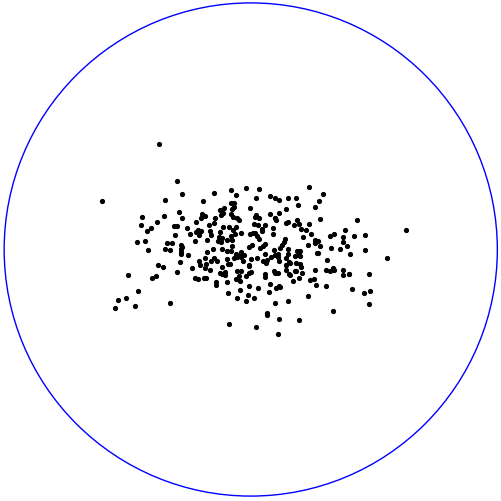}
		\includegraphics[width=0.20\linewidth]{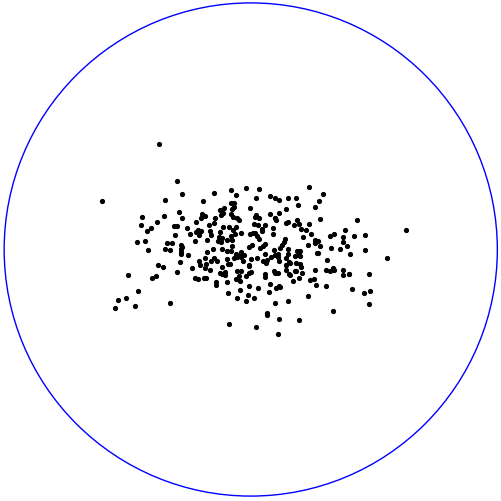}
		\includegraphics[width=0.20\linewidth]{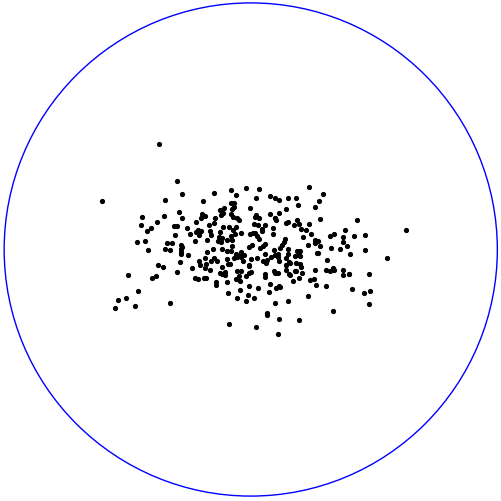}
		\includegraphics[width=0.20\linewidth]{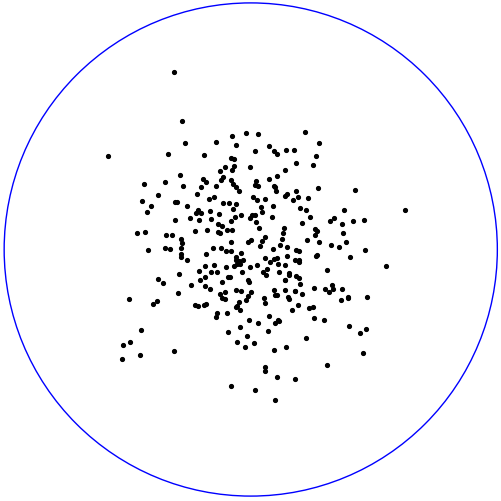}
 \caption{The columns display, in order, the dispersion, skewness, kurtosis, and spherical asymmetry data sets visualised in the Poincar\'e disk. The $i$th row, for $i\in\{1,2,3,4\}$ displays the data sets for $\nu=i-1$.}
	\label{fig:distributions}
\end{figure}

Figure \ref{fig:distributions} shows the 16 data sets embedded in the Poincar\'e disk. Table \ref{noextable1} gives $\delta_0(\beta)$, $\gamma_0(\beta)$ and $\kappa_0(\beta,\beta')$, defined in Section \ref{measures}, for each coordinate of $w_j^*$ from the respective data sets. These are imperfect as measures on hyperbolic space because they do not take curvature into account and depend on the coordinate system, but they do give a rough quantile-based feel. For the kurtosis measures $(\beta,\beta')=(0.2,0.8)$, and $\beta=0.5$ otherwise.

\begin{table}[!h]
    \centering
    \caption{Univariate measures of dispersion, skewness, and kurtosis for each of the two coordinates of the $w_j^*$'s from the relevant data sets.}
    {\small
\begin{tabular}{ |c||c|c|c| } 
\hline
$\nu$ & $\delta_0(0.5)$ & $\gamma_0(0.5)$ & $\kappa_0(0.2,0.8)$ \\
\hline
0 & 0.606, 2.845 & 0.485, 0.376 & 105.140, 79.614 \\ 
1 & 0.606, 1.422 & 0.360, 0.409 & 8.562, 7.933 \\ 
2 & 0.606, 0.711 & 0.191, 0.214 & 5.886, 5.467 \\ 
3 & 0.606, 0.356 & 0.045, 0.043 & 4.961, 4.607 \\
\hline
\end{tabular}
    }
    \label{noextable1}
\end{table}

We also performed the same experiments as in the main manuscript, with $\beta$ or $\beta'$ set to 0.98. The results are in Table \ref{extable2}, analogous to Table \ref{noextable2}, while Table \ref{extable1} is analogous to Table \ref{noextable1}. Our measures work just as well for extreme values of $\beta$, despite concerns about extreme geometric quantiles which need not follow the shape of the distribution for large $\beta$ nor be contained in the convex hull of the distribution, even when the support is compact.

\begin{table}[!h]
    \centering
    \caption{Univariate measures of dispersion, skewness, and kurtosis for each of the two coordinates of the $w_j^*$'s from the relevant data sets for extreme values of $\beta$ or $\beta'$.}
    {\small
\begin{tabular}{ |c||c|c|c| } 
\hline
$\nu$ & $\delta_0(0.98)$ & $\gamma_0(0.98)$ & $\kappa_0(0.2,0.98)$ \\
\hline
0 & 2.207, 8.828 & 0.876, 0.839 & 570.687, 427.777 \\ 
1 & 2.207, 4.414 & 0.830, 0.806 & 28.935, 26.713 \\ 
2 & 2.207, 2.207 & 0.525, 0.501 & 13.923, 12.916 \\ 
3 & 2.207, 1.104 & -0.029, -0.068 & 8.732, 8.103 \\
\hline
\end{tabular}
    }
    \label{extable1}
\end{table}

\begin{table}[!h]
    \centering
    \caption{Estimates for the seven measures for each $\nu\in\{0,1,2,3\}$ for extreme values of $\beta$ or $\beta'$.}
    {\small
\begin{tabular}{ |c||c|c|c|c|c|c|c| } 
\hline
$\nu$ & $\delta_1(0.98)$ & $\delta_2(0.98)$ & $\gamma_1(0.98)$ & $\lVert\gamma_2(0.98)\rVert$ & $\kappa_1(0.2,0.98)$ & $\kappa_2(0.2,0.98)$ & $\alpha(0.98)$\\
\hline
0 & 7.768 & 6.466 & 0.342 & 0.096 & 140.534 & 143.326 & 0.205 \\ 
1 & 5.510 & 5.407 & 0.310 & 0.086 & 30.769 & 33.394 & 0.098 \\ 
2 & 4.399 & 4.322 & 0.182 & 0.047 & 18.781 & 20.196 & 0.082 \\ 
3 & 4.055 & 3.709 & 0.019 & 0.002 & 14.305 & 15.070 & 0.051 \\
\hline
\end{tabular}
    }
    \label{extable2}
\end{table}

\begin{figure}[!t]
	\centering
		\includegraphics[width=0.23\linewidth]{dispersionquantiles0.5.png}
		\includegraphics[width=0.23\linewidth]{skewnessquantiles0.5.png}
		\includegraphics[width=0.23\linewidth]{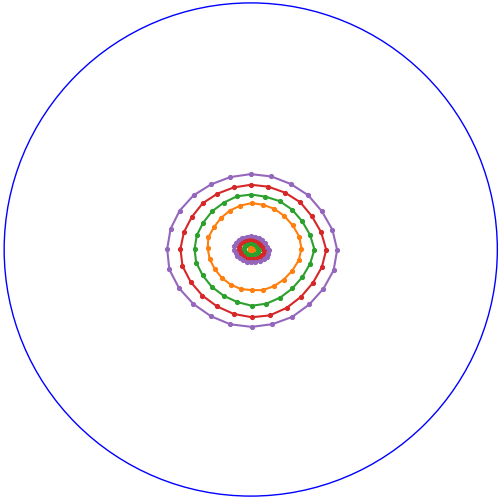}
		\includegraphics[width=0.23\linewidth]{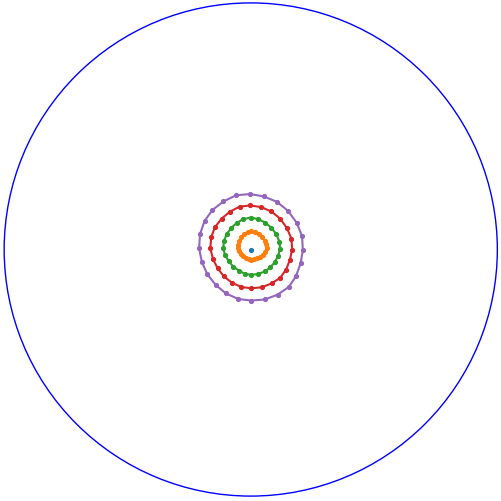}
		\includegraphics[width=0.23\linewidth]{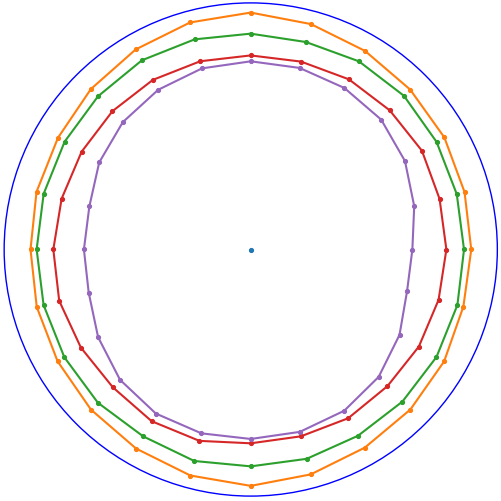}
		\includegraphics[width=0.23\linewidth]{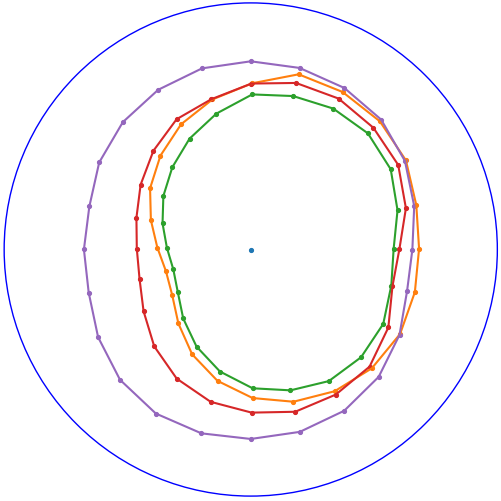}
		\includegraphics[width=0.23\linewidth]{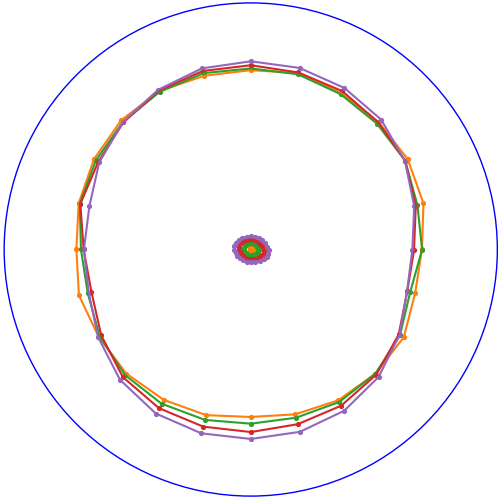}
		\includegraphics[width=0.23\linewidth]{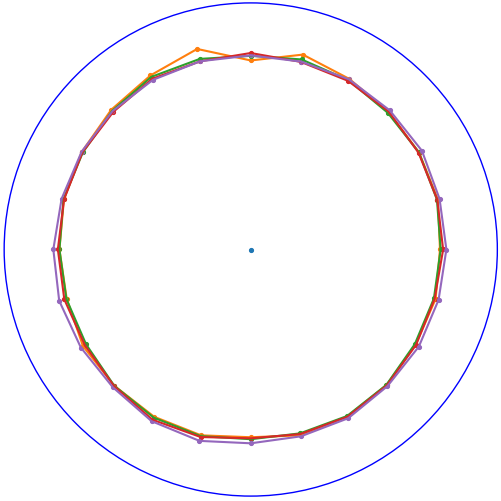}
 \caption{The four columns display isoquantile contours for, in order, the dispersion, skewness, kurtosis, and spherical asymmetry data sets visualised in the Poincar\'e disk. In the first row, $(\beta,\beta')=(0.2,0.8)$ for the kurtosis measures and $\beta=0.5$ otherwise; in the second row, $(\beta,\beta')=(0.2,0.98)$ for the kurtosis measures and $\beta=0.98$ otherwise. The orange, green, red, and purple contours are the contours for the $\nu=0,1,2,3$ data sets, respectively, while the blue point in the centre of the disk is their common median.}
	\label{fig:simquantilessupp}
\end{figure}

Figure \ref{fig:simquantilessupp} is an expansion on Figure \ref{fig:simquantiles}, containing the contours for all four characteristics and both extreme and non-extreme values of $\beta$ (or $\beta'$). The reductions in dispersion and skewness are clearly visible in the contours in the first two columns; the other two columns are harder to decipher but included for completeness.

%We can say that these (empirical) distributions serve our purposes with respect to $\nu$ thanks to the above rationale, and Figure \ref{fig:distributions} and Tables \ref{noextable1} and \ref{extable1} in the main paper.

\section{Proofs} \label{proofs}

		\subsection{Proof of Propositions~\ref{measurable} and \ref{equivariance}}\label{proof_meas}

  \begin{proof}[Proof of Proposition \ref{measurable}]
		If $x\neq p$ , $\angle_p(x,\xi)$ is continuous as a function of $x$ by Proposition II.9.2(1) of \citet{Bridson1999}, while $d(p,x)$ is continuous as a function of $x$ by the triangle inequality. $\rho(x,p;\beta,\xi)$ is, therefore, continuous as a function of $x$ on $M\backslash \{p\}$. Since $(1-\beta)d(p,x)\leq\rho(x,p;\beta,\xi)\leq(1+\beta)d(p,x)$ and the limits of both sides as $x$ approaches $p$ are 0, the map is also continuous at $x=p$. So, $\rho(x,p;\beta,\xi)$ is non-negative and continuous as a function of $x$ for all $x\in M$. 
	\end{proof}

 \begin{proof}[Proof of Proposition \ref{equivariance}]
It is clear from definitions that for any unit-speed geodesic $\gamma$ into $M$, $g\gamma$ is a unit-speed geodesic into $M'$ issuing from $g(\gamma(0))$, and that $\gamma_1$ and $\gamma_2$ are asymptotic geodesic rays in $M$ if and only if $g\gamma_1$ and $g\gamma_2$ are asymptotic geodesic rays in $M'$. So, $g\xi$ is an element of $\partial M'$. Scaled isometries also preserve Alexandrov angles between geodesics; that is, given geodesics $\gamma_3$ and $\gamma_4$ into $M$ issuing from the same point, $\angle(\gamma_3,\gamma_4)$=$\angle(g\gamma_3,g\gamma_4)$. This is because the Euclidean comparison triangle $\overline{\Delta}(g\gamma_3(0),g\gamma_3(\sigma t_3),g\gamma_4(\sigma t_4))$, where $t_3$ and $t_4$ are in the domains of $\gamma_3$ and $\gamma_4$, respectively, is similar to the comparison triangle $\overline{\Delta}(\gamma_3(0),\gamma_3(t_3),\gamma_4(t_4))$ scaled by a factor of $\sigma$. So, $\overline{\angle}_{g\gamma_3(0)}(g\gamma_3(\sigma t_3),g\gamma_4(\sigma t_4))=\overline{\angle}_{\gamma_3(0)}(\gamma_3(t_3),\gamma_4(t_4))$. By the definition of the Alexandrov angle, this implies $\angle(\gamma_3,\gamma_4)$=$\angle(g\gamma_3,g\gamma_4)$.
 
		All of this means that $\angle_{g(p)}(g(x),g\xi)=\angle_p(x,\xi)$. Therefore,
		\begin{eqnarray*}
			\rho(g(x),g(p);\beta,g\xi)&=&d'(g(p),g(x))+\beta d'(g(p),g(x))\mathrm{cos}(\angle_{g(p)}(g(x),g\xi)) \\
			&=&\sigma d(p,x)+\beta \sigma d(p,x)\mathrm{cos}(\angle_p(x,\xi))  \\
			&=&\sigma[d(p,x)+\beta d(p,x)\mathrm{cos}(\angle_p(x,\xi))] \\
			&=&\sigma\rho(x,p;\beta,\xi),
		\end{eqnarray*}
		from which the conclusion follows since $g$ is bijective.
	\end{proof}
  
		\subsection{A crucial proposition}\label{proof_prop3.3}
		The following result is crucial for the proofs in this paper. Proposition \ref{abc}(a) characterises $\xi_p$ asymptotically as the unit vector in $T_pM$ facing the ray $\gamma\in\xi$ as $t\rightarrow\infty$, which is very intuitive.
	
	\begin{proposition} \label{abc}
		Let $M$ be an $n$-dimensional Hadamard manifold and $TM$ its tangent bundle. (a) $\lim_{t\rightarrow\infty}\log_p(\gamma(t))/\lVert\log_p(\gamma(t))\rVert=\xi_p$, where $p\in M$ and $\gamma:[0,\infty)\rightarrow M$ is a member of the equivalence class $\xi$. (b) $\exp:TM\rightarrow M$ defined by $\exp(p,v)=\exp_p(v)$, and $(p,x)\mapsto (p,\log_p(x))$ are $C^\infty$, and the vector field on $M$ defined by $p\mapsto (p,\xi_p)$ is continuously differentiable. Hence, $(p,x)\mapsto\lVert \log_p(x)\rVert$ is $C^\infty$ at all $(p,x)\in M\times M$ except $p=x$, and $(p,x)\mapsto\langle \beta\xi_p,\log_p(x)\rangle$ is continuously differentiable at all $(p,x)$.
	\end{proposition}
		
		\begin{proof}
			
			(a) Note that $\gamma$ and $\gamma_1$, defined by $\gamma_1(t)=\exp_p(t\xi_p)$, are asymptotic. As noted in Section \ref{geom} of the main paper, $d(\gamma(t),\gamma_1(t))$ is monotonically non-increasing in $t$. Thus, for all $t\geq 0$,
			\begin{align*}
			t\bigg\lVert \frac{\log_p(\gamma(t))}{\lVert\log_p(\gamma(t))\rVert}-\xi_p\bigg\rVert=&\bigg\lVert (t-\lVert\log_p(\gamma(t))\rVert)\frac{\log_p(\gamma(t))}{\lVert\log_p(\gamma(t))\rVert}+\log_p(\gamma(t))-\log_p(\gamma_1(t))\bigg\rVert \\
   \leq&\lvert d(\gamma(0),\gamma(t))-d(p,\gamma(t))\rvert+\lVert \log_p(\gamma(t))-\log_p(\gamma_1(t))\rVert \\
   \leq& d(\gamma(0),p)+d(\gamma(t),\gamma_1(t)) \\
   \leq& 2d(\gamma(0),p),
			\end{align*}
			where norms are taken in the tangent space at $p$. The second inequality follows from the triangle inequality and the Cartan-Hadamard theorem. The result follows by dividing by $t$ and sending it to infinity.
			
			(b) Because geodesics corresponding to all $(p,v)\in TM$ can be defined on the entire real line, the geodesic flow on $TM$ (see Definition 2.4 in Ch. 3 of \citet{doCarmo1992}) defined by $\theta(t,p,v)=(\gamma_{p,v}(t),\dot{\gamma}_{p,v}(t))$, where $\gamma_{p,v}(t):=\exp_p(tv)$, can be defined on all of $\mathbb{R}\times TM$ and it is $C^\infty$ on this entire domain by the fundamental theorem on flows (Theorem 9.12 in \citet{Lee2012}). Then, since $\exp(p,v)=\pi(\theta(1,p,v))=\exp_p(v)$, where $\pi:TM\rightarrow M$ is the projection map, $\exp$ is also $C^\infty$.
			
			Since $\exp_p$ is a diffeomorphism between $T_pM\cong\mathbb{R}^n$ and $M$, the determinant of the differential of the $C^\infty$ map $(p,v)\mapsto (p,\exp(p,v))$ from $TM$ to $M\times M$ is
			\begin{align*}
			\det\begin{bmatrix}I & O \\ A & d\exp_p\end{bmatrix}=\det(I)\det(d\exp_p)\neq 0,
			\end{align*}
			where $d\exp_p$ is identified with its matrix representation in local coordinates and $A$ is some $n\times n$ matrix; thus, this map is a local diffeomorphism by the inverse function theorem. Since the inverse of the map exists and is precisely $(p,x)\mapsto (p,\log_p(x))$, it is a diffeomorphism.

   Finally, the continuous differentiability of the vector field was proven in Proposition 3.1 of \citet{Heintze1977}.
		\end{proof}
		
		\subsection{Proof of Proposition~\ref{basic}}\label{proof_prop3.4}
		\begin{proof}
			(a) Because $\rho(x,p^*;\beta,\xi)\in[(1-\beta)d(p^*,x),(1+\beta)d(p^*,x)]$ for all $(\beta,\xi)\in[0,1)\times\partial M$, 
   \begin{align*}
       G^{\beta,\xi}(p^*)&\leq E[(1+\beta)d(p^*,X)]\leq E\bigg[(1+\beta)\frac{\rho(x,p^*;\beta^*,\xi^*)}{1-\beta^*}\bigg]=\frac{1+\beta}{1-\beta^*}G^{\beta^*,\xi^*}(p^*)<\infty.
   \end{align*}
   
   So, the triangle inequality implies
			\begin{equation}\begin{aligned} \label{fin}
			G^{\beta,\xi}(p)&\leq(1+\beta)E[d(p,X)] \\
			&\leq(1+\beta)E[d(p,p^*)+d(p^*,X)]  \\
			&\leq(1+\beta)\bigg(d(p,p^*)+\frac{1}{1-\beta}E[\rho(X,p^*;\beta,\xi)]\bigg)  \\
			&\leq(1+\beta)\bigg(d(p,p^*)+\frac{1}{1-\beta}G^{\beta,\xi}(p^*)\bigg)  \\
			&<\infty.
			\end{aligned}\end{equation}
			So, both $G^{\beta,\xi}(p)$ and $E[d(p,X)]$ are finite for all $p\in M$. 
			
			For any sequence $\{p_j\}$ that converges to $p$, there exists some $J<\infty$ such that for all $j_1,j_2\geq J$, $d(p_{j_1},p_{j_2})<1$. Define a random variable $Y=(1+\beta)(1+\sum_{k=1}^J d(p_k,X))$. Then, $\rho(X,p_j;\beta,\xi)\leq Y$ for all $j<J$ and
			\begin{align*}
			\rho(X,p_j;\beta,\xi)&\leq(1+\beta)d(p_j,X)  \\
			&\leq(1+\beta)(d(p_j,p_J)+d(p_J,X)) \\
			&<(1+\beta)(1+d(p_J,X)) \\
			&\leq Y
			\end{align*}
			for all $j\geq J$. Thus, $\{\rho(X,p_j;\beta,\xi)\}_j$ is dominated by $Y$ and  $EY<\infty$ by (\ref{fin}). The map $p\mapsto\rho(x,p;\beta,\xi)$ is also known to be continuous as a function of its second argument (Proposition \ref{abc}(b)). So, $\rho(X,p_j;\beta,\xi)\rightarrow\rho(X,p;\beta,\xi)$ (almost) surely, and $G^{\beta,\xi}(p_j)\rightarrow G^{\beta,\xi}(p)$ by the dominated convergence theorem; thus, $G^{\beta,\xi}$ is continuous.
			
			(b) Letting $l:=\inf_{p\in M}G^{\beta,\xi}(p)$, there exists by the continuity of $G^{\beta,\xi}$ a sequence $\{q_j\}$ such that $G^{\beta,\xi}(q_j)$ goes to $l$. Then, $\{G^{\beta,\xi}(q_n)\}$ is a bounded sequence. By the triangle inequality, $d(q_1,q_j)\leq d(q_1,X)+d(X,q_j)\leq(\rho(X,q_1;\beta,\xi)+\rho(X,q_j;\beta,\xi))/(1-\beta)$, and taking expected values gives 
			\begin{equation} \label{bound}
			d(q_1,q_j)\leq \frac{G^{\beta,\xi}(q_1)+G^{\beta,\xi}(q_j)}{1-\beta}.
			\end{equation}
			This and the boundedness of $\{G^{\beta,\xi}(q_j)\}$ imply that $\{q_j\}$ is also bounded. Hadamard spaces are complete, so the bounded sequence $\{q_j\}$ has a subsequence $\{q_{j'}\}$ that converges to some $q'$. $\{G^{\beta,\xi}(q_{j'})\}$ converges to $l$, so by the continuity of $G^{\beta,\xi}$, $G^{\beta,\xi}(q')=l$ and the quantile set is nonempty. Since the quantile set is the preimage of the closed set $\{l\}$ under the continuous $G^{\beta,\xi}$, it is also closed. Replacing $q_1$ and $q_j$ in (\ref{bound}) with $q'$ and any other minimiser $m\in (G^{\beta,\xi})^{-1}(l)$, $d(q',m)\leq 2l/(1-\beta)$, so $(G^{\beta,\xi})^{-1}(l)$ is bounded too. Therefore, it is compact as connected Riemannian manifolds are proper.
		\end{proof}

		\subsection{Proof of Theorem~\ref{qgrad} }\label{proof_thm3.1}
		\begin{proof}
			Take an $h\in T_{q}M$ such that $\lVert h\rVert=1$ and let $\gamma_{q,h}(t)=\mathrm{exp}_{q}(th)$. We suppress the $q$ and $h$ parameters in the notation for $\gamma_{q,h}$. $\rho$ is differentiable as a function of its second argument on $M\backslash\{x\}$. So, if $x\neq q$,
			\begin{equation}\begin{aligned} \label{yes}
			\lim_{t\rightarrow 0+}\frac{\rho(x,\gamma(t);\beta,\xi)-\rho(x,\gamma(0);\beta,\xi)}{t}&=d\rho(x,\cdot;\beta,\xi)_{q}(h) \\
			&=\langle\nabla \rho(x,q;\beta,\xi),h\rangle.
			\end{aligned}\end{equation}
			If $x=q$, 
			\begin{equation}\begin{aligned} \label{no}
			\lim_{t\rightarrow 0+}\frac{\rho(x,\gamma(t);\beta,\xi)-\rho(x,\gamma(0);\beta,\xi)}{t}
			&=\lim_{t\rightarrow 0+}\frac{\lVert th\rVert+\langle \beta\xi_{\gamma(t)},\log_{\gamma(t)}(x)\rangle-0}{t}  \\
			&=\lim_{t\rightarrow 0+}\frac{t+\langle \beta\xi_{\gamma(t)},-t\gamma^\prime(t)\rangle}{t}  \\
			&=1-\langle \beta\xi_{q},h\rangle,
			\end{aligned}\end{equation}
			where we used continuity in the last equality.
			
			Take any sequence of positive numbers $t_1,t_2,\ldots$ that converges to 0. Fixing $j\in\mathbb{Z}_+$, if $x=\gamma(t')$ for some $t'\in(0,t_j)$,
			\begin{align*}
			&\frac{\lvert\rho(x,\gamma(t_j);\beta,\xi)-\rho(x,\gamma(0);\beta,\xi)\rvert}{t_j} \\
			&\leq~\frac{\lvert d(\gamma(t_j),\gamma(t'))+\langle \log_{\gamma(t_j)}(\gamma(t')),\xi_{\gamma(t_j)}\rangle\rvert+\lvert d(\gamma(0),\gamma(t'))+\langle \log_{\gamma(0)}(\gamma(t')),\xi_{\gamma(0)}\rangle\rvert}{t_j} \\
			&\leq~\frac{\vert t_j-t'\rvert\lVert h\rVert+\lvert t_j-t'\rvert\lVert h\rVert\lVert\xi_{\gamma(t_j)}\rVert+\vert t'\rvert\lVert h\rVert+\lvert t'\rvert\lVert h\rVert\lVert\xi_{\gamma(0)}\rVert}{t_j} \\
			&\leq~ 4;
			\end{align*}
			otherwise, by the mean value theorem, there exists some $t'\in(0,t_j)$ for which
			\begin{align*}
			\frac{\lvert\rho(x,\gamma(t_j);\beta,\xi)-\rho(x,\gamma(0);\beta,\xi)\rvert}{t_j} &=\lvert\langle\nabla \rho(x,\gamma(t');\beta,\xi),\dot\gamma(t')\rangle \rvert \\
			&\leq\lVert\nabla \rho(x,\gamma(t');\beta,\xi)\rVert \\
			&\leq\frac{\lVert-\log_{\gamma(t')}(x)\rVert}{\lVert \log_{\gamma(t')}(x)\rVert}+\lVert R(\gamma(t'),x)\rVert \\
			&=1+\lVert R(\gamma(t'),x)\rVert,
			\end{align*}
			where $R(p,x)$ is the gradient of $\langle\beta\xi_p,\log_p(x)\rangle$ as a function of $p$. Note that $\lVert R(p,x)\rVert$ is continuous as a function of $(p,x)$ by Proposition \ref{abc}(b). Then, $A:=\{(p,x):d(q,p)\leq \sup_{j\in\mathbb{Z}_+} t_j,~x\text{ is in the compact support of the distribution of $X$}\}$ is a compact subset of $M\times M$, and so $\sup_{(p,x)\in A}\lVert R(p,x)\rVert<\infty$. Thus, for all $j\in\mathbb{Z}_+$.
			\begin{equation}\begin{aligned} \label{bounded}
			\frac{\lvert\rho(X,\gamma(t_j);\beta,\xi)-\rho(X,\gamma(0);\beta,\xi)\rvert}{t_j}\leq 5+\sup_{(p,x)\in A}\lVert R(p,x)\rVert<\infty
			\end{aligned}\end{equation}
			almost surely.
			
			Therefore, we have the following, in which the first equality follows from (\ref{bounded}) and the bounded convergence theorem, and the second-to-last from (\ref{yes}) and (\ref{no}):
			\begin{equation}\begin{aligned}
			&\lim_{t\rightarrow 0+}\frac{G^{\beta,\xi}(\gamma(t))-G^{\beta,\xi}(\gamma(0))}{t}  \\
			&=E\bigg[\lim_{t\rightarrow 0+}\frac{\rho(X,\gamma(t);\beta,\xi)-\rho(X,\gamma(0);\beta,\xi)}{t}\bigg]  \\
			&=E\bigg[\lim_{t\rightarrow 0+}\frac{\rho(X,\gamma(t);\beta,\xi)-\rho(X,\gamma(0);\beta,\xi)}{t};X\neq q\bigg]  \\
			&\qquad\qquad+E\bigg[\lim_{t\rightarrow 0+}\frac{\rho(X,\gamma(t);\beta,\xi)-\rho(X,\gamma(0);\beta,\xi)}{t};X=q\bigg]  \\
			&=E[\langle\nabla \rho(X,q;\beta,\xi),h\rangle;X\neq q]+E[1-\langle \beta\xi_{q},h\rangle;X=q]  \\
			&=\langle E[\nabla \rho(X,q;\beta,\xi);X\neq q],h\rangle+P(X=q)(1-\langle \beta\xi_{q},h\rangle), 
			\end{aligned}\end{equation}
			and because $q$ is a minimiser, this quantity must be at least 0 for all $h\in T_qM$ such that $\lVert h\rVert=1$; that is, $\langle E[\nabla \rho(X,q;\beta,\xi);X\neq q]-P(X=q)\beta\xi_{q},-h\rangle\leq P(X=q)$
			for all such $h$, which is equivalent to the expression in the theorem.
		\end{proof}
		
		%\begin{proof}[Proof of Corollary~\ref{sampleqgrad}]
		%	This follows directly from Theorem \ref{qgrad} by letting $X$ be the random element that follows the empirical distribution determined by $X_1,\ldots,X_N$ and dividing the expression in Theorem \ref{qgrad} by $N$.
		%\end{proof}

		\subsection{Proof of Theorem~\ref{slln}}\label{proof_thm4.1}
		The proof of strong consistency is in part based on that of Theorem 2.3 in \citet{Bhattacharya2003}, but the proof of the following uniform law of large numbers is completely different.
	
	\begin{lemma} \label{sulln}
		\sloppy Let $M$ be a Hadamard manifold. For any compact $L\subset M$, $\lim_{N\rightarrow\infty}\sup_{p\in L}\lvert\hat{G}^{\beta,\xi}_N(p)-G^{\beta,\xi}(p)\rvert=0$ almost surely.
	\end{lemma}
	
		\begin{proof}
			$G^{\beta,\xi}$ is finite on all $M$ by Proposition \ref{basic}(a). Denote by $C(L)$ the vector space of continuous functions on $L$ equipped with the uniform (sup) norm (i.e., for $f\in C(L)$, $\lVert f\rVert_{\sup}:=\sup_{p\in L}\lvert f(p)\rvert$), which is a Banach space; because $L$ is a compact metric space, the Stone-Weierstrass theorem implies that $C(L)$ is separable. Then, the restrictions $\hat{G}^{\beta,\xi}_N|_L$ and $G^{\beta,\xi}|_L$ to $L$ are in $C(L)$ by Proposition \ref{basic}(a). Note that if $d(x,x')<\epsilon/(1+\beta)$, then 
			\begin{align*}
			\lVert\rho(x,\cdot;\beta,\xi)-\rho(x',\cdot;\beta,\xi)\rVert_{\sup}&\leq\sup_{p\in L}(\lvert d(p,x)-d(p,x')\rvert+\lvert\langle\beta\xi_p,\log_p(x)-\log_p(x')\rangle\rvert) \\
			&\leq\sup_{p\in L}(d(x,x')+\lVert\beta\xi_p\rVert\lVert\log_p(x)-\log_p(x')\rVert) \\
			&\leq \sup_{p\in L}(d(x,x')+\beta d(x,x')) \\
			&<\epsilon,
			\end{align*}
			where the second inequality follows from the triangle and Cauchy-Schwarz inequalities, and the third from the Cartan-Hadamard theorem. So, the function from $M$ to $C(L)$ defined by $x\mapsto\rho(x,\cdot;\beta,\xi)$ is continuous and hence measurable if $C(L)$ is equipped with its induced Borel $\sigma$-algebra. Therefore, $\rho(X,\cdot;\beta,\xi)|_L$, and hence each $\hat{G}^{\beta,\xi}_N|_L$, is a $C(L)$-valued random element. For any fixed $p'\in L$, noting $\mathrm{diam}(L):=\sup_{p,p_0\in L}d(p,p_0)<\infty$ by the compactness, and hence the boundedness, of $L\subset M$,
			\begin{align*}
			E[\lVert\rho(X,\cdot;\beta,\xi)\rVert_{\sup}]&=E\bigg[\sup_{p\in L}\rho(X,p;\beta,\xi)\bigg] \\
			&\leq E\bigg[\sup_{p\in L}(1+\beta)d(p,X)\bigg] \\
			&\leq E\bigg[\sup_{p\in L}(1+\beta)(d(p,p')+d(p',X))\bigg] \\
			&\leq E\bigg[(1+\beta)(\mathrm{diam}(L)+d(p',X))\bigg] \\
			%&\leq(1+\beta)\mathrm{diam}(L)+(1+\beta)E\bigg[\frac{\rho(X,p';\beta,\xi)}{1-\beta}\bigg] \\
			&=(1+\beta)\mathrm{diam}(L)+(1+\beta)G^{0,\xi}(p') \\
			&<\infty
			\end{align*}
			by Proposition \ref{basic}(a). Then, by the strong law of large numbers on separable Banach spaces (see \citet{Mourier1953} or Theorem 4.1.1 in \citet{Padgett1973}), $\lVert\hat{G}^{\beta,\xi}_N|_L-G^{\beta,\xi}|_L\rVert_{\sup}\rightarrow 0$ almost surely, from which the desired result follows. Note that the left side of this expression, which is equivalent to the expression in the lemma, is indeed a random variable because the metric function is continuous on metric spaces.
		\end{proof}

	\begin{lemma}\label{parta}		
		For any $\epsilon>0$, there exist some $\Omega_1\in\mathcal{F}$ and $N_1(\omega)<\infty$ for all $\omega\in\Omega_1$ such that $P(\Omega_1)=1$ and the sample $(\beta,\xi)$-quantile set of $X_1,\ldots,X_N$ is contained in $C^\epsilon=\{p\in M:d(p,q(\beta,\xi))<\epsilon\}$ for all $N\geq N_1(\omega)$.
	\end{lemma}
		
		\begin{proof}
			$G^{\beta,\xi}$ is finite on all $M$ by Proposition \ref{basic}(a). Recalling that $q(\beta,\xi)$ is nonempty by Proposition \ref{basic}(b), define $l$ to be $G^{\beta,\xi}(q)$ for any $q\in q(\beta,\xi)$, that is, $l=\min_{p\in M}G^{\beta,\xi}(p)$. We first show that there exist some compact $D\subset M$, $\Omega_2\in\mathcal{F}$ of $P$-measure 1 and an $N_2(\omega)<\infty$ for each $\omega\in\Omega_2$ for which $N\geq N_2(\omega)$ implies that 
			\begin{equation} \label{lp1}
			\hat{G}^{\beta,\xi}_N(p)>l+1
			\end{equation}
			for all $p\in M\backslash D$.
			
			For a fixed $q\in q(\beta,\xi)$,
			\begin{equation}\begin{aligned} \label{gq}
			\hat{G}^{\beta,\xi}_N(p)&\geq (1-\beta)\frac{1}{N}\sum_{i=1}^Nd(p,X_i)  \\
			&\geq (1-\beta)\frac{1}{N}\sum_{i=1}^N\lvert d(p,q)-d(q,X_i)\rvert \\
			&\geq (1-\beta)\frac{1}{N}\sum_{i=1}^N(d(p,q)-d(q,X_i)) \\
			&\geq (1-\beta)\bigg(d(p,q)-\frac{1}{1-\beta}\frac{1}{N}\sum_{i=1}^N\rho(X_i,q;\beta,\xi)\bigg)  \\
			&=(1-\beta)d(p,q)-\hat{G}^{\beta,\xi}_N(q).
			\end{aligned}\end{equation}
			By $G^{\beta,\xi}(q)=l<\infty$ and the almost sure convergence of $\hat{G}^{\beta,\xi}_N(q)$ to $G^{\beta,\xi}(q)$, there exist some $\Omega_2\in\mathcal{F}$ of $P$-measure 1 and an $N_2(\omega)<\infty$ for all $\omega\in\Omega_2$ for which $N\geq N_2(\omega)$ implies that $\hat{G}^{\beta,\xi}_N(q)<l+1$. Then, by (\ref{gq}), for $\Delta:=(2l+1)/(1-\beta)$, $d(p,q)>\Delta$, $\omega\in\Omega_2$ and $N\geq N_2(\omega)$ imply $\hat{G}^{\beta,\xi}_N(p)>l+1$. Then, defining $D:=\{p:d(p,q(\beta,\xi))\leq\Delta\}\supset{\{p:d(p,q)\leq\Delta\}}$, which is closed and bounded by Proposition \ref{basic}(b), $N\geq N_2(\omega)$ implies that $\hat{G}^{\beta,\xi}_N(p)>l+1$ for all $p\in M\backslash D$, so we have found the desired compact $D$.%Therefore $F(p)>l+1$ and $\hat{F}_N(p)$ are both greater than $l+1$ on $M\backslash D$ for all $N\geq N_5(\omega)$, where $\omega\in\Omega_5$.
			
			We now just need to show that there exist some $\theta(\epsilon)>0$, $\Omega_1\in\mathcal{F}$ of $P$-measure 1 and $N_1(\omega)<\infty$ for all $\omega\in\Omega_1$ such that $N\geq N_1(\omega), \omega\in\Omega_1$, implies
			\begin{equation}\begin{aligned} \label{final}
			&\hat{G}^{\beta,\xi}_N(p)\leq l+\frac{\theta(\epsilon)}{2} \qquad\text{ for all $p\in q(\beta,\xi)$ and} \\
			&\hat{G}^{\beta,\xi}_N(p)\geq l+\theta(\epsilon) \qquad\text{ for all $p\in M\backslash C^\epsilon$},
			\end{aligned}\end{equation}
			because then no minimiser of $\hat{G}^{\beta,\xi}_N$, that is, no element of the sample $(\beta,\xi)$-quantile set defined in the statement of the theorem, is in $M\backslash C^\epsilon$ if $\omega\in\Omega_1$ and $N\geq N_1(\omega)$, proving part (a) of the theorem. 
			
			$D_\epsilon:=D\cap(M\backslash C^\epsilon)$, being the intersection of a closed and bounded set and a closed one, is itself closed and bounded and hence compact because $M$. By continuity, $G^{\beta,\xi}$ attains its minimum $l_\epsilon:=\min_{p\in D_\epsilon}G^{\beta,\xi}(p)$ on the compact $D_\epsilon$, and since $D_\epsilon\cap q(\beta,\xi)$ is empty by definition, $l_\epsilon>l$, so there is some $\theta(\epsilon)\in(0,1)$ such that $l_\epsilon>l+2\theta(\epsilon)$. We now use Lemma \ref{sulln} for $L=D$, allowing us to find a set $\Omega_3\in\mathcal{F}$ and $N_3(\omega)<\infty$ for each $\omega\in\Omega_3$ such that $N\geq N_3(\omega)$ implies that $\sup_{p\in D}\lvert\hat{G}^{\beta,\xi}_N(p)-G^{\beta,\xi}(p)\rvert<\theta(\epsilon)/2$, and thus $\hat{G}^{\beta,\xi}_N(p)\leq G^{\beta,\xi}(p)+\theta(\epsilon)/2=l+\theta(\epsilon)/2$ for all $p\in q(\beta,\xi)\subset D$, satisfying the first part of (\ref{final}), and $\hat{G}^{\beta,\xi}_N(p)>G^{\beta,\xi}(p)-\theta(\epsilon)/2>l_\epsilon-\theta(\epsilon)>l+\theta(\epsilon)$ for all $p\in D_\epsilon\subset D$. This and (\ref{lp1}) imply that $\hat{G}^{\beta,\xi}_N(p)>l+\theta(\epsilon)$ for all $p\in D_\epsilon\cup(M\backslash D)\supset M\backslash C^\epsilon$ if $N\geq N_1(\omega)$ and $\omega\in\Omega_1$, where $\Omega_1=\Omega_2\cap\Omega_3$ and $N_1(\omega)=\max\{N_2(\omega),N_3(\omega)\}$, satisfying the second part of (\ref{final}).
		\end{proof}
		
		\begin{proof}[Proof of Theorem~\ref{slln}]
		The proposed convergences in part (b) hold for all $\omega\in\Omega_1$, and $P(\Omega_1)=1$, proving part (b).
		\end{proof}
		
		\subsection{Proofs for the results in Section \ref{asymp}}\label{proof_thm4.2}
		
		In this section, in addition to the already-introduced Euclidean inner product $\langle \cdot,\cdot\rangle_2$ and norm $\lVert \cdot\rVert_2$, and Riemannian inner product $\langle \cdot,\cdot\rangle$ and norm $\lVert \cdot\rVert$, we also make reference to the $\sup$, or $L_\infty$, norm of a vector with $\lVert\cdot\rVert_\infty$, and the Frobenius norm of a matrix with $\lVert\cdot\rVert_F$.
		
		The proof will require the following two lemmas.
		
		\begin{lemma} \label{lem}
			Let $M$ be an $n$-dimensional $(n\geq2)$ Hadamard manifold.
			\begin{itemize}
				\item[(a)] There exists a positive upper semi-continuous function $L:M\times M\rightarrow \mathbb{R}$ that satisfies
				\begin{align*}
				\lVert x-p\rVert_2=L(x,p)\lVert\log_p(x)\rVert.
				\end{align*}
				
				\item[(b)] There exists a positive upper semi-continuous function $S:M\times M\rightarrow \mathbb{R}$ that satisfies
				\begin{align*}
				\lvert D_{r'}\Psi^r(x,p;0,\xi)\rvert&\leq\frac{S(x,p)}{\lVert x-p\rVert_2},
				\end{align*}
				for all $r,r^\prime=1,\ldots,n$ whenever $x\neq p$.
			\end{itemize}
		\end{lemma}
		
		\begin{proof} (a) Consider the geodesic (according to the metric $g$) $\gamma_{p,x}:[0,1]\rightarrow M$ satisfying $\gamma_{p,x}(0)=p$, $\gamma_{p,x}(1)=x$. Denoting by $\zeta_{\min}(p)$ the smallest eigenvalue of $g_p$, $\zeta_{\min}$ is continuous because eigenvalues are the roots of polynomials whose coefficients are continuous functions of the entries of a matrix. Then define $\zeta_{\inf}(x,p):=\inf_{t\in[0,1]}\zeta_{\min}(\gamma_{p,x}(t))$. The length of the curve $\gamma_{p,x}$ measured in the Euclidean metric must be at least as large as the length of the straight line from $p$ to $x$, so
			\begin{align*}
			\sqrt{\zeta_{\inf}(x,p)}\lVert x-p\rVert_2&\leq\sqrt{\zeta_{\inf}(x,p)}\int_0^1 \lVert\dot\gamma_{p,x}(t)\rVert_2 dt \\
			&\leq\int_0^1\sqrt{\dot\gamma_{p,x}(t)^Tg_{\gamma_{p,x}(t)}\dot\gamma_{p,x}(t)}dt \\
			&=\int_0^1\lVert\dot\gamma_{p,x}(t)\rVert dt\\
			&=d(p,x) \\
			&=\lVert\log_p(x)\rVert,
			\end{align*}
			or
			\begin{align*}
			\frac{\lVert x-p\rVert_2}{\lVert\log_p(x)\rVert}\leq\frac{1}{\sqrt{\zeta_{\inf}(x,p)}}
			\end{align*}
			if $x\neq p$. Fix some $p'\in M$. By the continuity of $\zeta_{\min}$, there exists for any $\epsilon>0$ some $\delta>0$ for which $d(p,p')<\delta$ implies $\lvert 1/\sqrt{\zeta_{\min}(p)}-1/\sqrt{\zeta_{\min}(p')}\rvert<\epsilon$. All geodesic balls are convex on Hadamard manifolds, and because the image of $\gamma_{p,x}$ is compact, $\zeta_{\inf}(x,p)=\zeta_{\min}(\gamma_{p,x}(t))>0$ for some $t\in[0,1]$. Therefore, if $d(p,p')<\delta$ and $d(x,p')<\delta$, $d(p',\gamma_{p,x}(t))<\delta$ and $\lvert 1/\sqrt{\zeta_{\inf}(x,p)}-1/\sqrt{\zeta_{\min}(p')}\rvert<\epsilon$. This means that 
			\begin{equation}\begin{aligned} \label{uppersemi}
			\lim\sup_{(x,p)\rightarrow(p',p')}\frac{\lVert x-p\rVert_2}{\lVert\log_p(x)\rVert}\leq\lim_{(x,p)\rightarrow(p',p')}\frac{1}{\sqrt{\zeta_{\inf}(x,p)}}=\frac{1}{\sqrt{\zeta_{\min}(p')}}.
			\end{aligned}\end{equation}
			
			Then, $L(x,p)$ defined by  
			\begin{align*}
			L(x,p)=\begin{cases}
			\frac{\lVert x-p\rVert_2}{\lVert\log_p(x)\rVert}&~~\mbox{if $x\neq p$} \\ 
			\frac{1}{\sqrt{\zeta_{\min}(p)}} &~~ \mbox{if $x=p$}
			\end{cases}
			\end{align*}
			is continuous if $x\neq p$, while the upper semi-continuity at $(p,p)\in M\times M$ follows from (\ref{uppersemi}) and the continuity of $1/\sqrt{\zeta_{\min}(p)}$ on $M$.
			
			(b) The Euclidean gradient of $\lVert\log_p(x)\rVert$ as a function of $p$ is $\Psi(x,p;0,\xi)=-g_p\log_p(x)/\lVert\log_p(x)\rVert$ when $x\neq p$. Thus, the Euclidean Hessian matrix of the same function at $x\neq p$ is 
			\begin{align*}
			D\Psi(x,p;0,\xi)&=\frac{S^*(x,p)\lVert\log_p(x)\rVert-(g_p\log_p(x))(g_p\log_p(x))^T/\lVert\log_p(x)\rVert}{\lVert\log_p(x)\rVert^2} \\
			&=\frac{S^*(x,p)-(g_p\log_p(x))(g_p\log_p(x))^T/\lVert\log_p(x)\rVert^2}{\lVert\log_p(x)\rVert}
			\end{align*}
			for some $n\times n$ matrix function $S^*(x,p)$ that is $C^\infty$ on all of $M\times M$ (by the infinite differentiability of $g_p\log_p(x)$). Denoting the smallest and largest eigenvalues of $g_p$ by  $\zeta_{\min}(p)$ and $\zeta_{\max}(p)$, respectively, these are continuous on all of $M$ as functions of $p$, $\zeta_{\max}(p)\geq\zeta_{\min}(p)>0$ by positive definiteness,
			\begin{equation}\begin{aligned} \label{int1}
			\lVert\log_p(x)\rVert^2=\log_p(x)^Tg_p\log_p(x)\in\big[\zeta_{\min}(p)\lVert\log_p(x)\rVert_2^2,~\zeta_{\max}(p)\lVert\log_p(x)\rVert_2^2\big],
			\end{aligned}\end{equation}
			and
			\begin{equation}\begin{aligned}  \label{int2}
			&(g_p\log_p(x))^T(g_p\log_p(x))\\
   &=\log_p(x)^Tg_p^2\log_p(x)\in\big[\zeta_{\min}(p)^2\lVert\log_p(x)\rVert_2^2,~\zeta_{\max}(p)^2\lVert\log_p(x)\rVert_2^2\big].
			\end{aligned}\end{equation}
			Because of this, and because the absolute value of any entry of $S^*(x,p)$ is less than or equal to the Frobenius norm of $S^*(x,p)$ $\lVert S^*(x,p)\rVert_F$, and the absolute value of any entry of $(g_p\log_p(x))(g_p\log_p(x))^T$ is less than or equal to $(g_p\log_p(x))^T(g_p\log_p(x))$, 
			\begin{align*}
			\lvert D_{r'}\Psi^r(x,p;0,\xi)\rvert&\leq\frac{\lVert S^*(x,p)\rVert_F+\zeta_{\max}(p)^2\lVert\log_p(x)\rVert_2^2/\zeta_{\min}(p)\lVert\log_p(x)\rVert_2^2}{\lVert x-p\rVert_2}\frac{\lVert x-p\rVert_2}{\lVert\log_p(x)\rVert} \\
			&=\frac{(\lVert S^*(x,p)\rVert_F+\zeta_{\max}(p)^2/\zeta_{\min}(p))L(x,p)}{\lVert x-p\rVert_2} \\
			&=\frac{S(x,p)}{\lVert x-p\rVert_2},
			\end{align*}
			when $x\neq p$, where $L(x,p)$ is taken from (a) and $S(x,p)$ is a positive real function that is upper semi-continuous on all of $M\times M$.
		\end{proof}
		
		For $j\in\{2,\infty\}$, we use $u_{k,j}$ defined as
		\begin{align*}
		u_{k,j}(x,q,d)=\sup_{\tau:\lVert \tau-q\rVert_j\leq d}\lVert\Psi_k(x,\tau)-\Psi_k(x,q)\rVert_j,
		\end{align*}
		and
		\begin{align*}
		Z_{k,j}(x_1,\ldots,x_N,\tau,q)=\frac{\lVert\sum_{i=1}^N(\Psi_k(x_i,\tau)-\Psi_k(x_i,q)-\lambda_k(\tau)+\lambda_k(q))\rVert_j}{\sqrt{N}+N\lVert\lambda_k(\tau)\rVert_j}
		\end{align*}
		as in \citet{Huber1967}.
		
		\begin{lemma} \label{meas}
			Let $M$ be an $n$-dimensional $(n\geq2)$ Hadamard manifold and fix $q\in M$, $k\in\{1,\ldots,K\}$ and $j\in\{2,\infty\}$. Assume that $A_j:=\{\tau:\lVert \tau-q\rVert_j\leq d\}\subset M$.
			\begin{itemize}
				\item[(a)] $u_{k,j}(x,q,d)$ is measurable as a function of $x$.
				
				\item[(b)] Assuming $\lambda_k$ is continuous on $A_j$, $\sup_{\tau\in A_j}Z_{k,j}(x_1,\ldots,x_N,\tau,q)$ is measurable as a function of $x_1$,\ldots,$x_N$.
			\end{itemize}
		\end{lemma}
		
		\begin{proof}
			(a) Defining $R_k(x,p)$ to be the Riemannian gradient of $\langle\beta_k(\xi_k)_p,\log_p(x)\rangle$ as a function of $p$, $R_k$ is continuous on $M\times M$; see Proposition \ref{abc}(b). Then, $\Psi_k(x,\tau)=-g_\tau\log_\tau(x)/\lVert\log_\tau(x)\rVert+g_\tau R_k(x,\tau)$ when $x\neq\tau$, and noting that (\ref{no}) contains a proof that $R_k(x,p)=-\beta_k(\xi_k)_p$, then $\Psi_k(x,\tau)=g_\tau R_k(x,\tau)$ when $x=\tau$.
			
			The map defined by $(p,v)\mapsto(p,\exp_p(v))$ for $(p,v)\in TM$ has differential 
			\begin{align*}
			\begin{bmatrix}I & O \\ I & I\end{bmatrix}
			\end{align*}
			at $(x,0)\in TM$ since $d(\exp_x)_{0}=I$. Then its inverse $(p,b)\mapsto (p,\log_p(b))$ is $C^\infty$ by Proposition \ref{abc}(b) and has differential
			\begin{align*}
			\begin{bmatrix}I & O \\ -I & I\end{bmatrix}
			\end{align*}
			at $(x,x)$, so $d(\log_{\cdot}(x))_x=-I$.
			%where $\log(p,b):=\log_p(b)$.
			
			Then, by Taylor's theorem for multivariate functions, 
			\begin{equation}\begin{aligned}\label{mtaylors}
			\log_\tau(x)&=\log_x(x)+d(\log_{\cdot}(x))_x(\tau-x)+F(\tau)(\tau-x) \\
			&=(I-F(\tau))(x-\tau),
			\end{aligned}\end{equation}
			where $F:M\rightarrow \mathbb{R}^{n\times n}$ satisfies  $\lim_{\tau\rightarrow x}F(\tau)=0$. Therefore, for any $C^\infty$ path $c:[0,1]\rightarrow M$ satisfying $c(0)=x$, 
			\begin{equation}\begin{aligned} \label{lim}
			\frac{\log_{c(t)}(x)}{\lVert \log_{c(t)}(x)\rVert_2}&=\frac{(I-F(c(t))\frac{(c(0)-c(t))/t}{\lVert (c(0)-c(t))/t\rVert_2}}{\Big\lVert (I-F(c(t))\frac{(c(0)-c(t))/t}{\lVert (c(0)-c(t))/t\rVert_2}\Big\rVert_2} \\
			&\rightarrow \frac{(I-0)\frac{-\dot c(0)}{\lVert -\dot c(0)\rVert_2}}{\Big\lVert (I-0)\frac{-\dot c(0)}{\lVert -\dot c(0)\rVert_2}\Big\rVert_2} \\
			&=-\frac{\dot c(0)}{\lVert \dot c(0)\rVert_2}
			\end{aligned}\end{equation}
			as $t\rightarrow 0$.
			
			Suppose $j=2$. For any $p\in A_2$ and vector $v\in \mathbb{R}^n$, one can construct a $C^\infty$ path $c_2:[0,1]\rightarrow M$ satisfying $c_2(t)=p$ if and only if $t=0$, $c_2([0,1])\subset A_2$ and $\langle g_p\dot c_2(0),v\rangle_2=0$ in the following manner. Take any $w$ that is orthogonal (in the Euclidean inner product) to $v$. If $p$ is in the interior of $A_2$, or it is on the boundary $\partial A_2$ of $A_2$ and $g_p^{-1}w$ is not tangent to $\partial A_2$, which is diffeomorphic to $S^{n-1}$, at $p$, there exists a possibly negative $a$ that is sufficiently small in absolute value such that $p+ag_p^{-1}w\in A_2$; this is because if the Euclidean ray from $p$ in the direction of $g_p^{-1}w$ is not initially in $A_2$, then the Euclidean ray in the direction of $-g_p^{-1}w$ is. Then, $c_2(t)=p+tag_p^{-1}w$ satisfies the desiderata. On the other hand, if $p$ is on $\partial A_2$ and $g_p^{-1}w$ is tangent to $\partial A_2$ at $p$, $c_2(t)=\exp_p^{\partial A_2}(tg_p^{-1}w)$, where $\exp_p^{\partial A_2}$ is the exponential map at $p$ on $\partial A_2\cong S^{n-1}$. This $c_2$ also satisfies the desiderata.
			
			Suppose that $x\in A_2$ and let $p=x$ and $v=\Psi_k(x,x)-\Psi_k(x,q)$. Then, given (\ref{lim}), Lemma \ref{lem}(a), the first paragraph of this proof, and the properties of $c_2$,
			\begin{equation}\begin{aligned} \label{examp}
			&\lim\inf_{t\rightarrow 0}\lVert\Psi_k(x,c_2(t))-\Psi_k(x,q)\rVert_2^2 \\
			&=\lim\inf_{t\rightarrow 0}\bigg\lVert -g_{c_2(t)}\frac{\log_{c_2(t)}(x)}{\lVert\log_{c_2(t)}(x)\rVert}+g_{c_2(t)} R_k(x,c_2(t))-\Psi_k(x,q)\bigg\rVert_2^2 \\
			&\geq\lim\inf_{t\rightarrow 0}\bigg\lVert -g_{c_2(t)}\frac{\log_{c_2(t)}(x)}{\lVert\log_{c_2(t)}(x)\rVert}\bigg\rVert_2^2+\lim_{t\rightarrow 0}\lVert g_{c_2(t)} R_k(x,c_2(t))-\Psi_k(x,q)\rVert_2^2 \\
			&\qquad-2\lim\sup_{t\rightarrow 0}\bigg\lvert\bigg\langle-g_{c_2(t)}\frac{\log_{c_2(t)}(x)}{\lVert\log_{c_2(t)}(x)\rVert}, g_{c_2(t)} R_k(x,c_2(t))-\Psi_k(x,q)\bigg\rangle_2\bigg\rvert \\
			&\geq 0+\lVert g_xR_k(x,x)-\Psi_k(x,q)\rVert_2^2-2\lim\sup_{t\rightarrow 0}L(x,c_2(t)) \\
			&\qquad\times\bigg\lvert \bigg\langle -g_{c_2(t)}\frac{\log_{c_2(t)}(x)}{\lVert\log_{c_2(t)}(x)\rVert_2}, g_{c_2(t)} R_k(x,c_2(t))-\Psi_k(x,q)\bigg\rangle_2\bigg\rvert \\
			&=\lVert \Psi_k(x,x)-\Psi_k(x,q)\rVert_2^2-2\lim\sup_{t\rightarrow 0}L(x,c_2(t)) \\
			&\qquad\times\lim_{t\rightarrow 0}\bigg\lvert\bigg\langle-g_{c_2(t)}\frac{\log_{c_2(t)}(x)}{\lVert\log_{c_2(t)}(x)\rVert_2}, g_{c_2(t)} R_k(x,c_2(t))-\Psi_k(x,q)\bigg\rangle_2\bigg\rvert \\
			&\geq\lVert \Psi_k(x,x)-\Psi_k(x,q)\rVert_2^2 \\
			&\qquad-2L(x,x)\bigg\lvert\bigg\langle g_x\frac{\dot c_2(0)}{\lVert \dot c_2(0)\rVert_2}, g_x R_k(x,x)-\Psi_k(x,q)\bigg\rangle_2\bigg\rvert \\
			&=\lVert \Psi_k(x,x)-\Psi_k(x,q)\rVert_2^2.
			\end{aligned}\end{equation}
			Since $c_2(t)=x$ if and only if $t=0$, this proves that $\sup_{\tau\in A_2\backslash\{x\}}\lVert\Psi_k(x,\tau)-\Psi_k(x,q)\rVert_2\geq \lVert\Psi_k(x,x)-\Psi_k(x,q)\rVert_2$ if $x\in A_2$, and therefore that
			\begin{equation}\begin{aligned} \label{sup}
			u_{k,2}(x,q,d)=\sup_{\tau\in A_2\backslash\{x\}}\lVert\Psi_k(x,\tau)-\Psi_k(x,q)\rVert_2
			\end{aligned}\end{equation}
			for all $x\in M$. There exists a countable set $B:=\{b_1,b_2,\ldots\}\subset A_2\backslash \{x\}$ for which 
   \begin{equation}\begin{aligned} \label{sup2}
   \sup_{\tau\in B}\lVert\Psi_k(x,\tau)-\Psi_k(x,q)\rVert_2=\sup_{\tau\in A_2\backslash\{x\}}\lVert\Psi_k(x,\tau)-\Psi_k(x,q)\rVert_2,
   \end{aligned}\end{equation}
   and for each $l\in \mathbb{Z}_+$ a sequence of points in $A_2\backslash \{x\}$ with rational coordinates that converges to $b_l$; call the set of points in this sequence $B_l$. For all $x\in M$, $\lVert\Psi_k(x,\tau)-\Psi_k(x,q)\rVert_2$ is continuous as a function of $\tau$ on $A_2\backslash \{x\}$, so $\lVert\Psi_k(x,b_l)-\Psi_k(x,q)\rVert_2\leq \sup_{\tau\in B_l}\lVert\Psi_k(x,\tau)-\Psi_k(x,q)\rVert_2$, and therefore by (\ref{sup}) and (\ref{sup2}), $u_{k,2}(x,q,d)=\sup_{\tau\in B}\lVert\Psi_k(x,\tau)-\Psi_k(x,q)\rVert_2\leq\sup_{\tau\in \cup_{l=1}^\infty B_l}\lVert\Psi_k(x,\tau)-\Psi_k(x,q)\rVert_2\leq\sup_{\tau\in \mathbb{Q}^n\cap A_2}\lVert\Psi_k(x,\tau)-\Psi_k(x,q)\rVert_2$, since each $B_l$ is a subset of $\mathbb{Q}^n\cap A_2$. Clearly $u_{k,2}(x,q,d)\geq\sup_{\tau\in \mathbb{Q}^n\cap A_2}\lVert\Psi_k(x,\tau)-\Psi_k(x,q)\rVert_2$ also holds, so $u_{k,2}(x,q,d)=\sup_{\tau\in \mathbb{Q}^n\cap A_2}\lVert\Psi_k(x,\tau)-\Psi_k(x,q)\rVert_2$ for all $x\in M$. The desired conclusion follows for $j=2$ since the $\sup$ of countably many measurable functions is measurable.
			
			Now suppose $j=\infty$. Denote by $V$ the set of corners of the hypercube $A_\infty$ and by $\pi^1,\ldots,\pi^n:\mathbb{R}^n\rightarrow\mathbb{R}$ the projection functions defined by $\pi^{r'}(p^1,\ldots,p^n)=p^{r'}$ for $r'=1,\ldots,n$, which are continuous. Choose any $p\in A_\infty\backslash\{V\}$ and vector $v\in \mathbb{R}^n$. Let $r$ be any index for which $\lvert\pi^r(v)\rvert=\lVert v\rVert_\infty$ (that is, $r\in\arg\max_{r'\in\{1,\ldots,n\}}\lvert\pi^{r'}(v)\rvert$). The $(n-1)$-dimensional hyperplane $\{p+g_p^{-1}w:\langle w,e^r\rangle_2=0\}$, where $e^r$ is the $r$th standard unit vector, passes through $p\not\in V$ and hence, intersects the hypercube $A_\infty$ at some point $y\neq p$ (note that this might not be guaranteed if $p\in V$). Defining $c_\infty:[0,1]\rightarrow M$ by $c_\infty(t)=p+t(y-p)$, we have that $c_\infty(t)=p$ if and only if $t=0$, $c_\infty([0,1])\subset A_\infty$ and $\langle g_p\dot c_\infty(0),e^r\rangle_2=0$. Supposing $x\in A_\infty\backslash V$ and letting $p=x$ and $v=\Psi_k(x,x)-\Psi_k(x,q)$, by these properties of $c_\infty$, and again by (\ref{lim}), Lemma \ref{lem}(a) and the first paragraph of this proof,
			\begin{align*} 
			&\lim\inf_{t\rightarrow 0}\lVert\Psi_k(x,c_\infty(t))-\Psi_k(x,q)\rVert_\infty \\
			&=\lim\inf_{t\rightarrow 0}\bigg\lVert -g_{c_\infty(t)}\frac{\log_{c_\infty(t)}(x)}{\lVert\log_{c_\infty(t)}(x)\rVert}+g_{c_\infty(t)} R_k(x,c_\infty(t))-\Psi_k(x,q)\bigg\rVert_\infty \\
			&\geq\lim\inf_{t\rightarrow 0}\bigg\lvert\pi^r(g_{c_\infty(t)} R_k(x,c_\infty(t))-\Psi_k(x,q))+\pi^r\bigg(-g_{c_\infty(t)}\frac{\log_{c_\infty(t)}(x)}{\lVert\log_{c_\infty(t)}(x)\rVert}\bigg)\bigg\rvert \\
			&\geq\lim\inf_{t\rightarrow 0}\Bigg(\lvert\pi^r(g_{c_\infty(t)} R_k(x,c_\infty(t))-\Psi_k(x,q))\rvert-\bigg\lvert\pi^r\bigg(-g_{c_\infty(t)}\frac{\log_{c_\infty(t)}(x)}{\lVert\log_{c_\infty(t)}(x)\rVert}\bigg)\bigg\rvert\Bigg) \\
			&\geq \lim_{t\rightarrow 0}\lvert\pi^r(g_{c_\infty(t)} R_k(x,c_\infty(t))-\Psi_k(x,q))\rvert \\
			&\qquad-\lim\sup_{t\rightarrow 0}L(x,c_\infty(t))\bigg\lvert\pi^r\bigg(-g_{c_\infty(t)}\frac{\log_{c_\infty(t)}(x)}{\lVert\log_{c_\infty(t)}(x)\rVert_2}\bigg)\bigg\rvert \\
			&=\lvert\pi^r(g_x R_k(x,x)-\Psi_k(x,q))\rvert \\
			&\qquad-\lim\sup_{t\rightarrow 0}L(x,c_\infty(t))
			\lim_{t\rightarrow 0}\bigg\lvert\pi^r\bigg(-g_{c_\infty(t)}\frac{\log_{c_\infty(t)}(x)}{\lVert\log_{c_\infty(t)}(x)\rVert_2}\bigg)\bigg\rvert \\
			&\geq\lvert\pi^r(\Psi_k(x,x)-\Psi_k(x,q))\rvert-L(x,x)\bigg\lvert\pi^r\bigg(g_x\frac{\dot c_\infty(0)}{\lVert\dot c_\infty(0)\rVert_2}\bigg)\bigg\rvert \\
			&=\lVert\Psi_k(x,x)-\Psi_k(x,q)\rVert_\infty.
			\end{align*}
			Since $c_\infty(t)=x$ if and only if $t=0$, this proves that $\sup_{\tau\in A_\infty\backslash\{x\}}\lVert\Psi_k(x,\tau)-\Psi_k(x,q)\rVert_\infty\geq \lVert\Psi_k(x,x)-\Psi_k(x,q)\rVert_\infty$ if $x\in A_\infty\backslash V$.%, and therefore that 
			%\begin{align*}
			%u_{k,\infty}(x,q,d)=\sup_{\tau\in A_\infty\backslash\{x\}}\lVert\Psi_k(x,\tau)-\Psi_k(x,q)\rVert_\infty
			%\end{align*}
			%for all $x\in M\backslash V$. 
   Now an analogous argument to that of the $j=2$ case shows that for all $x\in M\backslash V$, $u_{k,\infty}(x,q,d)=\sup_{\tau\in \mathbb{Q}^n\cap A_\infty}\lVert\Psi_k(x,\tau)-\Psi_k(x,q)\rVert_\infty$. Then, noting that $V=\{p_1,\ldots,p_{2^n}\}$ consists of finitely many elements gives
			\begin{align*}
			&u_{k,\infty}(x,q,d) \\
			&=\bigg(\sup_{\tau\in \mathbb{Q}^n\cap A_\infty}\lVert\Psi_k(x,\tau)-\Psi_k(x,q)\rVert_\infty\bigg)I(x\not\in V)+\sum_{l=1}^{2^n}u_{k,\infty}(p_l,q,d)I(x=p_l),
			\end{align*}
			which is measurable since the $\sup$ of countably many measurable functions is measurable.
			
			(b) This proof is very similar to that of (a), with slight modifications. Fix $i'\in\{1,\ldots,N\}$.
			
			Dealing first with the $j=2$ case, suppose $x_{i'}\in A_2$ and define $c_2$ as described in the proof of (a) for $p=x_{i'}$ and $v=\sum_{i=1}^N(\Psi_k(x_i,x_{i'})-\Psi_k(x_i,q)-\lambda_k(x_{i'})+\lambda_k(q))$. Let $N':=\lvert\{i:x_i=x_{i'}, 1\leq i\leq N\}\rvert$ be the number of $x_i$ which equal $x_{i'}$. Recalling the continuity of $\lambda_k$ in $A_2$, we follow a similar logic to that of (\ref{examp}):
			\begin{align*}
			&\lim\inf_{t\rightarrow 0}\bigg\lVert\sum_{i=1}^N(\Psi_k(x_i,c_2(t))-\Psi_k(x_i,q)-\lambda_k(c_2(t))+\lambda_k(q))\bigg\rVert_2^2 \\
			&=\lim\inf_{t\rightarrow 0}\bigg\lVert\sum_{i:x_i=x_{i'}} -g_{c_2(t)}\frac{\log_{c_2(t)}(x_i)}{\lVert \log_{c_2(t)}(x_i)\rVert}+\sum_{i:x_i\neq x_{i'}} -g_{c_2(t)}\frac{\log_{c_2(t)}(x_i)}{\lVert \log_{c_2(t)}(x_i)\rVert} \\
			&\qquad+\sum_{i=1}^N(g_{c_2(t)}R_k(x_i,c_2(t))-\Psi_k(x_i,q)-\lambda_k(c_2(t))+\lambda_k(q))\bigg\rVert_2^2\\
			&\geq(N')^2\lim\inf_{t\rightarrow 0}\bigg\lVert -g_{c_2(t)}\frac{\log_{c_2(t)}(x_{i'})}{\lVert \log_{c_2(t)}(x_{i'})\rVert}\bigg\rVert_2^2+\lim_{t\rightarrow 0}\bigg\lVert \sum_{i:x_i\neq x_{i'}} -g_{c_2(t)}\frac{\log_{c_2(t)}(x_i)}{\lVert \log_{c_2(t)}(x_i)\rVert} \\
			&\qquad+\sum_{i=1}^N(g_{c_2(t)}R_k(x_i,c_2(t))-\Psi_k(x_i,q)-\lambda_k(c_2(t))+\lambda_k(q))\bigg\rVert_2^2 \\
			&\qquad-2N'\lim\sup_{t\rightarrow 0}\bigg\lvert\bigg\langle-g_{c_2(t)}\frac{\log_{c_2(t)}(x_{i'})}{\lVert\log_{c_2(t)}(x_{i'})\rVert},\sum_{i:x_i\neq x_{i'}} -g_{c_2(t)}\frac{\log_{c_2(t)}(x_i)}{\lVert \log_{c_2(t)}(x_i)\rVert}  \\
   &\qquad+\sum_{i=1}^N(g_{c_2(t)}R_k(x_i,c_2(t))-\Psi_k(x_i,q)-\lambda_k(c_2(t))+\lambda_k(q))\bigg\rangle_2\bigg\rvert \\
			&\geq\bigg\lVert \sum_{i=1}^N(\Psi_k(x_i,x_{i'})-\Psi_k(x_i,q)-\lambda_k(x_{i'})+\lambda_k(q))\bigg\rVert_2^2.
			\end{align*}
			Then, by the continuity of $\lambda_k$ in $A_2$,
			\begin{align*}
			&\lim\inf_{t\rightarrow 0} Z_{k,2}(x_1,\ldots,x_N,c_2(t),q) \\
   &=\frac{\lim\inf_{t\rightarrow 0}\lVert\sum_{i=1}^N(\Psi_k(x_i,c_2(t))-\Psi_k(x_i,q)-\lambda_k(c_2(t))+\lambda_k(q))\rVert_2}{\lim_{t\rightarrow 0}(\sqrt{N}+N\lVert\lambda_k(c_2(t))\rVert_2)} \\
			&\geq\frac{\lVert\sum_{i=1}^N(\Psi_k(x_i,x_{i'})-\Psi_k(x_i,q)-\lambda_k(x_{i'})+\lambda_k(q))\rVert_2}{\sqrt{N}+N\lVert\lambda_k(x_{i'})\rVert_2} \\
			&=Z_{k,2}(x_1,\ldots,x_N,x_{i'},q).
			\end{align*}
			This means that if $x_i\in A_2$, the $\sup$ of $Z_{k,2}(x_1,\ldots,x_N,\tau,q)$ for $\tau$ in some arbitrarily small neighbourhood of $x_{i'}$ sans $x_{i'}$ itself is at least as large as $Z_{k,2}(x_1,\ldots,x_N,x_{i'},q)$, and therefore, $\sup_{\tau\in A_2\backslash\{x_1,\ldots,x_N\}}Z_{k,2}(x_1,\ldots,x_N,\tau,q)\geq Z_{k,2}(x_1,\ldots,x_N,x_{i'},q)$. Thus,
			\begin{align*}
			\sup_{\tau\in A_2}\lVert Z_{k,2}(x_1,\ldots,x_N,\tau,q)\rVert_2=\sup_{\tau\in A_2\backslash\{x_1,\ldots,x_N\}}\lVert Z_{k,2}(x_1,\ldots,x_N,\tau,q)\rVert_2.
			\end{align*}
			The rest of the argument for the measurability of $\sup_{\tau\in A_2}Z_{k,2}(x_1,\ldots,x_N,\tau,q)$ is analogous to that for $u_{k,2}(x,q,d)$. The argument for the $j=\infty$ case is now easily seen by referring to earlier arguments.
			%Define $c_\infty$ using $p=x_{i'}$ and $v=\sum_{i=1}^N(\Psi_k(x_{i'},x_i)-\Psi_k(q,x_i)-\lambda_k(x_{i'})+\lambda_k(q))$, supposing that $x_{i'}\in A_\infty$. Then
			%\begin{align*}
			%&\lim\inf_{t\rightarrow 0}\bigg\lVert\sum_{i=1}^N(\Psi_k(c_\infty(t),x_i)-\Psi_k(q,x_i)-\lambda_k(c_\infty(t))+\lambda_k(q))\bigg\rVert_\infty \\
			%=&\lim\inf_{t\rightarrow 0}\bigg\lVert \sum_{i:x_i=x_{i'}} -g_{c_\infty(t)}\frac{\log_{c_\infty(t)}(x_i)}{\lVert \log_{c_\infty(t)}(x_i)\rVert}+\sum_{i:x_i\neq x_{i'}} -g_{c_\infty(t)}\frac{\log_{c_\infty(t)}(x_i)}{\lVert \log_{c_\infty(t)}(x_i)\rVert} \\
			%&\qquad+\sum_{i=1}^N(g_{c_\infty(t)}R(c_\infty(t),x_i)-\Psi_k(q,x)-\lambda_k(c_\infty(t))+\lambda_k(q))\bigg\rVert_\infty \\
			%\geq& \lim\inf_{t\rightarrow 0}\bigg\lvert\pi^r\bigg(\sum_{i:x_i\neq x_{i'}} -g_{c_\infty(t)}\frac{\log_{c_\infty(t)}(x_i)}{\lVert \log_{c_\infty(t)}(x_i)\rVert}+\sum_{i=1}^N(g_{c_\infty(t)}R(c_\infty(t),x_i)-\Psi_k(q,x)-\lambda_k(c_\infty(t))+\lambda_k(q))\bigg)\bigg\rvert \\
			%&\qquad-N'\lim\sup_{t\rightarrow 0}L(c_\infty(t),x)\bigg\lvert\pi^r\bigg(-g_{c_\infty(t)}\frac{\log_{c_\infty(t)}(x_i)}{\lVert \log_{c_\infty(t)}(x_i)\rVert}\bigg)\bigg\rvert \\
			%=&\bigg\lVert\sum_{i=1}^N(\Psi_k(x_{i'},x_i)-\Psi_k(q,x_i)-\lambda_k(x_{i'})+\lambda_k(q))\bigg\rVert_\infty,
			%\end{align*}
			%and therefore $\lim\inf_{t\rightarrow 0} Z_{k,2}(c_\infty(t),q,x_1,...,x_N)=Z_{k,2}(x_{i'},q,x_1,...,x_N)$, where we used the continuity of $\lambda_k$ in $A_\infty$. The rest of the argument is easy to see after referring to the arguments in the $j=2$ case and the $j=\infty$ case in (a).
		\end{proof}
		
		\begin{proof}[Proof of Theorem~\ref{clt}]
			This proof will extensively reference Section 4 of \citet{Huber1967}, in particular Theorem 3 and its Corollary; those results are also briefly described in Section 6.3 of \citet{Huber1981}. %We equip $\phi(M)$ with the standard Euclidean metric; that is, all norms $\lVert \cdot\rVert$ and inner products $\langle\cdot,\cdot\rangle_2$ in $\phi(M)$ are the standard Euclidean ones and not the Riemannian ones unless otherwise stated. 
			
			Fix a $k$ in $\{1,\ldots,K\}$. By absolute continuity in the neighbourhood $Q_k$ of (I), $P(X_i=X_j\in Q_k\text{ for some $i\neq j$})=0$, and $\hat{q}_{k,N}\in Q_k$ for all sufficiently large $N$ almost surely by Theorem \ref{slln}. Therefore, on the measurable set of probability 1 in which $X_i\neq X_j$ if $X_i\in Q_k$ for all $i,j\in\mathbb{Z}_+$ and $\hat{q}_{k,N}$ converges to $q_k$, $\sum_{i=1}^N I(X_i=\hat{q}_{k,N})\leq 1$ for all $N$ sufficiently large. Then by applying Theorem \ref{qgrad} to a random element distributed according to the empirical measure defined by the observations $X_1,...,X_N$, and noting that $g_q$ is symmetric and invertible for all $q\in\phi(M)$,
			\begin{align*}
			&\bigg\lVert\big(g_{\hat{q}_{k,N}}\big)^{-\frac{1}{2}}\frac{1}{\sqrt{N}}\sum_{i=1}^N\Psi_k(X_i,\hat{q}_{k,N})\bigg\rVert_2 \\
			&=\frac{1}{\sqrt{N}}\sqrt{\bigg(\sum_{i=1}^N\Psi_k(X_i,\hat{q}_{k,N})\bigg)^T\big(g_{\hat{q}_{k,N}}\big)^{-1}\bigg(\sum_{i=1}^N\Psi_k(X_i,\hat{q}_{k,N})\bigg)} \\
			&=\frac{1}{\sqrt{N}}\sqrt{\bigg(\big(g_{\hat{q}_{k,N}}\big)^{-1}\sum_{i=1}^N\Psi_k(X_i,\hat{q}_{k,N})\bigg)^Tg_{\hat{q}_{k,N}}\bigg(\big(g_{\hat{q}_{k,N}}\big)^{-1}\sum_{i=1}^N\Psi_k(X_i,\hat{q}_{k,N})\bigg)} \\
			&=\frac{1}{\sqrt{N}}\bigg\lVert\big(g_{\hat{q}_{k,N}}\big)^{-1}\sum_{i=1}^N\Psi_k(X_i,\hat{q}_{k,N})\bigg\rVert \\
			&\leq \frac{1}{\sqrt{N}}
			\end{align*}
			for all sufficiently large $N$ on the aforementioned set of probability 1. Therefore,
			\begin{align*}
			\big(g_{\hat{q}_{k,N}}\big)^{-\frac{1}{2}}\frac{1}{\sqrt{N}}\sum_{i=1}^N\Psi_k(X_i,\hat{q}_{k,N})\longrightarrow 0
			\end{align*}
			almost surely. This result and the almost sure convergence of
			\begin{align*}
			\big(g_{\hat{q}_{k,N}}\big)^{\frac{1}{2}}\longrightarrow \big(g_{q_k}\big)^{\frac{1}{2}},
			\end{align*}
			which follows from Theorem \ref{slln}, the continuous mapping theorem, and the continuity of the Riemannian metric, imply
			\begin{align*}
			\frac{1}{\sqrt{N}}\sum_{i=1}^N\Psi_k(X_i,\hat{q}_{k,N})\longrightarrow 0
			\end{align*}
			almost surely. This corresponds to equation (27) in \citet{Huber1967}.
			
			We now consider the four assumptions of Theorem 3 from \citet{Huber1967}. Assumptions (N-2) and (N-4) hold by hypothesis. As for assumption (N-1), letting the norm of interest in the original proof of \citet{Huber1967} be the standard $L_2$ norm $\lVert\cdot\rVert_2$, separability is only required to ensure that $u_{k,j}(X,q,d)$ and $\sup_{\tau:\lVert \tau-q_k\rVert_j\leq d_0}Z_{k,j}(X_1,\ldots,X_N,\tau,q)$, where $d>0$, $d_0>0$ and $q$ satisfy $\lVert q-q_k\rVert_j+d\leq d_0$, are random variables for $j=2$ and $j=\infty$ (note that $L_2$ and $L_\infty$ norms are equivalent and if (N-3), (N-4) and Lemma 3 in \citet{Huber1967} hold for any given norm, they hold for all equivalent norms). Lemma \ref{meas}(a) shows that this is the case for $u_{k,j}(X,q,d)$. As for $\sup_{\tau:\lVert \tau-q_k\rVert_j\leq d_0}Z_{k,j}(X_1,\ldots,X_N,\tau,q)$, since $\lVert\lambda_k(\tau)-\lambda_k(q)\rVert_j\leq E[\lVert \Psi_k(X,\tau)-\Psi_k(X,q)\rVert_j]\leq E[u_{k,j}(X,q,d)]$ for $\tau$ satisfying $\lVert \tau-q\rVert_j\leq d$, the continuity of $\lambda_k$ in $\{\tau:\lVert\tau-q_k\rVert_j\leq d_0\}$, and hence, the conclusion of Lemma \ref{meas}(b), follows from condition (ii) of assumption (N-3) under the Euclidean norm. Thus, the conclusion of Theorem 3 of \citet{Huber1967} follows once we demonstrate the three conditions in assumption (N-3). To do so, we extensively use the fact that for any matrix $A$ and vector $v$ that are conformable,
			\begin{equation}\begin{aligned} \label{Av}
			\lVert Av\rVert_2\leq \lVert A\rVert_F\lVert v\rVert_2,
			\end{aligned}\end{equation}
			which follows from the Cauchy-Schwarz inequality and the fact that each component of $Av$ is the dot product of $v$ and a row in $A$. Finally, $T_k:M\times M\rightarrow \mathbb{R}$ defined so that $T_k(x,q)$ is the Frobenius norm of the $n\times n$ Euclidean Hessian matrix of $\langle\beta\xi_q,\log_q(x)\rangle$ as a function of $q$, is continuous on all of $M\times M$ by the given $C^2$ condition, and
			\begin{align*}
			\lvert D_{r'}\Psi_k^r(x,q)\rvert&\leq\frac{S(x,q)}{\lVert x-q\rVert_2}+T_k(x,q),
			\end{align*}
			for all $r,r'=1,\ldots,n$ by Lemma \ref{lem}(b).
			
			Note that for any of the three conditions, if it holds for a given $d_0>0$, it will also hold for any smaller positive value, so it is sufficient to show the three conditions hold separately for some values of $d_0$ that are not necessarily equal.
			
			(i):  If (III) is true for a certain value of $d_1$, it is also true for any smaller positive value of $d_1$. Therefore, assume without loss of generality that $\{q:\lVert q-q_k\rVert_2\leq 2d_1\}$ is contained in $Q_k\cap U$. Fix $r$ and $r'$. Then, (III), Lemma \ref{lem}, and the boundedness of the density, which we will call $f$, on $Q_k$, imply
			\begin{align*}
			B_1:&=E\bigg[\sup_{q:\lVert q-q_k\rVert\leq d_1}\lvert D_{r'}\Psi_k^r(X,q)\rvert^{1.5};X\not\in \bar{U}\bigg]<\infty, \\
			B_2:&=\sup_{(x,q):\lVert x-q_k\rVert_2\geq 2d_1, x\in \bar{U},\lVert q-q_k\rVert_2\leq d_1} \lvert D_{r'}\Psi_k^r(x,q)\rvert<\infty, \\
			B_3:&=\sup_{(x,q):\lVert x-q_k\rVert_2\leq 2d_1,\lVert q-q_k\rVert_2\leq d_1} S(x,q)<\infty, \\
			B_4:&=\sup_{(x,q):\lVert x-q_k\rVert_2\leq 2d_1,\lVert q-q_k\rVert_2\leq d_1} T_k(x,q)<\infty, \\
			B_5:&=\sup_{x:\lVert x-q_k\rVert_2\leq 2d_1} f(x)<\infty,
			\end{align*}
			since upper semi-continuous functions attain their maxima on compact sets. For non-negative real numbers $a$ and $b$, $(a+b)^{1.5}=(a+b)\sqrt{a+b}\leq(a+b)(\sqrt{a}+\sqrt{b})=a^{1.5}+ab^{0.5}+a^{0.5}b+b^{1.5}$; so 
			\begin{align*}
			&\bigg(\frac{S(x,q)}{\lVert X-q\rVert_2}+T_k(x,q)\bigg)^{1.5} \\
   &\leq\frac{S(x,q)^{1.5}}{\lVert X-q\rVert_2^{1.5}}+\frac{S(x,q)T_k(x,q)^{0.5}}{\lVert X-q\rVert_2}+\frac{S(x,q)^{0.5}T_k(x,q)}{\lVert X-q\rVert_2^{0.5}}+T_k(x,q)^{1.5}
			\end{align*}
			when $x\neq q$, which we use below. Let us first prove the existence and continuity of the function $q\mapsto E[D_{r'}\Psi_k^r(X,q)]$ on $\{q:\lVert q-q_k\rVert_2\leq d_1\}$. For any $q$ in this set,
			\begin{equation}\begin{aligned} \label{long}
			&E\big[\lvert D_{r'}\Psi_k^r(X,q)\rvert^{1.5}\big] \\
			&=E\big[\lvert D_{r'}\Psi_k^r(X,q)\rvert^{1.5};X\not\in \bar{U}\big] \\
			&\qquad+E\big[\lvert D_{r'}\Psi_k^r(X,q)\rvert^{1.5};\lVert X-q_k\rVert_2\geq 2d_1, X\in \bar{U}\big] \\
			&\qquad+E\big[\lvert D_{r'}\Psi_k^r(X,q)\rvert^{1.5};\lVert X-q_k\rVert_2< 2d_1\big] \\
			&\leq E\bigg[\sup_{q:\lVert q-q_k\rVert\leq d_1}\lvert D_{r'}\Psi_k^r(X,q)\rvert^{1.5};X \not\in \bar{U}\bigg] \\
			&\qquad+E\bigg[\sup_{q:\lVert q-q_k\rVert\leq d_1}\lvert D_{r'}\Psi_k^r(X,q)\rvert^{1.5};\lVert X-q_k\rVert_2\geq 2d_1, X\in \bar{U}\bigg] \\
			&\qquad+E\bigg[\bigg(\frac{S(X,q)}{\lVert X-q\rVert_2}+T_k(X,q)\bigg)^{1.5};\lVert X-q_k\rVert_2\leq 2d_1\bigg] \\
			&\leq B_1+B_2^{1.5}+B_3^{1.5}E\bigg[\frac{1}{\lVert X-q\rVert_2^{1.5}};\lVert X-q_k\rVert_2\leq 2d_1\bigg] \\ 
			&\qquad+B_3B_4^{0.5}E\bigg[\frac{1}{\lVert X-q\rVert_2};\lVert X-q_k\rVert_2\leq 2d_1\bigg] \\
			&\qquad+B_3^{0.5}B_4E\bigg[\frac{1}{\lVert X-q\rVert_2^{0.5}};\lVert X-q_k\rVert_2\leq 2d_1\bigg]+B_4^{1.5}.
			\end{aligned}\end{equation}
			Now, for $\nu>0$,
			\begin{equation}\begin{aligned} \label{longer}
			E\bigg[\frac{1}{\lVert X-q\rVert_2^\nu};\lVert X-q_k\rVert_2\leq 2d_1\bigg]&\leq B_5\int I(\lVert x-q_k\rVert_2\leq 2d_1)\frac{1}{\lVert x-q\rVert_2^\nu} dx \\
			&=B_5\int I(\lVert y+q-q_k\rVert_2\leq 2d_1)\frac{1}{\lVert y\rVert_2^\nu} dy \\
			&\leq B_5\int I(\lVert y\rVert_2\leq 3d_1)\frac{1}{\lVert y\rVert_2^\nu} dy.
			\end{aligned}\end{equation}
			In the standard $n$-dimensional spherical coordinate substitution of $y=(y^1,\ldots,y^n)^T$ with $(R,\theta_1,\ldots,\theta_{n-1})^T$, $R=\lVert y\rVert_2$, $y^1=R\sin(\theta_1)\cdots\sin(\theta_{n-1})$ and $y^r=R\sin(\theta_1)\cdots\sin(\theta_{n-r})\cos(\theta_{n-r+1})$ for $r=2,\ldots,n$, giving $dy=dy^1\cdots dy^n=R^{n-1}\sin^{n-2}(\theta_1)\cdots\sin(\theta_{n-2})dRd\theta_1\ldots d\theta_{n-1}$. Therefore,
			\begin{equation}\begin{aligned} \label{short1}
			&\int I(\lVert y\rVert_2\leq 3d_1)\frac{1}{\lVert y\rVert_2^\nu} dy  \\
			&=\int_0^{2\pi}\int_0^\pi\ldots\int_0^\pi\int_0^{3d_1} \frac{R^{n-1}\sin^{n-2}(\theta_1)\ldots\sin(\theta_{n-2})}{R^\nu}dRd\theta_1\ldots d\theta_{n-1}  \\
			&\leq \int_0^{2\pi}\int_0^\pi\ldots\int_0^\pi\int_0^{3d_1} \frac{R^{n-1}}{R^\nu}dRd\theta_1\ldots d\theta_{n-1}  \\
   %&=\frac{(2d_1)^{n-\nu}}{n-\nu}\times \int_0^\pi \sin(\theta_1)d\theta_1\ldots\int_0^\pi\sin(\theta_{n-2})d\theta_{n-2}\times 2\pi  \\
			&=\frac{(3d_1)^{n-\nu}}{n-\nu}\times \pi^{n-2}\times 2\pi \\
			&<\infty
			\end{aligned}\end{equation}
			if $n>\nu$. Combining (\ref{long}), (\ref{longer}), and (\ref{short1}), we see that $E[\lvert D_{r'}\Psi_k^r(X,q)\rvert^{1.5}]$ is finite for all $q$ satisfying $\lVert q-q_k\rVert_2\leq d_1$, and moreover it is bounded above by a quantity independent of $q$. Therefore all elements of $\{D_{r'}\Psi_k^r(X,q)\}_{q:\lVert q-q_k\rVert_2\leq d_1}$ are integrable and this family is $L^{1.5}$-bounded, and therefore uniformly integrable. Additionally, for all $x\neq q$, $D_{r'}\Psi_k^r(x,\tau)\rightarrow D_{r'}\Psi_k^r(x,q)$ as $\tau\rightarrow q$ by the continuity of $q\mapsto D_{r'}\Psi_k^r(x,q)$ at $q\neq x$, and therefore $D_{r'}\Psi_k^r(X,\tau)\rightarrow D_{r'}\Psi_k^r(X,q)$ almost surely for $q\in \{q:\lVert q-q_k\rVert_2\leq d_1\}\subset Q_k$. With uniform integrability, this implies that the map $q\mapsto E[D_{r'}\Psi_k^r(X,q)]$ is continuous on $\{q:\lVert q-q_k\rVert_2\leq d_1\}$ by Vitali's convergence theorem, as desired; this is true for all $r,r'=1,\ldots,n$. 
			
			Notice in particular that we have also proven $D_{r'}\Psi_k^r(X,q_k)$ is integrable, as claimed in the statement of the theorem. Since the $n\times n$ matrix $D\Psi_k(x,q_k)$ is symmetric, $E[D\Psi_k(X,q_k)]$ is symmetric and, by (IV), nonsingular. Therefore, it has real, non-zero eigenvalues. Call the eigenvalue with the smallest absolute value $\zeta$. By the conclusion of the above paragraph, there exists some $d_2\in (0,d_1)$ for which
			\begin{equation}\begin{aligned} \label{unif}
			\sum_{r'=1}^n\sum_{r=1}^n \sup_{\tau:\lVert \tau-q_k\rVert_2\leq d_2}\Big\lvert E[D_{r'}\Psi_k^r(X,\tau)]-E[D_{r'}\Psi_k^r(X,q_k)]\Big\rvert<\frac{\lvert\zeta\rvert}{2}.
			\end{aligned}\end{equation}
			
			For any $q$ satisfying $\lVert q-q_k\rVert_2\leq d_1$, denote by $c_{q_k,q}$ the straight line segment in Euclidean space for which $c_{q_k,q}(0)=q_k$ and $c_{q_k,q}(1)=q$. Notice that the same argument used to prove the continuity of $q\mapsto E[D_{r'}\Psi_k^r(X,q)]$ on $\{q:\lVert q-q_k\rVert_2\leq d_1\}$ also shows the continuity of $q\mapsto E[\lvert D_{r'}\Psi_k^r(X,q)\rvert]$ on this set. Thus, 
			\begin{align*}
			\int_0^1 E[\lvert D_{r'}\Psi_k^r(X,c_{q_k,q}(t))\rvert]dt\leq \sup_{\tau:\lVert\tau-q_k\rVert_2\leq d_1}E[\lvert D_{r'}\Psi_k^r(X,\tau)\rvert]<\infty.
			\end{align*}
			Since $D_{r'}\Psi_k^r$ is 0 on $\{(x,p):x=p\}$ and continuous elsewhere in $\phi(M)\times\phi(M)$, it is measurable, so $(\omega,t)\mapsto D_{r'}\Psi_k^r(X(\omega),c_{q_k,q}(t))$ is too, and we can use Tonelli's and Fubini's theorems:
\begin{align}\label{fubini}
			\int_0^1 E[D_{r'}\Psi_k^r(X,c_{q_k,q}(t))]dt=E\Big[\int_0^1 D_{r'}\Psi_k^r(X,c_{q_k,q}(t))dt\Big].
			\end{align}		
			
			For any $r=1,\ldots,n$ and $\tau$, denote by $\lambda_k^r(\tau)$ and $\Psi_k^r(X,\tau)$ the $r$th elements of $\lambda_k(\tau)$ and $\Psi_k(X,\tau)$, respectively. Then, 
			\begin{equation}\begin{aligned}\label{ftc}
			&\lambda_k^r(q)-\lambda_k^r(q_k) \\
			&=E[\Psi_k^r(X,c_{q_k,q}(1))-\Psi_k^r(X,c_{q_k,q}(0))] \\
			&=E\bigg[\int_0^1 [D_1\Psi_k^r(X,c_{q_k,q}(t)),\ldots,D_n\Psi_k^r(X,c_{q_k,q}(t))]c_{q_k,q}'(t)dt\bigg] \\
			&=\bigg[E\Big[\int_0^1D_1\Psi_k^r(X,c_{q_k,q}(t))dt\Big],\ldots, E\Big[\int_0^1D_n\Psi_k^r(X,c_{q_k,q}(t))dt\Big]\bigg](q-q_k) \\
			&=\bigg[\int_0^1E[D_1\Psi_k^r(X,c_{q_k,q}(t))]dt,\ldots, \int_0^1E[D_n\Psi_k^r(X,c_{q_k,q}(t))]dt\bigg](q-q_k).
			\end{aligned}\end{equation}
			The second equality follows from the fundamental theorem of calculus, the chain rule, and the continuity and differentiability of $t\mapsto \Psi_k^r(x,c_{q_k,q}(t))$ on $[0,1]$ and $(0,1)$, respectively, for all $x$ outside the image of $c_{q_k,q}$ (the probability of this image containing $X$ is 0). The last equality follows from (\ref{fubini}).
						
			Then, (\ref{ftc}), (\ref{Av}), the dominance of the $L_1$ norm over the $L_2$ norm, and the continuity of $q\mapsto E[D_{r'}\Psi_k^r(X,q)]$ at $q=q_k$ imply that
			\begin{equation}\begin{aligned} \label{hess}
			&\lVert \lambda_k(q)-\lambda_k(q_k)-E[D\Psi_k(X,q_k)](q-q_k) \rVert_2/\lVert q-q_k\Vert_2 \\
			&=\bigg\lVert \bigg[\int_0^1E[D_{r'}\Psi_k^r(X,c_{q_k,q}(t))]-E[D_{r'}\Psi_k^r(X,q_k)]dt\bigg]_{r,r'=1,\ldots,n} \\
			&\qquad \times (q-q_k) \bigg\rVert_2/\lVert q-q_k\Vert_2\\
			&\leq \bigg\lVert \bigg[\int_0^1E[D_{r'}\Psi_k^r(X,c_{q_k,q}(t))]-E[D_{r'}\Psi_k^r(X,q_k)]dt\bigg]_{r,r'=1,\ldots,n}\bigg\rVert_F\\
			&\leq \sum_{r'=1}^n\sum_{r=1}^n\bigg\lvert\int_0^1E[D_{r'}\Psi_k^r(X,c_{q_k,q}(t))]-E[D_{r'}\Psi_k^r(X,q_k)]dt\bigg\rvert\\
			&\leq \sum_{r'=1}^n\sum_{r=1}^n\int_0^1\Big\lvert E[D_{r'}\Psi_k^r(X,c_{q_k,q}(t))]-E[D_{r'}\Psi_k^r(X,q_k)]\Big\rvert dt\\
			&\leq \sum_{r'=1}^n\sum_{r=1}^n\sup_{\tau:\lVert \tau-q_k\rVert_2\leq\lVert q-q_k\rVert_2}\Big\lvert E[D_{r'}\Psi_k^r(X,\tau)]-E[D_{r'}\Psi_k^r(X,q_k)]\Big\rvert\\
			&\rightarrow 0
			\end{aligned}\end{equation}
			as $q\rightarrow q_k$. Together with (\ref{unif}), this means that $\lVert q-q_k\rVert_2<d_2$ implies
			\begin{align*}
			&\lVert \lambda_k(q)\rVert_2 \\
			&=\lVert  E[D\Psi_k(X,q_k)](q-q_k)+\lambda_k(q)-\lambda_k(q_k)-E[D\Psi_k(X,q_k)](q-q_k)\rVert_2 \\
			&>\big\lVert\big[E[D\Psi_k(X,q_k)]\big](q-q_k) \big\rVert_2-\frac{\lvert \zeta\rvert}{2}\lVert q-q_k\Vert_2 \\
			&\geq\lvert \zeta\rvert\lVert q-q_k\Vert_2-\frac{\lvert \zeta\rvert}{2}\lVert q-q_k\Vert_2 \\
			&=\frac{\lvert \zeta\rvert}{2}\lVert q-q_k\Vert_2,
			\end{align*}
			proving that (i) holds with $d_2$ in the role of $d_0$.
			
			(ii): As in the proof for (i), $d_1$ can be arbitrarily close to 0. Therefore, assume without loss of generality that $d_1<1$, and, as before, that $\{q:\lVert q-q_k\rVert_2\leq 2d_1\}$ is contained in $Q_k\cap U$. Then, fixing $q$ and $d\geq 0$ so that $\lVert q-q_k\rVert_2+d\leq d_1$ is satisfied, $\lVert x-q\rVert_2< 2d$ implies $\lVert x-q_k\rVert_2< 2d_1$.
			\begin{equation}\begin{aligned} \label{sum}
			&E[u_{k,2}(X,q,d)] \\
			&=E[u_{k,2}(X,q,d);\lVert X-q\rVert_2 < 2d] \\
			&\qquad+E[u_{k,2}(X,q,d);\lVert X-q\rVert_2\geq 2d, \lVert X-q_k\rVert_2<2d_1] \\
			&\qquad+E[u_{k,2}(X,q,d);\lVert X-q_k\rVert_2\geq 2d_1,X\in \bar{U}]+E[u_{k,2}(X,q,d);X\not\in \bar{U}].
			\end{aligned}\end{equation}
			Noting that when $\tau\neq x$, $\Psi_k(x,\tau)=-g_\tau\log_\tau(x)/\lVert\log_\tau(x)\rVert+T_k^*(x,\tau)$ for some continuous function $T_k^*:M\times M\rightarrow \mathbb{R}^n$, defining $\zeta_{\min}(\tau)$ and $\zeta_{\max}(\tau)$ as the smallest and largest eigenvalues, respectively, of $g_\tau$ and recalling (\ref{int1}) and (\ref{int2}) in the proof of Lemma \ref{lem},
			\begin{equation}\begin{aligned} \label{ref}
			\lVert\Psi_k(x,\tau)\rVert_2\leq\frac{\lVert g_\tau\log_\tau(x)\rVert_2}{\lVert\log_\tau(x)\rVert}+\lVert T_k^*(x,\tau)\rVert_2\leq\frac{\zeta_{\max}(\tau)}{\sqrt{\zeta_{\min}(\tau)}}+\lVert T_k^*(x,\tau)\rVert_2
			\end{aligned}\end{equation}
			when $x\neq\tau$, so
			\begin{align*}
			&u_{k,2}(x,q,d) \\
			&\leq 2\sup_{\tau:\lVert \tau-q_k\rVert_2\leq d_1}\lVert\Psi_k(x,\tau)\rVert_2 \\
			&\leq 2\max\bigg\{\sup_{\tau:\lVert \tau-q_k\rVert_2\leq d_1}\bigg(\frac{\zeta_{\max}(\tau)}{\sqrt{\zeta_{\min}(\tau)}}+\lVert T_k^*(x,\tau)\rVert_2\bigg),\sup_{\tau:\lVert \tau-q_k\rVert_2\leq d_1}\lVert g_\tau(\xi_\tau)_k\rVert_2 \bigg\}.
			\end{align*}
			In the last line, each of the terms whose supremum is being taken is continuous as a function of $(x,\tau)$, and therefore, attains some maximum on $\{(x,\tau):x\in \bar{U},\lVert \tau-q_k\rVert_2\leq d_1\}$ that does not depend on $q$ or $d$; so $u_{k,2}(x,q,d)\leq C_1$ on $\{x:x\in \bar{U}\}$ for some finite $C_1$. Then, letting $C_2$ be the constant for which a ball of radius $R$ in $\mathbb{R}^n$ has volume $C_2R^n$ and recalling the definition of $B_5$ at the beginning of the proof for (i),
			\begin{equation}\begin{aligned} \label{sum1}
			E[u_{k,2}(X,q,d);\lVert X-q\rVert_2 < 2d]\leq B_5C_1C_2(2d)^n\leq 2^nB_5C_1C_2d
			\end{aligned}\end{equation}
   since $d\leq d_1<1$.
			
			Denote by $c_{q,\tau}:[0,1]\rightarrow M$ the straight line segment in Euclidean space for which $c_{q,\tau}(0)=q$ and $c_{q,\tau}(1)=\tau$. If $x$ does not lie on the image of $c_{q,\tau}$, $(\Psi_k(x,\tau)-\Psi_k(x,q))^T\Psi_k(x,c_{q,\tau}(t))$ is continuously differentiable as a function of $t$, and therefore, by the mean value theorem there exists some $t'\in(0,1)$ for which 
			\begin{equation}\begin{aligned} \label{help}
			&\lVert \Psi_k(x,\tau)-\Psi_k(x,q)\rVert_2^2 \\
			&=\big(\Psi_k(x,\tau)-\Psi_k(x,q)\big)^T\Psi_k(x,c_{q,\tau}(1)) \\
			&\qquad-\big(\Psi_k(x,\tau)-\Psi_k(x,q)\big)^T\Psi_k(x,c_{q,\tau}(0)) \\
			&=\big(\Psi_k(x,\tau)-\Psi_k(x,q)\big)^T(D\Psi_k(x,c_{q,\tau}(t'))(\tau-q)) \\
			&\leq\lVert \Psi_k(x,\tau)-\Psi_k(x,q)\rVert_2\lVert\tau-q\rVert_2\lVert D\Psi_k(x,c_{q,\tau}(t'))\rVert_F \\
			&\leq\lVert \Psi_k(x,\tau)-\Psi_k(x,q)\rVert_2\lVert\tau-q\rVert_2\sum_{r'=1}^n\sum_{r=1}^n\lvert D_{r'}\Psi_k^r(x,c_{q,\tau}(t'))\rvert.
			\end{aligned}\end{equation}
			Here, the second-to-last line follows from the Cauchy-Schwarz inequality and (\ref{Av}), and the last from the fact that the $L_2$ norm is less than or equal to the $L_1$ norm. Therefore, dividing through by $\lVert \Psi_k(x,\tau)-\Psi_k(x,q)\rVert_2$, if $\lVert x-q\rVert_2>d$, 
			\begin{equation}\begin{aligned} \label{uxqd}
			u_{k,2}(x,q,d)&\leq \sup_{\tau,\tau':\lVert \tau-q\rVert_2\leq d,\lVert \tau'-q\rVert_2\leq d}\lVert\tau-q\rVert_2\sum_{r'=1}^n\sum_{r=1}^n\lvert D_{r'}\Psi_k^r(x,\tau')\rvert \\
			&\leq d \sum_{r'=1}^n\sum_{r=1}^n\sup_{\tau':\lVert \tau'-q\rVert_2\leq d}\lvert D_{r'}\Psi_k^r(x,\tau')\rvert  \\
			&\leq d \sum_{r'=1}^n\sum_{r=1}^n\sup_{\tau':\lVert \tau'-q_k\rVert_2\leq d_1}\lvert D_{r'}\Psi_k^r(x,\tau')\rvert.
			\end{aligned}\end{equation}
			
			Using the second line in (\ref{uxqd}), the constant $B_5$, Lemma \ref{lem}, constants 
			\begin{align*}
			C_3&:=\sup_{(x,\tau'):\lVert x-q_k\rVert_2\leq 2d_1,\lVert \tau'-q_k\rVert_2\leq d_1}S(x,\tau')<\infty \\
			C_4&:=\sup_{(x,\tau'):\lVert x-q_k\rVert_2\leq 2d_1,\lVert \tau'-q_k\rVert_2\leq d_1}T_k(x,\tau')<\infty,
			\end{align*}
			and two substitutions, first of $x$ with $y=(y^1,\ldots,y^n)^T$ achieved by shifting the space by $-q$, and then, the standard $n$-dimensional spherical coordinate substitution of $y$ with $(R,\theta_1,\ldots,\theta_{n-1})^T$, detailed in the proof for (i), which gives \\$dy^1\cdots dy^n=R^{n-1}\sin^{n-2}(\theta_1)\cdots\sin(\theta_{n-2})dRd\theta_1\cdots d\theta_{n-1}$,
			\begin{equation}\begin{aligned} \label{sum2a}
			&E[u_{k,2}(X,q,d);\lVert X-q\rVert_2\geq 2d, \lVert X-q_k\rVert_2< 2d_1] \\
			%\leq& E[u_{k,2}(X,q,d);\lVert X-q\rVert\geq 2d, \lVert X-q_k\rVert\leq d_2] \\
      &\leq  d\sum_{r'=1}^n\sum_{r=1}^n E\bigg[\sup_{\lVert \tau'-q\rVert_2\leq d}\lvert D_{r'}\Psi_k^r(X,\tau')\rvert;\lVert X-q\rVert_2\geq 2d,\lVert X-q_k\rVert_2< 2d_1\bigg]  \\
			&\leq  B_5d\sum_{r'=1}^n\sum_{r=1}^n\int I(\lVert x-q\rVert_2\geq 2d, \lVert x-q_k\rVert_2< 2d_1) \\
			&\qquad\times\sup_{\lVert \tau'-q\rVert_2\leq d}\lvert D_{r'}\Psi_k^r(x,\tau')\rvert dx  \\
			&\leq n^2B_5d\int I(2d\leq\lVert x-q\rVert_2< 4d_1)\sup_{\lVert \tau'-q\rVert_2\leq d}\bigg(\frac{C_3}{\lVert x-\tau'\rVert_2}+C_4\bigg) dx  \\
			&\leq n^2B_5\bigg(C_2C_4(4d_1)^n \\
			&\qquad+C_3\int I(2d\leq\lVert x-q\rVert_2< 4d_1)\sup_{\lVert \tau'-q\rVert_2\leq d}\frac{1}{\lVert x-\tau'\rVert_2} dx\bigg)d \\
			&=n^2B_5\bigg(4^nC_2C_4d_1^n+C_3\int I(2d\leq\lVert y\rVert_2< 4d_1)\sup_{\lVert \tau'\rVert_2\leq d}\frac{1}{\lVert y-\tau'\rVert_2} dy\bigg)d  \\
			&=n^2B_5\bigg(4^nC_2C_4d_1^n+C_3\int_0^{2\pi}\int_0^\pi\ldots \\
			&\qquad \ldots\int_0^\pi\int_{2d}^{4d_1} \frac{R^{n-1}\sin^{n-2}(\theta_1)\ldots\sin(\theta_{n-2})}{R-d} dRd\theta_1\ldots d\theta_{n-1}\bigg)d  \\
&\leq n^2B_5\bigg(4^nC_2C_4d_1^n+C_3\int_0^{2\pi}\int_0^\pi\ldots\int_0^\pi\int_{2d}^{4d_1} \frac{R^{n-1}}{R-d} dRd\theta_1\ldots d\theta_{n-1}\bigg)d  \\
			&= n^2B_5\bigg(4^nC_2C_4d_1^n+C_3(2\pi)\pi^{n-2}\int_{2d}^{4d_1} \frac{R^{n-1}}{R-d} dR\bigg)d.
			\end{aligned}\end{equation}
			Now, if $d=0$, the integral in the last line becomes 
			\begin{equation}\begin{aligned} \label{dzero}
			\int_0^{4d_1}R^{n-2}dR=\frac{(4d_1)^{n-1}}{n-1}<\infty;
			\end{aligned}\end{equation}
			otherwise, it becomes 
			\begin{equation}\begin{aligned}
			&\int_{2d}^{4d_1} \frac{(R-d+d)^{n-1}}{R-d} dR \\
			&=\int_{2d}^{4d_1}\sum_{j=0}^{n-1}{\binom{n-1}{j}}d^{n-j-1}(R-d)^{j-1} dR  \\
			&=d^{n-1}(\log(4d_1-d)-\log d)+\sum_{j=1}^{n-1}{\binom{n-1}{j}}\frac{d^{n-j-1}(4d_1-d)^j-d^{n-1}}{j}. 
			\end{aligned}\end{equation}
			This is continuous in $d$ when $d>0$, and it converges to $(4d_1)^{n-1}/(n-1)$ when $d\rightarrow 0$ since $\lim_{d\rightarrow 0}d^{n-1}\log d=0$, coinciding with (\ref{dzero}). Therefore, the integral in the last line of (\ref{sum2a}) is continuous as a function of $(q,d)$ (in fact, it does not depend on $q$ at all), and since $\{(q,d):\lVert q-q_k\rVert_2+d\leq d_1,d\geq 0\}$ is compact, this integral attains its maximum (call it $C_5<\infty$) on this set. So, (\ref{sum2a}) gives
			\begin{equation}\begin{aligned} \label{sum2}
			&E[u_{k,2}(X,q,d);\lVert X-q\rVert_2\geq 2d, \lVert X-q_k\rVert_2<2d_1] \\
			&\leq n^2B_5(4^nC_2C_4d_1^n+2C_3C_5\pi^{n-1})d.
			\end{aligned}\end{equation}
			
			Using (\ref{uxqd}) and constants 
			\begin{align*}
			C_6&:=\max_{r,r':1\leq r,r'\leq n}\sup_{(x,\tau'):\lVert x-q_k\rVert_2\geq 2d_1,x\in \bar{U},\lVert \tau'-q_k\rVert_2\leq d_1} \lvert D_{r'}\Psi_k^r(x,\tau')\rvert<\infty, \\
			C_7:&=\max_{r,r':1\leq r,r'\leq n}E\bigg[\sup_{\tau':\lVert \tau'-q_k\rVert_2\leq d_1}\lvert D_{r'}\Psi_k^r(X,\tau')\rvert;X\not\in\bar{U}\bigg],
			\end{align*}
			gives
			\begin{equation}\begin{aligned} \label{sum3}
			&E[u_{k,2}(X,q,d);\lVert X-q_k\rVert_2\geq 2d_1,X\in \bar{U}] \\
			&\leq d\sum_{r'=1}^n\sum_{r=1}^nE\bigg[\sup_{\tau':\lVert \tau'-q_k\rVert_2\leq d_1}\lvert D_{r'}\Psi_k^r(X,\tau')\rvert;\lVert X-q_k\rVert_2\geq 2d_1,X\in \bar{U}\bigg] \\
			&\leq n^2C_6d
			\end{aligned}\end{equation}
			and
			\begin{equation}\begin{aligned} \label{sum4}
			E[u_{k,2}(X,q,d);X\not\in \bar{U}] &\leq d\sum_{r'=1}^n\sum_{r=1}^nE\bigg[\sup_{\tau':\lVert \tau'-q_k\rVert_2\leq d_1}\lvert D_{r'}\Psi_k^r(X,\tau')\rvert;X\not\in \bar{U}\bigg] \\
			&\leq n^2C_7d.
			\end{aligned}\end{equation}
			
			The results (\ref{sum}), (\ref{sum1}), (\ref{sum2}), (\ref{sum3}), and (\ref{sum4}) show that (ii) holds with $d_1$ in the role of $d_0$. %An examination of (\ref{sum2a}), (\ref{dzero}), (\ref{sum3}) and (\ref{sum4}) also reveals that, by letting $d\rightarrow 0$, each $E[\lvert D_{r'}\Psi_k^r(X,q)\rvert]\leq B_5(4^nC_2C_4d_1^n+2C_3\pi^{n-1}(4d_1)^{n-1}/(n-1))+C_6+C_7$, and so is finite for all $q$ that satisfy $\lVert q-q_k\rVert_2\leq d_1$, in particular $q=q_k$.
			
			(iii): This case is similar to that of (ii). Make the same assumptions about $d_1$ as in (ii), and fix $q$ and $d\geq 0$ so that $\lVert q-q_k\rVert_2+d\leq d_1$.
			\begin{equation}\begin{aligned} \label{summ}
			&E[u_{k,2}(X,q,d)^2] \\
			&=E[u_{k,2}(X,q,d)^2;\lVert X-q\rVert_2 < 2d] \\
			&\qquad+E[u_{k,2}(X,q,d)^2;\lVert X-q\rVert_2\geq 2d, \lVert X-q_k\rVert_2<2d_1] \\
			&\qquad+E[u_{k,2}(X,q,d)^2;\lVert X-q_k\rVert_2\geq 2d_1,X\in \bar{U}] \\
			&\qquad+E[u_{k,2}(X,q,d)^2;X\not\in \bar{U}].
			\end{aligned}\end{equation}
			Here, we deal with the second of the four summands above last. Using the constants defined in the proofs for (i) and (ii),
			\begin{equation}\begin{aligned} \label{summ1}
			E[u_{k,2}(X,q,d)^2;\lVert X-q\rVert_2 < 2d]\leq B_5C_1^2C_2(2d)^n\leq 2^nB_5C_1^2C_2d.
			\end{aligned}\end{equation}
			Dividing the second-to-last line of (\ref{help}) by $\lVert \Psi_k(x,\tau)-\Psi_k(x,q)\rVert_2$ and squaring, if $\lVert x-q\rVert_2>d$, 
			\begin{equation}\begin{aligned} \label{uxqd2}
			u_{k,2}(x,q,d)^2&\leq \sup_{\tau,\tau':\lVert \tau-q\rVert_2\leq d,\lVert \tau'-q\rVert_2\leq d}\lVert\tau-q\rVert_2^2\sum_{r'=1}^n\sum_{r=1}^n\lvert D_{r'}\Psi_k^r(x,\tau')\rvert^2 \\
   &\leq d^2 \sum_{r'=1}^n\sum_{r=1}^n\sup_{\tau':\lVert \tau'-q\rVert_2\leq d}\lvert D_{r'}\Psi_k^r(x,\tau')\rvert^2 \\
			&\leq d \sum_{r'=1}^n\sum_{r=1}^n\sup_{\tau':\lVert \tau'-q_k\rVert_2\leq d_1}\lvert D_{r'}\Psi_k^r(x,\tau')\rvert^2,
			\end{aligned}\end{equation}
			so, 
			\begin{equation}\begin{aligned} \label{summ2}
			&E[u_{k,2}(X,q,d)^2;\lVert X-q_k\rVert_2\geq 2d_1,X\in\bar{U}] \\
			& \leq d\sum_{r'=1}^n\sum_{r=1}^nE\bigg[\sup_{\tau':\lVert \tau'-q_k\rVert_2\leq d_1}\lvert D_{r'}\Psi_k^r(X,\tau')\rvert^2;\lVert X-q_k\rVert_2\geq 2d_1,X\in\bar{U}\bigg] \\
			& \leq n^2C_6^2d,
			\end{aligned}\end{equation}
			and defining a constant
			\begin{align*}
			C_8:&=\max_{r,r':1\leq r,r'\leq n}E\bigg[\sup_{\tau':\lVert \tau'-q_k\rVert_2\leq d_1}\lvert D_{r'}\Psi_k^r(X,\tau')\rvert^2;X\not\in\bar{U}\bigg]
			\end{align*}
			that is finite by (III),
			\begin{equation}\begin{aligned} \label{summ3}
			&E[u_{k,2}(X,q,d)^2;X\not\in\bar{U}] \\
			&\leq d\sum_{r'=1}^n\sum_{r=1}^nE\bigg[\sup_{\tau':\lVert \tau'-q_k\rVert_2\leq d_1}\lvert D_{r'}\Psi_k^r(X,\tau')\rvert^2;X\not\in\bar{U}\bigg] \\
			&\leq n^2C_8d.
			\end{aligned}\end{equation}
			
			For the remaining part, making the same substitutions as in (ii),
			\begin{equation}\begin{aligned} \label{summ4a}
			&E[u_{k,2}(X,q,d)^2;\lVert X-q\rVert_2\geq 2d, \lVert X-q_k\rVert_2< 2d_1] \\
			&\leq  B_5d\sum_{r'=1}^n\sum_{r=1}^n\int I(\lVert x-q\rVert_2\geq 2d, \lVert x-q_k\rVert_2< 2d_1) \\
			&\qquad\times \sup_{\lVert \tau'-q\rVert_2\leq d}\lvert D_{r'}\Psi_k^r(X,\tau')\rvert^2 dx  \\
			&\leq n^2B_5d\int I(2d\leq\lVert x-q\rVert_2< 4d_1)\sup_{\lVert \tau'-q\rVert_2\leq d}\bigg(\frac{C_3^2}{\lVert x-\tau'\rVert_2^2} \\
			&\qquad+\frac{2C_3C_4}{\lVert x-\tau'\rVert_2}+C_4^2\bigg) dx  \\
			&\leq n^2B_5\bigg(C_2C_4^2(4d_1)^n \\
			&\qquad+2C_3C_4\int I(2d\leq\lVert x-q\rVert_2< 4d_1)\sup_{\lVert \tau'-q\rVert_2\leq d}\frac{1}{\lVert x-\tau'\rVert_2} dx \\
			&\qquad+C_3^2\int I(2d\leq\lVert x-q\rVert_2< 4d_1)\sup_{\lVert \tau'-q\rVert_2\leq d}\frac{1}{\lVert x-\tau'\rVert_2^2} dx\bigg)d \\
			&=n^2B_5\bigg(C_2C_4^2(4d_1)^n+2C_3C_4C_5(2\pi)\pi^{n-2}  \\
   &\qquad+C_3^2\int I(2d\leq\lVert y\rVert_2< 4d_1)\sup_{\lVert \tau'\rVert_2\leq d}\frac{1}{\lVert y-\tau'\rVert_2^2} dy\bigg)d \\
			&\leq n^2B_5\bigg(C_2C_4^2(4d_1)^n+2C_3C_4C_5(2\pi)\pi^{n-2} \\
			&\qquad+C_3^2(2\pi)\pi^{n-2}\int_{2d}^{4d_1} \frac{R^{n-1}}{(R-d)^2} dR\bigg)d.
			\end{aligned}\end{equation}
			Assume for now that $n\geq 3$. Then, $d=0$ gives 
			\begin{equation}\begin{aligned} 
			\int_0^{4d_1}R^{n-3}dR=\frac{(4d_1)^{n-2}}{n-2}<\infty; 
			\end{aligned}\end{equation}
			if $d>0$, 
			\begin{equation}\begin{aligned}
			&\int_{2d}^{4d_1} \frac{R^{n-1}}{(R-d)^2} dR \\
			&=\int_{2d}^{4d_1}\sum_{j=0}^{n-1}{\binom{n-1}{j}}d^{n-j-1}(R-d)^{j-2} dR  \\
			&=\bigg(d^{n-2}-\frac{d^{n-1}}{4d_1-d}\bigg)+(n-1)d^{n-2}(\log(4d_1-d)-\log d) \\
			&\qquad+\sum_{j=2}^{n-1}{\binom{n-1}{j}}\frac{d^{n-j-1}(4d_1-d)^{j-1}-d^{n-2}}{j-1}. 
			\end{aligned}\end{equation}
			Again, this is continuous in $d$ when $d>0$, and it converges to $(4d_1)^{n-2}/(n-2)$ when $d\rightarrow 0$. Therefore, by an argument analogous to that for (ii), (\ref{summ4a}) gives
			\begin{equation}\begin{aligned} \label{summ4}
			&E[u_{k,2}(X,q,d)^2;\lVert X-q\rVert_2\geq 2d, \lVert X-q_k\rVert_2<2d_1] \\
			& \leq n^2B_5\bigg(C_2C_4^2(4d_1)^n+4C_3C_4C_5\pi^{n-1}+2C_3^2C_9\pi^{n-1}\bigg)d
			\end{aligned}\end{equation}
			for some finite constant $C_9$. Thus, we have the desired result when $n\geq 3$.
			
			To deal with the case of $n=2$, note that (iii) is a stronger condition than is required; as \citet{Huber1967} mentioned immediately after stating (iii), $E[u_{k,2}(X,q,d)^2]=o(\lvert \log d\rvert^{-1})$ as $d\rightarrow 0$ is sufficient. Since $d=o(\lvert\log d\rvert^{-1})$ as $d\rightarrow 0$ and (\ref{summ}), (\ref{summ1}), (\ref{summ2}), (\ref{summ3}), and (\ref{summ4a}) are still true when $n=2$, it is sufficient to show 
			\begin{align*}
			d\int_{2d}^{4d_1} \frac{R^{n-1}}{(R-d)^2} dR=d\log(4d_1-d)-\frac{d^2}{4d_1-d}-d\log d+d=o(\lvert\log d\rvert^{-1})
			\end{align*}
			as $d\rightarrow 0$, or equivalently,
			\begin{align*}
			\lim_{d\rightarrow 0}\bigg(d\log d\log(4d_1-d)-\frac{d^2\log d }{4d_1-d}-d(\log d)^2+d\log d\bigg)\longrightarrow 0.
			\end{align*}
			It is known that $d\log d\rightarrow 0$ as  $d\rightarrow 0$, and by L'H\^opital's rule,
			\begin{align*}
			\lim_{d\rightarrow 0}d(\log d)^2=\lim_{d\rightarrow 0}\frac{(\log d)^2}{1/d}=\lim_{d\rightarrow 0}\frac{2\log d/d}{-1/d^2}=-2\lim_{d\rightarrow 0}d\log d=0,
			\end{align*}
			so we are done.
			
			Having demonstrated that (i), (ii), and (iii) hold, we apply Theorem 3 in \citet{Huber1967}:
			\begin{align*}
			\frac{1}{\sqrt{N}}\sum_{i=1}^N\Psi_k(X_i,q_k)+\sqrt{N}\lambda_k(\hat{q}_{k,N})\rightarrow 0
			\end{align*}
			in probability. Since this is true for each $k$, we simply concatenate to conclude that
			\begin{align*}
			\frac{1}{\sqrt{N}}\sum_{i=1}^N(\Psi_1(X_i,q_1)^T,\ldots,\Psi_K(X_i,q_K)^T)^T+\sqrt{N}(\lambda_1(\hat{q}_{1,N}),\ldots,\lambda_K(\hat{q}_{K,N}))^T\rightarrow 0
			\end{align*}
			in probability. Finally, (\ref{hess}), in the proof for (i), implies $\lambda_k$ is differentiable at $q_k$, and its derivative is $E[D\Psi_k(X,q_k)]$. Then, by condition (IV), we apply the Corollary in \citet{Huber1967}, and since, as observed earlier, each $E[D\Psi_k(X,q_k)]$ is symmetric, the theorem follows.
		\end{proof}
		
		\begin{proof}[Proof of Corollary \ref{bdd}]
			\sloppy By the boundedness of the support, $\lambda_k(q_k)=g_{q_k} E[\nabla \rho(X,q_k;\beta_k,\xi_k);X\neq q_k]=0$ by Theorem \ref{qgrad}. Referring to (\ref{ref}), the map $x\mapsto \lVert\Psi_k(x,q_k)\rVert_2$ is bounded on the compact support of $X$, so condition (II) in Theorem \ref{clt} holds. Additionally, there exists a sufficiently large neighbourhood $U$ of $q_k$  which is bounded in the Riemannian metric and for which $P(X\not\in \bar{U})=0$, so condition (III) in Theorem \ref{clt} holds too.
		\end{proof}

                \begin{proof}[Proof of Proposition \ref{speed}]
    Letting $\tau=\hat{q}_{k,N}$ and $x=q_k$ in (\ref{mtaylors}),
    \begin{align*}
        \sqrt{N}d(\hat{q}_{k,N},q_k)&=\sqrt{N}\lVert (I-F(\hat{q}_{k,N}))(q_k-\hat{q}_{k,N})\rVert \\
        &\leq\lVert \sqrt{N}(q_k-\hat{q}_{k,N})\rVert+\lVert F(\hat{q}_{k,N})\sqrt{N}(q_k-\hat{q}_{k,N})\rVert.
    \end{align*}
\sloppy On one hand, $\lVert F(\hat{q}_{k,N})\sqrt{N}(q_k-\hat{q}_{k,N})\rVert=(\sqrt{N}(q_k-\hat{q}_{k,N}))^TF(\hat{q}_{k,N})^Tg_{\hat{q}_{k,N}}F(\hat{q}_{k,N})\sqrt{N}(q_k-\hat{q}_{k,N}))^{1/2}\rightarrow 0$ in probability as $N\rightarrow\infty$ by the consistency of $\hat{q}_{k,N}$, the definition of $F$, the continuity of the Riemannian metric, Slutsky's theorem and the continuous mapping theorem; thus it is $o_p(1)$. Similarly, letting $Y$ be the random $n$-vector to which $\sqrt{N}(\hat{q}_{k,N}-q_k)$ converges weakly, $\lVert \sqrt{N}(q_k-\hat{q}_{k,N})\rVert=(\sqrt{N}(q_k-\hat{q}_{k,N}))^Tg_{\hat{q}_{k,N}}\sqrt{N}(q_k-\hat{q}_{k,N}))^{1/2}$ converges weakly to $(Y^Tg_{q_k}Y)^{1/2}$; thus it is $O_p(1)$. Since $O_p(1)+o_p(1)=O_p(1)$, the result follows.
\end{proof}

\subsection{Remarks on Theorem \ref{clt}}\label{remarks}

\begin{remark} \label{lsc}
    We can show that $\sup_{q:\lVert q-q_k\rVert_2\leq d_1}\lvert D_{r'}\Psi_k^r(X,q)\rvert^2 I(X\not\in \bar{U})$ from (III) is indeed a random variable and its expected value is therefore well-defined as follows.
    
    Let $\Theta_1$ and $\Theta_2$ be measurable subsets of $M$, $\Theta_2$ be compact, and $f:M\times M\rightarrow\mathbb{R}$ be a function such that for any fixed $x\in \Theta_1$, $q\mapsto f(x,q)$ is continuous on $\Theta_2$ and for any fixed $q\in \Theta_2$, $x\mapsto f(x,q)$ is continuous on $\Theta_1$. For a fixed $x=x_0\in \Theta_1$, there exists some $q_0\in \Theta_2$ for which $f(x_0,q_0)=\sup_{q:q\in \Theta_2}f(x_0,q)$. Hence $\lim\inf_{x\rightarrow x_0}\sup_{q:q\in \Theta_2}f(x,q)\geq\lim\inf_{x\rightarrow x_0}f(x,q_0)=f(x_0,q_0)=\sup_{q:q\in \Theta_2}f(x_0,q)$; thus, $\sup_{q:q\in \Theta_2}f(x,q)I(x\in \Theta_1)$ is lower semi-continuous as a function of $x$ on $\Theta_1$, making it measurable. Functions of the type $\sup_{q:q\in \Theta_2}f(x,q)$ appear repeatedly in the proof of Theorem \ref{clt}, and we can rest assured that their integrals are well-defined.
\end{remark}

\begin{remark} \label{conditions}
Conditions (I) and (IV) are intuitive. We touch on the uniqueness of quantiles in Remark \ref{unique}, and the boundedness of the density in a neighbourhood of $q_k$ is achieved if, for example, $X$ has a continuous density on $M$. The nonsingularity of each $E[D\Psi_k(X,q_k)]$ is necessary because the limiting distribution is defined using $\Lambda^{-1}$. The rest simply require some expectation to be finite. (II) requires that the norm of the gradient at the quantile has a finite second moment. In (III), $\{q:\lVert q-q_k\rVert_2\leq d_1\}\subset U$ is needed to prevent the expected value from being infinite. Despite its mysterious appearance, (III) is not difficult to achieve; for example, in the Euclidean case using the identity map as the chart, when $p\neq x$ the Hessian matrix $D\Psi_k(x,p)$ is $I_n/\lVert x-p\rVert_2-(x-p)(x-p)^T/\lVert x-p\rVert_2^3$, so letting $U=\{q:\lVert q-q_k\rVert_2\leq d\}$ for any $d_1$ and $d$ such that $d>d_1$,
\begin{align*}
    &E\bigg[\sup_{q:\lVert q-q_k\rVert_2\leq d_1}\lvert D_{r'}\Psi_k^r(X,q)\rvert^2;X\not\in \bar{U}\bigg] \\
    &\leq E\bigg[\sup_{q:\lVert q-q_k\rVert_2\leq d_1}\mathrm{tr}(D\Psi_k(x,q)^TD\Psi_k(x,q));X\not\in \bar{U}\bigg] \\
    &= E\bigg[\sup_{q:\lVert q-q_k\rVert_2\leq d_1}\frac{n-1}{\lVert x-q\rVert_2^2};X\not\in \bar{U}\bigg]\leq\frac{n-1}{d-d_1}<\infty
\end{align*}
for all $k=1,\ldots K$; that is, (III) is attained for any random vector $X$.
\end{remark}

\begin{remark}
     Since the work of \citet{Bhattacharya2005}, there has been much research concerning the asymptotic normality of sample Fr\'echet means on non-Euclidean metric spaces, including that of \citet{Bhattacharya2017}. %, \citet{Mattingly2023a}, \citet{Mattingly2023b}). We highlight two in particular: \citet{Bhattacharya2017} further generalised the work of \citet{Bhattacharya2005} to partially $C^\infty$ metric spaces, and \citet{Eltzner2019} developed a so-called `smeary' central limit theorem for Fr\'echet means that holds under very general conditions on manifolds. 
There are also a few theorems dealing with asymptotic normality for minimisers of more general classes of loss functions, and one may wonder whether the results of this section can be easily proven with any of these. Theorem 6 of \citet{Huckemann2011} concerns loss functions that satisfy a smoothness condition; due to the singularity of $\rho(x,\cdot;\beta,\xi)$ at $p=x$, this theorem only works for quantiles if the population quantile is outside the support of $X$. Theorem \ref{clt} and its corollary, our main results here, do not have this strict condition. \citet{Brunel2023} demonstrated consistency and asymptotic normality for minimisers of convex loss functions on Riemannian manifolds, but as noted in Remark \ref{unique}, these do not apply in general to quantiles on Hadamard manifolds. 

Theorem 2.11 of \citet{Eltzner2019} is general, but it does not seem applicable here. Firstly, Assumption 2.4(ii) of that theorem, in particular the finiteness of $E[\dot \rho(X)^2]$, does not seem to follow from our conditions, without significantly strengthening (II). That theorem also requires $G^{\beta,\xi}$ to be sufficiently smooth in a neighbourhood of $q(\beta,\xi)$. Verifying this, especially when $n=2$, under the conditions of Theorem \ref{clt} is not simple as not only is $\rho(x,\cdot;\beta,\xi)$ not differentiable at $x$, but its second derivative may explode near $x$: In $\mathbb{R}^n$ the Hessian matrix at $p\neq x$ is $\lvert x-p\rvert^{-1}(I-\lvert x-p\rvert^{-2}(x-p)(x-p)^T)$. Thus, our results are not at all trivial corollaries, if corollaries at all, of this theorem.

Another problem is that \citet{Eltzner2019} used $\log_\mu$ as a parameter in their chart; in contrast, this section is valid for arbitrary charts. Use of $\log_\mu$ is not desirable because it requires knowledge of $\mu$, defeating one of the main purposes of asymptotic normality; that is, inference for an unknown parameter.
 \end{remark}

\subsection{Proofs of Theorem \ref{thm5.1}}\label{proof_thm5.1}

\begin{proposition} \label{hyper}
For any $y\in \mathbb{H}^n$ and unit vector in $T_y\mathbb{H}^n$, let $\xi$ be the unique point in $\partial\mathbb{H}^n$ such that this vector equals $\xi_y$. The radial field on $\mathbb{H}^n$ is $C^\infty$ and given by 
\begin{align*}
\xi_p=\frac{y+\xi_y+\langle p,y+\xi_y\rangle_{\textbf{M}} p}{\lVert y+\xi_y+\langle p,y+\xi_y\rangle_{\textbf{M}} p\rVert_{\textbf{M}}},
\end{align*}
where $p\in\mathbb{H}^n$.
	\end{proposition}
\begin{proof}
		Note that $\log_p(\exp_y(t\xi_y))/\lVert\log_p(\exp_y(t\xi_y))\rVert_{\textbf{M}}\rightarrow\xi_p$ as $t\rightarrow \infty$ by Proposition \ref{abc}(a), and since
		\begin{align*}
		\frac{\log_p(\exp_y(t\xi_y))}{\lVert\log_p(\exp_y(t\xi_y))\rVert_{\textbf{M}}}&=\frac{\cosh(t)y+\sinh(t)\xi_y+\langle p,\cosh(t)y+\sinh(t)\xi_y\rangle_{\textbf{M}} p}{\lVert\cosh(t)y+\sinh(t)\xi_y+\langle p,\cosh(t)y+\sinh(t)\xi_y\rangle_{\textbf{M}} p\rVert_{\textbf{M}}} \\
  &=\frac{y+\tanh(t)\xi_y+\langle p,y+\tanh(t)\xi_y\rangle_{\textbf{M}} p}{\lVert y+\tanh(t)\xi_y+\langle p,y+\tanh(t)\xi_y\rangle_{\textbf{M}} p\rVert_{\textbf{M}}},
		\end{align*}
		the expression for $\xi_p$ follows. The field defined by $\xi_p$ is clearly $C^\infty$ as a function of $p$ so long as the denominator does not equal 0. As remarked upon above, $\langle\cdot,\cdot\rangle_{\textbf{M}}$ is positive definite when restricted to $T_p\mathbb{H}^n$, so the denominator equals 0 if and only if the numerator does, which can only happen if $y+\xi_y$ is some scalar multiple of $p$ since $p\neq 0$. Suppose $y+\xi_y=ap$ for $a\in\mathbb{R}$. This is only possible if $a=0$ since $\langle y+\xi_y,y+\xi_y\rangle_{\textbf{M}}=0$ and $\langle ap,ap\rangle_{\textbf{M}}=-a^2$, so $y+\xi_y=0$. This is impossible because $\langle y,-y\rangle_{\textbf{M}}=1$ and $\langle y,\xi_y\rangle_{\textbf{M}}=0$; therefore the denominator cannot be 0.
	\end{proof}
		
		\begin{proof}[Proof of Theorem \ref{thm5.1}]
			Given an open interval $I\subset\mathbb{R}$ that contains 0 and a $C^\infty$ family of geodesics $\{\gamma_s\}_{s\in I}$, let $e_1,\ldots,e_{n-1}$ be parallel unit vector fields along $\gamma_0$ such that at each $t$, $\dot\gamma_0(t)/\lVert\dot\gamma_0(t)\rVert,e_1(t),\ldots,e_{n-1}(t)$ form an orthonormal basis for $T_{\gamma_0(t)}M$. Define $J_0(t):=-\dot\gamma_0(t)/\lVert\dot\gamma_0(t)\rVert$, $J_1(t):=-t\dot\gamma_0(t)/\lVert\dot\gamma_0(t)\rVert$,
			\begin{align*}
			J_{2r}(t)&:=\cosh(\sqrt{-\kappa}\lVert \dot\gamma_0(t)\rVert t)e_r(t),\\
			J_{2r+1}(t)&:=\sinh(\sqrt{-\kappa}\lVert \dot\gamma_0(t)\rVert t)e_r(t),
			\end{align*}
			for $r=1,\ldots,n-1$. The Jacobi equation on $M$ is equivalent to
			\begin{align*}
			\frac{D^2}{dt^2}J^\top(t)&=0, \\
			\frac{D^2}{dt^2}J^\perp(t)+\kappa\lVert\dot\gamma_0(0)\rVert^2 J^\perp(t)&=0,
			\end{align*}
			where $J^\top$ and $J^\perp$ are the components of $J$ that are tangential and orthogonal components, respectively, to $\dot\gamma_0$, and thus, $J_0$,\ldots,$J_{2n-1}$ are Jacobi fields along $\gamma_0$. Thus, since the space of Jacobi fields along a geodesic is of dimension $2n$, all Jacobi fields along $\gamma_0$ are of the form
			\begin{equation}\begin{aligned} \label{jacobi}
			&J(t)=-(A+B t)\frac{\dot\gamma_0(t)}{\lVert\dot\gamma_0(0)\rVert}+\cosh(\sqrt{-\kappa}\lVert\dot\gamma_0(0)\rVert t)w_1(t) \\
			&\qquad\qquad\qquad+\sinh(\sqrt{-\kappa}\lVert\dot\gamma_0(0)\rVert t)w_2(t),
			\end{aligned}\end{equation}
			where $A$ and $B$ are constants, and $w_1$ and $w_2$ are parallel vector fields along $
			\gamma_0$ such that $w_1(0)$ and $w_2(0)$ are orthogonal to $\dot\gamma_0(0)$.
			
			For any $v\in T_pM$, take an $\epsilon>0$ and a $C^\infty$ path $c:[-\epsilon,\epsilon]\rightarrow M$ such that $c(0)=p$, $\dot c(0)=v$ and $x\not\in c([-\epsilon,\epsilon])$, and any unit-speed geodesic ray $c^*:[0,\infty)\rightarrow M$ in the equivalence class $\xi$ such that the images of $c$ and $c^*$ are disjoint and $x\not\in c^*([-\epsilon,\epsilon])$. Define $y_m:=c^*(m)$ for each non-negative integer $m$. Then, $\zeta_m:[-1,1]\rightarrow\mathbb{R}$, defined by
			\begin{align*}
			\zeta_m(s)=\bigg\langle \frac{\log_{c(s)}(y_m)}{d(c(s),y_m)},\log_{c(s)}(x)\bigg\rangle
			\end{align*}
			is continuously differentiable by Proposition \ref{abc}(b), and
			\begin{equation}\begin{aligned} \label{zeroth}
			\zeta_m'(s_0)&=\bigg\langle\frac{\log_{c(s)}(y_m)}{d(c(s),y_m)},\frac{D}{ds}\log_{c(s)}(x)\bigg\rangle\bigg|_{s=s_0} \\
			&\qquad\qquad+\bigg\langle\frac{D}{ds}\frac{\log_{c(s)}(y_m)}{d(c(s),y_m)},\log_{c(s)}(x)\bigg\rangle\bigg|_{s=s_0}
			\end{aligned}\end{equation}
			for any $s_0\in[-\epsilon,\epsilon]$.
			
			Define a $C^\infty$ family of geodesics $\{\gamma_s^x\}_{s\in [-\epsilon,\epsilon]}$ by $\gamma_s^x(t)=\exp_x(t\cdot\log_x(c(s)))$. Note, for each $s_0\in [-\epsilon,\epsilon]$, that $J_{s_0}^x$ defined by $J_{s_0}^x(t)=(\partial/\partial s)\gamma_s^x(t)|_{s=s_0}$ is a Jacobi field along $\gamma_{s_0}^x$, and $\log_{c(s_0)}(x)=-(\partial/\partial t)\gamma_{s_0}^x(t)|_{t=1}$. Then, by the symmetry of the covariant derivative,
			\begin{equation}\begin{aligned}\label{switch1}
			\frac{D}{ds}\log_{c(s)}(x)\bigg|_{s=s_0}&=-\frac{D}{\partial s}\bigg(\frac{\partial}{\partial t}\gamma_s^x(t)\bigg|_{t=1}\bigg)\bigg|_{s=s_0} \\
   &=-\frac{D}{\partial t}\bigg(\frac{\partial}{\partial s}\gamma_s^x(t)\bigg|_{s=s_0}\bigg)\bigg|_{t=1} \\
   &=-\frac{D}{dt}J_{s_0}^x(t)\bigg|_{t=1}.
			\end{aligned}\end{equation}
			
			Applying (\ref{jacobi}) to the aforementioned Jacobi fields results in
			\begin{align*}
			J_{s_0}^x(t)&=-(A_{s_0}^x+B_{s_0}^x t)\frac{\dot\gamma_{s_0}^x(t)}{d(c(s_0),x)}+\cosh(\sqrt{-\kappa} d(c(s_0),x)t)w_{s_0,1}^x(t) \\
			&\qquad\qquad+\sinh(\sqrt{-\kappa} d(c(s_0),x)t)w_{s_0,2}^x(t).
			\end{align*}
			For each $s_0\in [-\epsilon,\epsilon]$, because $J_{s_0}^x(0)=(\partial/\partial s)\gamma_s^x(0)|_{s=s_0}=(d/ds)x=0$, we know that $A_{s_0}^x$ and $w_{s_0,1}^x(0)$, and hence, $w_{s_0,1}^x(t)$ for all $t$, are 0. On the other hand, $J_{s_0}^x(1)=(\partial/\partial s)\gamma_s^x(1)|_{s=s_0}=(d/ds)\exp_x(\log_x(c(s)))|_{s=s_0}=(d/ds)c(s)|_{s=s_0}=\dot c(s_0)$, so $B_{s_0}^x$ is $\langle \dot c(s_0),-\dot\gamma_{s_0}^x(t)/d(c(s_0),x)\rangle=\langle \dot c(s_0),\log_{c(s_0)}(x)/d(c(s_0),x)\rangle$ and 
			\begin{align*}
			&\sinh(\sqrt{-\kappa}d(c(s_0),x))w_{s_0,2}^x(1) \\
			&=\dot c(s_0)-\langle \dot c(s_0),\log_{c(s_0)}(x)/d(c(s_0),x)\rangle\log_{c(s_0)}(x)/d(c(s_0),x), 
			\end{align*}
			the projection of $\dot c(s_0)$ onto the orthogonal complement of $\log_{c(s_0)}(x)/d(c(s_0),x)$.
			
			In addition, because the covariant derivative of a parallel vector field along a $C^\infty$ curve is 0, 
			\begin{equation*}
			\frac{D}{dt}J_{s_0}^x(t)=-B_{s_0}^x\frac{\dot\gamma_{s_0}^x(t)}{d(c(s_0),x)}+\sqrt{-\kappa}d(c(s_0),x)\cosh(\sqrt{-\kappa}d(c(s_0),x) t)w_{s_0,2}^x(t).
			\end{equation*}
			The above result, the results in the preceding paragraph, and (\ref{switch1}) imply
			\begin{equation}\begin{aligned}\label{first}
			&\bigg\langle\frac{\log_{c(s)}(y_m)}{d(c(s),y_m)},\frac{D}{ds}\log_{c(s)}(x)\bigg\rangle\bigg|_{s=s_0} \\
			&=\Bigg\langle\frac{\log_{c(s_0)}(y_m)}{d(c(s_0),y_m)},-\Bigg(\bigg\langle \dot c(s_0),\frac{\log_{c(s_0)}(x)}{d(c(s_0),x)}\bigg\rangle\frac{\log_{c(s_0)}(x)}{d(c(s_0),x)}  \\
			&\qquad+\sqrt{-\kappa}d(c(s_0),x)\coth(\sqrt{-\kappa}d(c(s_0),x))\bigg(\dot c(s_0)-\bigg\langle \dot c(s_0),\frac{\log_{c(s_0)}(x)}{d(c(s_0),x)}\bigg\rangle \\
   &\qquad\times\frac{\log_{c(s_0)}(x)}{d(c(s_0),x)}\bigg)\Bigg)\Bigg\rangle \\
			&=-\Bigg\langle\frac{\log_{c(s_0)}(y_m)}{d(c(s_0),y_m)},(1-\sqrt{-\kappa}d(c(s_0),x)\coth(\sqrt{-\kappa}d(c(s_0),x))) \\
			&\qquad\times\bigg\langle \dot c(s_0),\frac{\log_{c(s_0)}(x)}{d(c(s_0),x)}\bigg\rangle\frac{\log_{c(s_0)}(x)}{d(c(s_0),x)} \\
			&\qquad+\sqrt{-\kappa}d(c(s_0),x)\coth(\sqrt{-\kappa}d(c(s_0),x))\dot c(s_0)\Bigg\rangle.
			\end{aligned}\end{equation}
			
			We can similarly conclude about another $C^\infty$ family of geodesics $\{\gamma_s^m\}_{s\in [-\epsilon,\epsilon]}$ defined by $\gamma_s^m(t)=\exp_y(t\cdot\log_y(c(s))/d(c(s),y_m))$ as follows: for each $s_0\in [-\epsilon,\epsilon]$, $J_{s_0}^m$ defined by $J_{s_0}^m(t)=(\partial/\partial s)\gamma_{s_0}^m(t)$, is a Jacobi field along $\gamma_{s_0}^m$, 
			\[
			\frac{\log_{c(s_0)}(y_m)}{d(c(s_0),y_m)}=-\frac{\partial}{\partial t}\gamma_{s_0}^m(t)\bigg|_{t=d(c(s_0),y_m)},
			\]
			\[
			\frac{D}{ds}\frac{\log_{c(s)}(y_m)}{d(c(s),y_m)}\bigg|_{s=s_0}=-\frac{D}{dt}J_{s_0}^m(t)\bigg|_{t=d(c(s_0),y_m)},
			\]
			\[
			J_{s_0}^m(t)=-(A_{s_0}^m+B_{s_0}^m t)\dot\gamma_{s_0}^m(t)+\cosh(\sqrt{-\kappa} t)w_{s_0,1}^m(t)+\sinh(\sqrt{-\kappa}t)w_{s_0,2}^m(t),
			\]
			and $A_{s_0}^m=0$, $w_{s_0,1}^m=0$, $B_{s_0}^m=(1/d(c(s_0),y_m))\langle v,\log_{c(s_0)}(y_m)/d(c(s_0),y_m)\rangle$, 
   \begin{align*}
   &\sinh(\sqrt{-\kappa}d(c(s_0),y_m))w_{s_0,2}^m(d(c(s_0),y_m)) \\
   &=\dot c(s_0)-\bigg\langle \dot c(s_0),\frac{\log_{c(s_0)}(y_m)}{d(c(s_0),y_m)}\bigg\rangle\frac{\log_{c(s_0)}(y_m)}{d(c(s_0),y_m)},
  \end{align*}
			\begin{equation*}
			\frac{D}{dt}J_{s_0}^m(t)=-B_{s_0}^m\dot\gamma_{s_0}^m(t)+\sqrt{-\kappa}\cosh(\sqrt{-\kappa} t)w_{s_0,2}^m(t);
			\end{equation*}
			consequently
			\begin{equation}\begin{aligned}\label{second}
			&\bigg\langle\frac{D}{ds}\frac{\log_{c(s)}(y_m)}{d(c(s),y_m)},\log_{c(s)}(x)\bigg\rangle\bigg|_{s=s_0} \\
			&=-\Bigg\langle\bigg(\frac{1}{d(c(s_0),y_m)}-\sqrt{-\kappa}\coth(\sqrt{-\kappa}d(c(s_0),y_m))\bigg) \\
			&\qquad\times\bigg\langle \dot c(s_0),\frac{\log_{c(s_0)}(y_m)}{d(c(s_0),y_m)}\bigg\rangle\frac{\log_{c(s_0)}(y_m)}{d(c(s_0),y_m)} \\
			&\qquad+\sqrt{-\kappa}\coth(\sqrt{-\kappa}d({c(s_0)},y_m))\dot c(s_0),\log_{c(s)}(x)\Bigg\rangle.
			\end{aligned}\end{equation}
			Substituting (\ref{first}) and (\ref{second}) into (\ref{zeroth}) and rearranging, we obtain
			\begin{equation}\begin{aligned} \label{zetaprime}
			&\zeta_m'(s) \\
			&=-(1-\sqrt{-\kappa}d(c(s),x)\coth(\sqrt{-\kappa}d(c(s),x))) \\
			&\qquad\bigg\langle\frac{\log_{c(s)}(y_m)}{d(c(s),y_m)},\frac{\log_{c(s)}(x)}{d(c(s),x)}\bigg\rangle\bigg\langle\frac{\log_{c(s)}(x)}{d(c(s),x)},\dot c(s)\bigg\rangle \\
   &\qquad-\sqrt{-\kappa}d(c(s),x)\coth(\sqrt{-\kappa}d(c(s),x))\bigg\langle\frac{\log_{c(s)}(y_m)}{d(c(s),y_m)},\dot c(s)\bigg\rangle \\
			&\qquad -\bigg(\frac{1}{d(c(s),y_m)}-\sqrt{-\kappa}\coth(\sqrt{-\kappa}d(c(s),y_m))\bigg) \\
			&\qquad\times\bigg\langle\frac{\log_{c(s)}(y_m)}{d(c(s),y_m)},\log_{c(s)}(x)\bigg\rangle\times\bigg\langle\frac{\log_{c(s)}(y_m)}{d(c(s),y_m)},\dot c(s)\bigg\rangle \\
			&\qquad-\sqrt{-\kappa}\coth(\sqrt{-\kappa}d(c(s),y_m))\bigg\langle\log_{c(s)}(x),\dot c(s)\bigg\rangle.
			\end{aligned}\end{equation}
			By the Cauchy-Schwarz inequality and the proof of Proposition \ref{abc}(a), 
			\begin{align*}
			\bigg\lvert\bigg\langle\frac{\log_{c(s)}(y_m)}{d(c(s),y_m)},\dot c(s)\bigg\rangle-\bigg\langle\xi_{c(s)},\dot c(s)\bigg\rangle\bigg\rvert&=\bigg\lvert\bigg\langle\frac{\log_{c(s)}(y_m)}{d(c(s),y_m)}-\xi_{c(s)},\dot c(s)\bigg\rangle\bigg\rvert \\
			&\leq\bigg\lVert\frac{\log_{c(s)}(y_m)}{d(c(s),y_m)}-\xi_{c(s)}\bigg\rVert\lVert\dot c(s)\rVert \\
			&\leq \frac{2d(c(s),y_0)}{m}\lVert\dot c(s)\rVert;
			\end{align*}
			by the continuity, and hence the boundedness, on $[-\epsilon,\epsilon]$ of $d(c(s),y_0)$ and $\lVert\dot c(s)\rVert$ as functions of $s$, $\langle\log_{c(s)}(y_m)/d(c(s),y_m),\dot c(s)\rangle$ converges uniformly as a function of $s$ on $[-\epsilon,\epsilon]$ to $\langle\xi_{c(s)},\dot c(s)\rangle$ as $m\rightarrow \infty$. Similarly, $\langle\log_{c(s)}(y_m)/d(c(s),y_m),\log_{c(s)}(x)/d(c(s),x)\rangle$ and $\langle\log_{c(s)}(y_m)/d(c(s),y_m),\log_{c(s)}(x)\rangle$ converge uniformly as functions of $s\in[-\epsilon,\epsilon]$ to $\langle\xi_{c(s)},\log_{c(s)}(x)/d(c(s),x)\rangle$ and $\langle\xi_{c(s)},\log_{c(s)}(x)\rangle$, respectively, as $m\rightarrow \infty$. Also, $d(c(s),y_m)\geq d(y_1,y_m)-d(c(s),y_0)=m-d(c(s),y_0)$, and thus, $1/d(c(s),y_m)$ and $\coth(\sqrt{-\kappa}d(c(s),y_m))$ converge uniformly to $0$ and $1$, respectively, on $[-\epsilon,\epsilon]$. Finally, the product of uniformly convergent sequences of bounded functions is uniformly convergent to the product of the limits, and each factor of each summand in (\ref{zetaprime}) is continuous as a function of $s$ on $[-\epsilon,\epsilon]$, and hence, bounded. Therefore, the results in this paragraph imply that $\zeta_m':[-\epsilon,\epsilon]\rightarrow\mathbb{R}$ converges uniformly as $m\rightarrow\infty$ to the function of $s$ on $[-\epsilon,\epsilon]$ defined by
			\begin{align*}
			&-(1-\sqrt{-\kappa}d(c(s),x)\coth(\sqrt{-\kappa}d(c(s),x)))\bigg\langle\xi_{c(s)},\frac{\log_{c(s)}(x)}{d(c(s),x)}\bigg\rangle\bigg\langle\frac{\log_{c(s)}(x)}{d(c(s),x)},\dot c(s)\bigg\rangle  \\
			&\qquad-\sqrt{-\kappa}d(c(s),x)\coth(\sqrt{-\kappa}d(c(s),x))\bigg\langle\xi_{c(s)},\dot c(s)\bigg\rangle \\
			&\qquad+\sqrt{-\kappa}\bigg\langle\xi_{c(s)},\log_{c(s)}(x)\bigg\rangle\bigg\langle\xi_{c(s)},\dot c(s)\bigg\rangle-\sqrt{-\kappa}\bigg\langle\log_{c(s)}(x),\dot c(s)\bigg\rangle \\
			&=\Bigg\langle-\Bigg(\sqrt{-\kappa}d(c(s),x)\bigg(\coth(\sqrt{-\kappa}d(c(s),x))-\bigg\langle\xi_{c(s)},\frac{\log_{c(s)}(x)}{d(c(s),x)}\bigg\rangle\bigg)\xi_{c(s)}  \\
			&\qquad+\bigg((1-\sqrt{-\kappa}d(c(s),x)\coth(\sqrt{-\kappa}d(c(s),x)))\bigg\langle\xi_{c(s)},\frac{\log_{c(s)}(x)}{d(c(s),x)}\bigg\rangle \\
   &\qquad+\sqrt{-\kappa}d(c(s),x)\bigg)\frac{\log_{c(s)}(x)}{d(c(s),x)}\Bigg),\dot c(s)\Bigg\rangle.
			\end{align*}		
			Then, because $\zeta_m$ converges pointwise to $\zeta:[-\epsilon,\epsilon]\rightarrow\mathbb{R}$, defined by $\zeta(s)=\langle\xi_{c(s)},\log_{c(s)}(x)\rangle$, by Proposition \ref{abc}(a), $\zeta':[-\epsilon,\epsilon]\rightarrow\mathbb{R}$ defined by the above expression is the derivative of $\zeta$ by the differentiable limit theorem. Therefore, the gradient of $\langle\xi_p,\log_p(x)\rangle$ with respect to $p=c(0)$ is
			\begin{align*}
			&-\sqrt{-\kappa}d(p,x)\bigg(\coth(\sqrt{-\kappa}d(p,x))-\bigg\langle\xi_{p},\frac{\log_{p}(x)}{d(p,x)}\bigg\rangle\bigg)\xi_{p}  \\
			&\quad-\bigg((1-\sqrt{-\kappa}d(p,x)\coth(\sqrt{-\kappa}d(p,x)))\bigg\langle\xi_{p},\frac{\log_{p}(x)}{d(p,x)}\bigg\rangle+\sqrt{-\kappa}d(p,x)\bigg)\frac{\log_{p}(x)}{d(p,x)}.
			\end{align*}
			Recalling that this whole process required $x\not\in c([-\epsilon,\epsilon])$, the above expression does not apply when $p=x$. However, $\lim_{t\rightarrow 0}t\coth(t)=1$ means that $\lim_{p\rightarrow x}\sqrt{-\kappa}d(p,x)\coth(\sqrt{-\kappa}d(p,x))=1$, which together with $\lVert\log_p(x)/d(p,x)\rVert\leq 1$ implies that the coefficients of $\xi_p$ and $\log_p(x)/d(p,x)$ in the above expression approach $-1$ and $0$, respectively, as $p\rightarrow x$. Because Proposition \ref{abc}(b) implies the continuity of $p\mapsto\nabla\langle\xi_p,\log_p(x)\rangle$, this means that $\nabla \langle\xi_p,\log_p(x)\rangle$ at $p=x$ is $-\xi_p$.
		\end{proof}
	\end{appendix}

 \bibliographystyle{apalike}
\bibliography{references}

\end{document}